\documentclass[times]{amsart}

\usepackage{amssymb,amsmath,amsthm,stmaryrd}
\usepackage[colorlinks,plainpages,backref]{hyperref}
\usepackage{mathrsfs}
\usepackage[neveradjust]{paralist}
\usepackage[all,cmtip]{xy}
\xyoption{dvips}
\usepackage{epic,eepic}
\usepackage{url,epsfig,psfrag,graphicx,subfigure}
\usepackage{float,wrapfig}
\usepackage{tikz}
\usepackage{tikz-cd}
\usepackage{adjustbox}

\restylefloat{figure}

\theoremstyle{definition}
\newtheorem{ntn}{Notation}[section]

\newtheorem{dfn}[ntn]{Definition}

\newtheorem{cnv}[ntn]{Convention}
\newtheorem*{cnvv}{Notational conventions}

\theoremstyle{plain}
\newtheorem{lem}[ntn]{Lemma}
\newtheorem{prp}[ntn]{Proposition}
\newtheorem{thm}[ntn]{Theorem}
\newtheorem{cor}[ntn]{Corollary}
\newtheorem{cnj}[ntn]{Conjecture}

\theoremstyle{remark}

\newtheorem{rmk}[ntn]{Remark}
\newtheorem{exa}[ntn]{Example}

\newcommand{\boldzero}{{\mathbf{0}}}
\newcommand{\boldone}{{\mathbf{1}}}
\newcommand{\bolda}{{\mathbf{a}}}
\newcommand{\boldb}{{\mathbf{b}}}
\newcommand{\boldc}{{\mathbf{c}}}

\newcommand{\bolds}{{\mathbf{s}}}
\newcommand{\boldu}{{\mathbf{u}}}
\newcommand{\boldv}{{\mathbf{v}}}

\newcommand{\mathboldq}{{\boldsymbol{q}}}
\newcommand{\bolddel}{{\boldsymbol \partial}}
\newcommand{\boldeps}{{\boldsymbol \varepsilon}}
\newcommand{\boldxi}{{\boldsymbol \xi}}
\newcommand{\boldt}{{\boldsymbol{t}}}
\newcommand{\boldx}{{\boldsymbol{x}}}
\newcommand{\boldy}{{\boldsymbol{y}}}

\newcommand{\del}{\partial}

\newcommand{\de}{{\mathrm d}}

\newcommand{\eps}{\varepsilon}

\newcommand{\ideal}[1]{{\langle#1\rangle}}
\newcommand{\into}{\hookrightarrow}
\newcommand{\iso}{\simeq}

\renewcommand{\to}{\longrightarrow}

\newcommand{\minus}{\smallsetminus}

\newcommand{\calD}{\mathscr{D}}
\newcommand{\calE}{{\mathscr{E}}}

\newcommand{\calG}{\mathscr{G}}
\newcommand{\calI}{\mathscr{I}}
\newcommand{\calH}{\mathscr{H}}
\newcommand{\calJ}{\mathscr{J}}
\newcommand{\calL}{\mathscr{L}}
\newcommand{\calM}{\mathscr{M}}
\newcommand{\calN}{\mathscr{N}}
\newcommand{\calO}{\mathscr{O}}

\newcommand{\calS}{\mathscr{S}}

\newcommand{\calK}{\mathcal{K}}

\newcommand{\frakd}{{\mathfrak{d}}}
\newcommand{\frakm}{{\mathfrak{m}}}
\newcommand{\frakp}{{\mathfrak{p}}}

\renewcommand{\AA}{\mathbb{A}}
\newcommand{\CC}{\mathbb{C}}

\newcommand{\NN}{\mathbb{N}}
\newcommand{\PP}{\mathbb{P}}
\newcommand{\QQ}{\mathbb{Q}}
\newcommand{\RR}{\mathbb{R}}
\newcommand{\TT}{\mathbb{T}}
\newcommand{\VV}{\mathbb{V}}
\newcommand{\ZZ}{\mathbb{Z}}

\newcommand{\pointl}{{\mathfrak l}}
\newcommand{\pointq}{{\mathfrak q}}
\newcommand{\pointt}{{\mathfrak t}}
\newcommand{\pointw}{{\mathfrak w}}
\newcommand{\pointx}{{\mathfrak x}}
\newcommand{\pointy}{{\mathfrak y}}
\newcommand{\pointz}{{\mathfrak z}}

\renewcommand{\bar}{\overline}

\DeclareMathOperator{\ch}{ChV}
\DeclareMathOperator{\cc}{ChC}

\DeclareMathOperator{\conv}{conv}

\DeclareMathOperator{\Div}{\textup{Div}}

\DeclareMathOperator{\erf}{erf}

\DeclareMathOperator{\FL}{FL}
\DeclareMathOperator{\FS}{FS}
\DeclareMathOperator{\gr}{gr}
\DeclareMathOperator{\Gr}{Gr}

\DeclareMathOperator{\HM}{HM}
\DeclareMathOperator{\Hom}{Hom}
\DeclareMathOperator{\id}{id}
\DeclareMathOperator{\IC}{IC}
\DeclareMathOperator{\IH}{IH}
\DeclareMathOperator{\image}{im}

\DeclareMathOperator{\MHM}{MHM}

\DeclareMathOperator{\PD}{PD}

\DeclareMathOperator{\Proj}{Proj}
\DeclareMathOperator{\qdeg}{qdeg}
\DeclareMathOperator{\QDM}{QDM}
\DeclareMathOperator{\Radon}{RT}
\DeclareMathOperator{\Radonshort}{R}
\DeclareMathOperator{\rk}{rk}

\DeclareMathOperator{\Sing}{Sing}

\DeclareMathOperator{\Sol}{Sol}
\DeclareMathOperator{\sRes}{sRes}
\DeclareMathOperator{\Spec}{Spec}

\DeclareMathOperator{\supp}{supp}
\DeclareMathOperator{\tdeg}{tdeg}

\DeclareMathOperator{\Var}{Var}
\DeclareMathOperator{\vol}{vol}
\DeclareMathOperator{\pr}{pr}

\newcommand{\an}{{\mathrm{an}}}
\newcommand{\can}{{\mathrm{can}}}
\newcommand{\cst}{{\mathrm{cst}}}
\newcommand{\Hodge}{{\mathrm{Hodge}}}
\newcommand{\irr}{{\mathrm{irr}}}
\newcommand{\loc}{{\mathrm{loc}}}
\newcommand{\ord}{{\mathrm{ord}}}
\newcommand{\rh}{{\mathrm{rh}}}
\newcommand{\tw}{{\mathrm{tw}}}

\input Rothdn.fd
\input RoyalIn.fd

\newcommand{\sideA}{{\ensuremath{\boldsymbol{\mathcal A}}}}
\newcommand{\sideB}{{\ensuremath{\boldsymbol{\mathcal B}}}}
\newcommand{\Rees}{{R}}
\newcommand{\shRees}{{\mathscr R}}

\def\schluss{\hfill\ensuremath{\Diamond}}

\def\comment#1{}

\begin{document}

\title{Algebraic aspects of hypergeometric differential equations}

\author[T.~Reichelt]{Thomas Reichelt} \address{Universit\"at Heidelberg,
  Mathematisches Institut, Im Neuenheimer Feld 205, 69120 Heidelberg,
  Germany} \email{treichelt@mathi.uni-heidelberg.de}

\author[M.~Schulze]{Mathias Schulze}
\address{Technische Universit\"at Kaiserslautern, Fachbereich Mathematik,
Postfach 3049, 67653 Kaiserslautern, Germany}
\email{mschulze@mathematik.uni-kl.de}

\author[C.~Sevenheck]{Christian Sevenheck}
\address{Technische Universit\"at
  Chemnitz, Fakult\"at f\"ur Mathematik, 09107 Chemnitz, Germany}
\email{christian.sevenheck@mathematik.tu-chemnitz.de}

\author[U.~Walther]{Uli Walther}
\address{Department of Mathematics, Purdue University, 150 N
  University St., West Lafayette, IN~47907, USA}
\email{walther@math.purdue.edu}

\thanks{TR was supported by DFG Emmy-Noether-Fellowship RE 3567/1-2.}
\thanks{CS acknowledges partial support by the DFG grant SE1114/5-2.}
\thanks{UW was supported in part by Simons Collaboration Grant for
  Mathematicians \#580839}

\begin{abstract}
We review some classical and modern aspects of hypergeometric
differential equations, including $A$-hypergeometric systems of
Gel{$'$}fand, Graev, Kapranov and Zelevinsky. Some recent advances in
this theory, such as Euler-Koszul homology, rank jump phenomena,
irregularity questions and Hodge theoretic aspects are discussed with
more details. We also give some applications of the theory of
hypergeometric systems to toric mirror symmetry.
\end{abstract}

\maketitle

\tableofcontents

\section{Introduction}
\label{sec:Introduction}

\begin{cnvv}
  We use Italic letters $M$ for rings, variables and modules;
  calligraphic letters $\calD$ for sheaves; Roman letters $\FL$ for
  functors; Gothic letters for prime ideals $\frakp$ and points
  $\pointx$ of spaces.

  Lattice elements $\bolda$ are in Roman bold; coordinate sets $\boldt$ and
  other sets of functions or operators $\bolddel$ in Italic bold.
  \schluss\end{cnvv}

\subsection{Hypergeometric functions}

The study of hypergeometric functions started more than two centuries
ago and formed a important part of the work of Euler and Gau\ss.
A power series
\[
f(z)=\sum_{i=0}^\infty a_i z^i/i!
\]
is \emph{hypergeometric} if the quotient $a_{i+1}/a_i$ of consecutive
coefficients is a rational function in $i$.
Traditional convention dictates that the exponential function is
regarded as the standard hypergeometric function (to $a_{i+1}/a_i$
constant); this ``explains''
the choice of $a_i/i!$ over $a_i$ as series coefficient.
Further examples include Bessel, Airy, trigonometric and (higher)
logarithmic as well as all other special functions, and the
hypergeometric functions that express roots of algebraic equations \cite{Sturmfels-eqnsolving}.

The continuing interest in hypergeometric functions stems to some
extent from the fact that they are often solutions
to very appealing linear differential equations taken from
physics. For example, the Bessel functions $J_{\pm r}(x)$ of the first
kind
arise as solutions to a linear second order equation
that shows up in heat and
electromagnetic propagation in a cylinder, vibrations of circular
membranes, and more generally when solving the Helmholtz or Laplace
equation. Indeed, such connections to physics through differential
equations prompted the first
studies of (specific) hypergeometric functions. However, hypergeometric functions also appear in many other parts of mathematics: as we will see soon, each time an action of an algebraic torus on a space is observed, one can expect to find some differential equation of hypergeometric type connected to this situation. The abundance of  toric varieties in geometry explains why there are so many different interesting hypergeometric functions.
We discuss in Section \ref{sec:MirrorSym} below one prominent case where hypergeometric differential equations prove to be useful: the so-called mirror symmetry phenomenon for certain smooth toric varieties. Other recent applications that are beyond the scope of this article include the holonomic gradient method in algebraic statistics (\cite{HibiTakayama}) or Feynman integral computations in quantum field theory (\cite{QFT0},\cite{QFT1},\cite{QFT2},\cite{QFT3}).

As it turns out, it is exactly the type of differential equation
satisfied by a function that determines whether the function should be
considered as hypergeometric, since these force the right kind of
recursions on the series.  The most successful approach to generalize
hypergeometric differential equations to several variables was
initiated by Gel{$'$}fand, Graev, Kapranov and Zelevinsky in the 1980s, and
some of the features of this theory form the topic of this
article. We start with some motivating examples.

\begin{exa}[The error function, part I]\label{exa-error}
The (Gau\ss) \emph{error function} $\erf(x)$ is defined by
\[
\erf(z)=\frac{2}{\sqrt\pi}\int_0^z\exp(-t^2)\,\de t.
\]
While this integral cannot be solved in closed form, it can be
developed into a convergent Taylor series
\begin{eqnarray}\label{eqn-erf-hyp-eqn}
\erf(z)=\frac{2}{\sqrt{\pi}}z\sum_{i=0}^\infty a_i\frac{(-z^2)^i}{i!}
\end{eqnarray}
where $a_i=1/(2i+1)$,
so that
\[
\erf(z)=\frac{2z}{\sqrt\pi}
 \left(1-\frac{z^2}{3}+\frac{(z^2)^2}{10}-\frac{(z^2)^3}{42}-
 \frac{(z^2)^4}{216}+ \frac{(z^2)^5}{1320}\mp\cdots\right)
\]
is hypergeometric.
\schluss\end{exa}

The univariate hypergeometric functions are classified by the rational
function $a_{i+1}/a_i$. More precisely, suppose that
$a_{i+1}/a_i=P(i)/Q(i)$ where $P,Q\in\CC[i]$ are monic with
$P=\prod_{j=1}^p(i+\alpha_j)$ and $Q=\prod_{j=1}^q(i+\beta_j)$. Then
the \emph{univariate hypergeometric function associated to $P,Q$} is
\begin{eqnarray}\label{eqn-pFq}
_pF_q(\alpha_1,\ldots,\alpha_p;\beta_1,\ldots,\beta_q;z)&=&\sum_{i=0}^\infty
\frac{a_i z^i}{i!}
\end{eqnarray}
where $a_0=1$ and
\[
\frac{a_{i+1}}{a_i}=\frac{(i+\alpha_1)(i+\alpha_2)\cdots(i+\alpha_p)}{(i+\beta_1)(i+\beta_2)\cdots(i+\beta_q)}.
\]
\begin{exa}[The error function, part II]
  It follows from \eqref{eqn-erf-hyp-eqn} that
  $\erf(z)$ is, up to the factor $2z/\sqrt{\pi}$, equal to
  ${}_1F_1(1/2;3/2;-z^2)$, where
      \[
    {}_1F_1(1/2;3/2;z)=1+\frac{z}{3}+\frac{z^2}{10}+\frac{z^3}{42}+
    \frac{z^4}{216}+ \frac{z^5}{1320}+\cdots
    \]
    is the \emph{Kummer confluent function} which encodes all intrinsic
    analytic and combinatorial properties of $\erf(x)$ and, with
    $\theta_z=z\frac{d}{dz}$, satisfies the differential equation
    \begin{eqnarray}\label{eqn-Kummer}
    \theta_z(\theta_z+1/2)\bullet(f)-z(\theta_z-1/2)\bullet(f)=0.
    \end{eqnarray}
    The particular shape of this equation will be used in the next section
    for a conversion process from univariate hypergeometric functions to
    $A$-hypergeometric ones.
\schluss\end{exa}

In the following example we document how hypergeometric functions
arise naturally from differential forms with parameters.  The
computation was apparently already known to Kummer; compare
\cite{BrieskornKnoerrer} for details. In modern terms, it represents
the birth of the notion of a variation of Hodge structures.

\begin{exa}[Hypergeometry and Hodge filtrations]\label{exa-elliptic}

The equation $f_z=0$ with
\[
f_z(u,v)=v^2-u(u-1)(u-z)
\]
defines for each $z\in\CC\minus\{0,1\}$ a smooth curve $E_z$ over
$\CC$. Its projective closure $\bar E_z\subseteq \PP^2_\CC$ meets the
line at infinity in a single point and is smooth as long as
$z\not\in\{ 0,1,\infty\}$. The natural projection from $E_z$ to $\CC$
via ``forgetting $v$'' is generically $2:1$ and branches at $0,1,z$;
the induced map $\bar E_z\to\PP^1_\CC$ also branches at infinity.

The differential $1$-form $\omega_z:=\de u/v$ is everywhere
holomorphic and nowhere zero on $\bar E_z$; the existence of this
``form of the first kind'' in Riemann's language makes the elliptic
curve $\bar E_z$ a Calabi--Yau manifold in modern terms. The ``form of
the second kind'' $\omega'_z:=\omega_z/(u-z)$ has a unique pole, at
$u=z$, at which it is residue-free. Considering $v=v(u,z)$ as dependent
variable and writing $\omega_z, \omega'_z$ in terms of $u$ and $z$,
one notes that $\frac{\del}{\del
  z}(\omega_z)=\frac{1}{2}\omega'_z$, and (compare especially
\cite[Page 685]{BrieskornKnoerrer})
\[
\frac{\del}{\del z}(\omega'_z)=\frac{3\de u}{4v(u-z)^2}=\underbrace{\frac{1}{4z(1-z)}}_{p(z)}\omega_z
+\underbrace{\frac{-1+2z}{z(1-z)}}_{q(z)}\omega'_z
+\de\left(\frac{v}{2(u-z)^2z(1-z)}\right),
\]
the differential on the right being taken in $u,v$ with $z$ constant
(and noting that on $E$ one has $\de(u(u-1)(u-z))=2v\,\de v$).

Let $\lambda\in H_1(\bar E_z;\ZZ)\iso \ZZ\oplus\ZZ$ and set
$I_1(\lambda)=\int_\lambda\omega_z$ and
$I_2(\lambda)=\int_\lambda\omega'_z$, multi-valued functions on $\bar
E_z$ defined via elliptic integrals.  The differential equations for
$\omega_z,\omega'_z$ imply (compare \cite[Lemma
  12]{BrieskornKnoerrer}) that $I_1(\lambda)$ and $I_2(\lambda)$ are
solutions to
\begin{equation}\label{eqn-forms}
f''-qf'=pf,
\end{equation}
with singularities at $0,1$ and $\infty$.  It is the special case
$1=2a=2b=c$ of the general Gau\ss\ hypergeometric differential
equation
\[
f''+\frac{c-(a+b+1)t}{z(1-z)}f'=\frac{ab}{z(1-z)}f
\]
with solution space basis given by Gau\ss' hypergeometric functions
\begin{eqnarray*}
  F_1&=&\sum_{n=0}^\infty \frac{[a]_n[b]_n}{[c]_n}\frac{z^n}{n!},\\
F_2&=&-\sqrt{-1}\sum_{n=0}^\infty
\frac{[a]_n[b]_n}{[c]_n}\frac{(1-z)^n}{n!},
\end{eqnarray*}
which have singularities at $0,\infty$ and $1,\infty$ respectively.

Suppose $\lambda_z,\lambda'_z$ are the standard basis (the minimal
geodesics) for the first homology group of the torus $\bar E_z$.  Then
two elementary (but non-trivial) computations reveal:
\begin{asparaenum}
\item analytic continuation of the solution space basis $F=(F_1,F_2)^T$
  around the points $z=0$ and $z=1$ corresponds to multiplication of
  $F$ by $M_0=\begin{pmatrix}1&0\\-2&1\end{pmatrix}$ and
  $M_1=\begin{pmatrix}1&2\\0&1\end{pmatrix}$ respectively;
\item the map
  \begin{eqnarray*}
    \pi\colon \PP^2_\CC\minus \{(1,0,0),(1,1,0),(0,0,1)\}&\to&
    \PP^1_\CC,\\
    wu(u-w)&\mapsfrom&z_0,\\
    wv^2-u^3-u^2w&\mapsfrom&z_1,
  \end{eqnarray*}
  is a bundle with fiber $\bar E_{z_1/z_0}$ that admits an
  Ehresmann connection. In particular, the cohomology classes of
  the fibers allow parallel transport. The induced vector
  bundle with fiber $H_1(\bar E_z;\ZZ)=\ZZ \lambda_z+\ZZ\lambda'_z$
  admits a monodromy action, lifting the loops around $z=(0,1)$
  and $z=(1,1)$. Analysis of the geometry of $\pi$ shows that this
  monodromy is given  again by the actions of $M_1$ and
  $M_2$ respectively.
\end{asparaenum}
More abstractly, the $D$-module on the base of $\pi$ corresponding to
the (derived) direct image (compare Notation \ref{ntn-f}) of the
structure sheaf on the source of $\pi$, also known as the
\emph{Gau\ss\--Manin system}, has monodromy action via $M_1,M_2$.

On the complement of the points $0, 1, \infty$ this $D_z$-module is a
vector bundle with a flat connection. The fibers of this vector bundle
are the cohomology groups $H^1(\overline{E}_{z_1 /z_0};\CC)$. This
vector bundle is actually a variation of pure Hodge structures of
weight $1$ where the $(1,0)$-part is generated by the differential
form $\omega_z$, the variation of this $(1,0)$-subbundle being described
by \eqref{eqn-forms}.

It follows that, up to scalars, $I_1(\lambda_z)=F_1(z)$,
$I_2(\lambda_z)=F_2(z)$. In particular, the
ratio $\tau(z)=I_1(\lambda_z)/I_2(\lambda_z)$ is the modulus of the
elliptic curve in the sense that the fiber over $z$ is isomorphic to
the quotient of $\CC$ by $\ZZ+\sqrt{-1}\tau \cdot \ZZ$.

We will take up the discussion of Hodge structures associated to more general univariate hypergeometric operators (see equation \eqref{eqn-gauss-diffeq} below) later in Section \ref{sec:Hodge} (see page \pageref{page:SimpsonFedorov}).
\schluss\end{exa}

\subsection{From univariate to GKZ and back}\label{subsec-dim-red}

In the 1980s, the Russian school around I.M.~Gel{$'$}fand found a universal
way of encoding univariate hypergeometric functions by way of certain
systems of PDEs that arise from an integer matrix $A$ and complex
parameter vector $\beta$. We start with the general definition and
then explain how univariate hypergeometric functions arise as
solutions of these $D$-modules.

\begin{ntn}
  In the first three sections of this article,
  \[
  A=(\bolda_1,\ldots,\bolda_n)\in\ZZ^{d\times n}
  \]
  denotes an integer matrix with $d$ rows and $n$ columns. In the last
  two sections, $A$ will still be integer, but at least sometimes of
  size $(d+1)\times (n+1)$.
\schluss\end{ntn}

For convenience, we place the following constraints on the matrix $A$;
they make concise statements possible, or at least easier to make.
\begin{cnv}[Standard assumptions on $A$]\label{cnv-basicA}
  With $A$ as above, $A$ spans a semigroup
  \[
  \NN A:=\sum_{j=1}^n \NN\bolda_j \subseteq \ZZ A
  \]
  inside $\ZZ^d$.
  Throughout we assume that
  \begin{itemize}
  \item the group $\ZZ A$ generated by $A$ agrees with $\ZZ^d$ ($A$ is
    \emph{full});
  \item the semigroup $\NN A$ contains no units besides $\boldzero$
    ($A$ is \emph{pointed}). We note that pointedness of $A$ is
    equivalent to the existence of a group homomorphism from $\ZZ^d$
    to $\ZZ$ that is positive on every $\bolda_j$.
  \end{itemize}
\schluss\end{cnv}

We now give the definition of the main character of our story.
\begin{dfn}[$A$-hypergeometric system, \cite{GGZ87}]\label{def:GKZ}
  Fix $A\in \ZZ^{d\times n}$ as in Convention \ref{cnv-basicA} and
  choose $\beta\in\CC^d$. Let
  \[
  D_A:=\CC[\boldx]\ideal{\bolddel}
  \]
  be the $n$-th Weyl algebra over $\CC$.
  Here $\boldx=x_1,\ldots,x_n, \bolddel=\del_1,\ldots,\del_n$,
  and $\del_j$ is identified with the partial differentiation
  operator $\frac{\del}{\del x_j}$. We also let
  \[
  R_A:=\CC[\bolddel]\subseteq D_A
  \]
  denote the polynomial subring.

  Letting $\theta_j$ stand for
  $x_j\del_j$, the
  \emph{Euler operator} $E_i$ is
  \[
  E_i=\sum_{j=1}^n a_{i,j}\theta_j.
  \]
  For each $\boldu\in\ZZ^n$ in the kernel of $A$ its \emph{box
    operator} is
  \[
  \Box_\boldu=\bolddel^{\boldu_+}-\bolddel^{\boldu_-},
  \]
  where $(\boldu_+)_j=\max\{0,\boldu_j\}$ and
  $(\boldu_-)_j=\max\{0,-\boldu_j\}$. The \emph{toric ideal} $I_A$ is
  the $R_A$-ideal generated by all $\Box_\boldu$ with
  $\boldu\in\ker A$. Finally, the \emph{hypergeometric ideal} and
  \emph{module} to $A,\beta$ are
  \[
  H_A(\beta):=D_A(I_A,\{E_i-\beta_i\}_1^d),\qquad\qquad
  M_A(\beta):=D_A/H_A(\beta).
  \]
\schluss\end{dfn}
Before we embark on a general discussion of these modules we wish to
distinguish two special subclasses that will play a lead role.
\begin{dfn}\label{dfn-homogeneous}
  The matrix $A$ is \emph{homogeneous} if the following equivalent
  properties are satisfied:
  \begin{itemize}
  \item there is a group homomorphism from $\ZZ^d$ to $\ZZ$ that
    sends every $\bolda_j$ to $1\in\ZZ$;
  \item the vector $(1,1,\ldots,1)$ is in the row span of $A$;
  \item the ideal $I_A$ is standard graded and thus defines a
    projective variety inside projective $(n-1)$-space.
  \end{itemize}
\schluss\end{dfn}
\begin{dfn}\label{dfn-normal}
  The semigroup $\NN A$ is \emph{saturated} if $\NN A$ agrees with the
  intersection of $\ZZ A$ with the cone $\RR_{\geq 0} A$ spanned by
  the columns of $A$ viewed as elements of $\RR^n=\ZZ^n\otimes_\ZZ
  \RR$.
\schluss\end{dfn}

In a series of articles, including \cite{GGZ87,GKZ89,GKZ90},
I.M.~Gel{$'$}fand and his collaborators M.~Graev, M.~Kapranov and
A.~Zelevinsky developed the basic theory of these systems of linear
PDEs. The initial motivation came from Aomoto type integrals
\begin{equation}\label{eqn-hyp-int}
  Y(\beta;\boldx)=\int_Ct^\beta\exp\left(\sum_{i=1}^nx_it^{\bolda_i}\right)\frac{\de
    t_1}{t_1}\cdots\frac{\de t_d}{t_d}
\end{equation}
depending on a complex parameter vector $\beta\in\CC^d$, It is not
hard to verify that a hypergeometric function defined by the integral
\eqref{eqn-hyp-int} is annihilated by both the Euler operators and the
box operators \cite{GKZ90,Adolphson-duke94} but it took a decade to arrive at
the general formulation given here.

\medskip

It turns out that every univariate hypergeometric function arises
as a solution of an $A$-hypergeometric system; we sketch next the steps to
construct the proper $A,\beta$.
The general hypergeometric univariate differential equation is
\begin{eqnarray}\label{eqn-gen-hyp-eqn}
\prod_{v_j>0} \prod_{\ell=0}^{v_j-1} (v_j \theta_z +c_j -l)
&=& z \cdot \prod_{v_j<0} \prod_{\ell=0}^{|v_j|-1} (v_j \theta_z +c_j -l).
\end{eqnarray}

It is elementary, but not always trivial, to bring a differential
equation derived from a series expansion of a hypergeometric function
into this shape; it may require changes of variables in $z$. Note that
${}_pF_q(\alpha;\beta;z)$ is a solution to the special form
\begin{eqnarray}\label{eqn-gauss-diffeq}
\theta_z\prod_{j=1}^q(\theta_z+\beta_j-1)&=&z\cdot\prod_{j=1}^p(\theta_z+\alpha_j)
\end{eqnarray}
as one can see from applying the two operators to the power series
\eqref{eqn-pFq}.

Let $\boldv$ and $\boldc$ be the vectors with entries $v_j$ and $c_j$
respectively. For ${}_2F_1$ (equal to the function $F_1$ in Example
\ref{exa-elliptic}), $\boldv=(1,1,-1,-1)$ while for the Kummer
confluent function ${}_1F_1$, $\boldv=(1,1,-1)$.

Now, in order to manufacture $A$ and $\beta$ from equation
\eqref{eqn-gen-hyp-eqn}, choose an integral matrix $A$ such that
$\ZZ\cdot\boldv=\ker A$ and set $\beta=A\cdot\boldc$. Then the
solutions of $H_A(\beta)$ (in other words, the functions annihilated
by every operator in this left ideal) ``contain the solutions to
\eqref{eqn-gen-hyp-eqn}'' in the following sense.

\begin{exa}[The GKZ-system to the Kummer confluent function]
Consider the system of partial differential equations
\begin{eqnarray}
\label{eqn-erf-1}
\left(1\theta_1+\phantom{+1\theta_2+}1\theta_3\right)\bullet(u)&=&(-1/2)u\\
\label{eqn-erf-2}
\left(\phantom{+1\theta_1+}1\theta_2+1\theta_3\right)\bullet(u)&=&(0)u\\
\label{eqn-erf-3}
\left(\del_1\del_2-\del_3\right)\bullet(u)&=&0
\end{eqnarray}
in $x_1,x_2,x_3$.
This is the $A$-hypergeometric system to
\begin{equation}\label{1}
A=\begin{pmatrix}1&0&1\\0&1&1\end{pmatrix},\qquad
 \beta=\begin{pmatrix}-1/2\\0\end{pmatrix},
\end{equation}
since $v=(1,1,-1)$ is the
$\ZZ$-kernel of $A$.

Equation~(\ref{eqn-erf-1}) forces any solution $u$ to be homogeneous
(and of degree $-1/2$) under the grading that attaches the weights
$(1,0,1)$ to $(x_1,x_2,x_3)$. Similarly, Equation~(\ref{eqn-erf-2})
asserts that $u$ is homogeneous of weight zero if
$(x_1,x_2,x_3)\mapsto(0,1,1)$. It follows that one can write
\[
u(x_1,x_2,x_3)=x_1^ax_2^bx_3^c\, g(x_1x_2/x_3)
\]
where the monomial $x_1^ax_2^bx_3^c$ is of bi-degree $(-1/2,0)$, and
$g$ is a \emph{univariate} function.

Set $z=x_1x_2/x_3$ and write
\[
g(z)=\sum_{i=0}^\infty g_iz^i.
\]
Enforcing the vanishing of $\del_1\del_2-\del_3$ on $u(x_1,x_2,x_3)$
as suggested by Equation~(\ref{eqn-erf-3})
implies the recurrence relations
\[
(c-i)g_i=(a+i+1)(b+i+1)g_{i+1}
\]
for all $i$, and the starting condition
\[
\del_1\del_2\bullet(x_1^ax_2^b)=0.
\]
For $a=0$, observing that $x_1^ax_2^bx_3^c$ is of bi-degree $(-1/2,0)$,
we infer $b=-c=1/2$ and thus the recurrence is
\[
(-1/2-i)g_i=(i+1)(1/2+i+1)g_{i+1},
\]
showing that $g(z)$ essentially agrees with the Kummer confluent
function.
\schluss\end{exa}

\begin{exa}[GKZ-system to ${}_2F_1$]
Take the
equation \eqref{eqn-gauss-diffeq} with $p=q=2$ and
$\boldc=(1,c,a,b)$. Then
 $\boldv=(1,1,-1,-1)$ and the matrix
$A$ can be chosen as
\[
A=\begin{pmatrix}1&1&1&1\\1&0&0&1\\0&1&0&1\end{pmatrix},
\]
so that $\beta=A\cdot\boldc=(c-1,-a,-b)$.
The three Euler operators $\{\sum_{j=1}^4
a_{i,j}\theta_j-\beta_i\}_{i=1}^3$ annihilate each solution, so
every monomial $x^\boldu$ in the power series expansion of every
solution to the $A$-hypergeometric system must satisfy the three
conditions
\begin{eqnarray*}
  \left(u_1+u_2+u_3+u_4\right)&=&\beta_1;\\
  \left(u_1\phantom{{}+u_2+u_3}+u_4\right)&=&\beta_2;\\
  \left(\phantom{u_1+{}}u_2\phantom{{}+u_3}+u_4\right)&=&\beta_3.
\end{eqnarray*}
For a monomial $x^\boldu$, we call $A\cdot\boldu\in\ZZ A$ the
\emph{$A$-degree of $x^\boldu$}.  Then, every solution
$u(x_1,x_2,x_3,x_4)$ can be written as a univariate function $g$ in
$\frac{x_1x_4}{x_2x_3}$, multiplied by a monomial of $A$-degree
$\beta$. As in the previous example, one can use the fact that
$\Box_\boldv$ kills $u$ to show that $g$ satisfies the Gau\ss\ hypergeometric
differential equation.
\schluss\end{exa}

Of course, the kernel of $A$ being $\ZZ\cdot\boldv$ means that
$A\in\ZZ^{(n-1)\times n}$ and $I_A=(\Box_\boldv)$ is principal.  On
the other hand, the $A$-hypergeometric paradigm also encodes
multivariate hypergeometric series of higher rank (namely $n-d$) when
$d<n-1$. The solutions to $H_A(\beta)$ use $n$ variables and satisfy
$d$ homogeneities, so that effectively they are functions in $n-d$
independent quantities. Some aspects of the translation between the
two setups is discussed in \cite{BMW-torusInv}. The advantage of the
$A$-hypergeometric point of view is that it allows hypergeometric
functions to be studied with methods coming from algebraic geometry,
commutative algebra, and the theory of torus actions. We describe in
the following sections some of the advances and some of the new
problems that have been created through these new techniques.

\subsection{Solutions}
\label{subsec:Sols}

While we do not focus very much on solutions of $A$-hypergeometric
systems in this survey, it is only fair to indicate to some extent the
development of the understanding of their solution space over time. We also refer the reader to Remark \ref{rem:SolsIrreg} below, where we list and discuss some more references, after having explained issues like irregularity and slopes of hypergeometric systems.

Classically, functions were considered as hypergeometric if they could
be developed into a hypergeometric series. They typically arose from
specific differential equations and the hypergeometricity was a
consequence of the recurrence relations that came out of the
differential equation. While introducing $A$-hypergeometric systems,
Gel{$'$}fand and his collaborators Graev, Kapranov and Zelevinsky developed
a similar paradigm for the multi-variable homogeneous case, see
Definition \ref{dfn-homogeneous}. With setup as in
Section \ref{sec:EK-Rank}, so $A\cdot\gamma =\beta$ and $L_A$ the
kernel of $A$, the series
\[
\sum_{\bolda\in L_A}\boldx^{\gamma+\bolda}/\prod_{1\le j\le
n}\Gamma(\gamma_j+a_j+1)
\]
formally is a solution of $H_A(\beta)$. Assuming a certain amount of
genericity for $\gamma$ (such as  \emph{non-resonance}, see Definition
\ref{dfn-resonance})
the article \cite{GKZ89} also finds that the regions of convergence of these
series contain an open cone of the same shape as $(\RR_{\geq 0})^n$.

The series approach to solving differential equations of
hypergeometric type was then taken further by Sturmfels, Saito and
Takayama in their book \cite{SST00} through the technique of Gr\"obner
bases. As part of this mechanism, triangulations arise. The connection
between certain special solution series on one side and  and
triangulations on the other appears
already in \cite{GKZ89}.
In the homogeneous normal case (see Definition \ref{dfn-homogeneous})
it can be used to count the number of solutions
as the simplicial volume of
the convex hull of the columns of $A$;
\cite{SST00} provides various generalizations.

The first functions that were identified as hypergeometric were the
$\Gamma$-type integrals $\int t^a(1-t)^b(1-zt)^c\de t$ of Euler for the Gau\ss\
hypergeometric function. In \cite{GKZ90},
the authors consider integrals
\[
\int_\sigma\boldt^\beta \prod P_i(\boldt)^{\alpha_i}\de t_1\ldots
\de t_d
\]
where $P_i(\boldt)$ are Laurent polynomials and the integrals are
functions in the coefficients of the polynomials $P_i$. Here, $\sigma$ is a
$k$-cycle; in the Euler integrals $\sigma$ is a curve. Gel{$'$}fand,
Kapranov and Zelevinsky show that the above integrals are
$A$-hypergeometric and under suitable conditions span the solution
space. This approach generalizes Aomoto's integrals on complements of
generic hyperplane arrangements \cite{Aomoto-77}, a source of
inspiration in the search for the right definition of
$A$-hypergeometric systems.

There has always been a strong trend towards the study of ``special''
hypergeometric systems, namely those for which the solution space is
spanned by special classes of functions. This starts with Gau{\ss}'
observation \cite[page 125, Formel I.-V.]{Gauss-WerkeIII} that some
parameter choices in the Gau\ss\ hypergeometric differential equation
yield algebraic solutions. Kummer in \cite{Kummer-CrelleBd15-1}, Riemann,
and Gauss \cite[page 207]{Gauss-WerkeIII} developed tools to search
for other such instances. Then Schwarz constructed his famous list
\cite{Schwarz-Crelle1873} of the Euler--Gau\ss\ hypergeometric differential
equations whose solution space is spanned by algebraic functions. The
case of all $_{p}F_{p-1}$ was dealt with much later by Beukers and Heckman in
\cite{BeukersHeckman-Inventiones89} as part of their study of the
monodromy. For irreducible such equations with real parameters
$\alpha_1,\ldots,\alpha_p,\beta_1,\ldots,\beta_{p-1}$ set $
\beta_p=1$. Their exponentials on the unit circle are \emph{interlaced}
provided that the images of $\alpha_i$ and $\beta_j$ are encountered
alternatingly on a trip around the unit circle. Then
\cite{BeukersHeckman-Inventiones89} shows that interlacing is
equivalent to the solution space of the differential equation being
spanned by algebraic functions. Other cases were characterized in
\cite{Sasaki76,CohenWolfart92} (Appell--Lauricella $F_D$),
\cite{Kato-F2,Kato-F4} (Appell $F_2,F_4$).

For saturated irreducible
 homogeneous $A$-hypergeometric systems $M_A(\beta)$ with rational
 $\beta$, Beukers discovered the
 following fact about the number of algebraic solutions.
 Let $C_{A,\beta}=(\beta+\ZZ A)\cap (\RR_{\geq 0}A)$ and
 consider it as a module over the semigroup $\NN A$. Let
 $\sigma_A(\beta)$ be the number of generators of $C_{A,\beta}$ over
 $\NN A$. Then, Beukers shows in \cite{Beukers-Inventiones10} that
 $\sigma_A(\beta)$ never exceeds the volume of $A$, and equality of
 $\sigma_A(k\beta)=\vol(A)$ for all $1\le k\le D$ coprime to the
 least common denominator $D$ of $\beta_1,\ldots,\beta_d$ happens
 precisely when the solution space is spanned by algebraic
 functions. We remark that
 irreducibility is linked to non-resonance (compare Definition
 \ref{dfn-resonance}) by
 \cite{Beukers-resonance,Saito-compositio11,SchulzeWalther-resgkz}.

The story for inhomogeneous (\emph{i.e.}, confluent) systems is more complicated, both
theoretically and algorithmically. Since the solutions do not need to
lie in the Nilsson ring, a systematic search in the sense of
\cite{SST00} using Gr\"obner bases is not possible. Nonetheless, in
\cite{EsterovTakeuchi-AJM15} an idea of Adolphson
\cite{Adolphson-duke94} is completed that casts
solutions of non-resonant $A$-hypergeometric systems as integrals
\[
 \int_{\gamma^z}\exp{\left(\sum_{j=1}^nx_j\boldt^{\bolda_j}\right)}t_1^{c_1-1}\cdots
 t_d^{c_d-1}\de t_1\ldots \de t_d.
\]
Here, $\gamma$ is a continuous family of real $d$-dimensional
topological cycles in the torus, on which the integrand decays rapidly
at infinity in the sense of Hien \cite{Hien-Inventiones09}. This was
also already studied in the
context of integrals from hyperplane arrangements by
\cite{KimuraHaraokaTakano92}.

\section{Torus action and Euler-Koszul complex}
\label{sec:EK-Rank}

In this section we start exploring algebraic properties of the system
$H_A(\beta)$ by introducing a homological tool from \cite{MMW05} that
has proved to be very successful: the Euler--Koszul complex. It has
been used to study the number of solutions, their monodromy, and
several other aspects. We refer to the start of Subsection
\ref{subsec-dim-red} for basic notations and assumptions regarding $A$.

\subsection{Torus action and $A$-grading}\label{subsec-torusaction}

Given a $D_A$-module $Q$, its \emph{Fourier--Laplace transform}
$\widehat{Q}$ is
equal
to $Q$ as a $\CC$-vector space and carries a $\widehat{D}_A :=
\CC[\boldxi]\langle \bolddel \rangle$ structure given by
\begin{equation}\label{eq:FLconcrete}
\xi_j \cdot m := \del_{x_j} \cdot m, \qquad  \del_{\xi_j} \cdot m  :=  - x_j \cdot m,
\end{equation}
for any $m\in Q$. See \eqref{eq:FL} for a functorial
description, and compare Subsection \ref{subsec-weight} for a related construct,
the Fourier--Sato transform.

The polynomial ring $R_A$ is naturally identified with the coordinate
ring $\CC[\boldxi] $ of the Fourier--Laplace dual space $\widehat{\CC}^n$ of $\CC^n$.  The matrix
$A$ defines an algebraic action
\[
\TT\times \widehat{\CC}^n \to \widehat{\CC}^n
\]
of the \emph{$d$-torus}
\[
\TT:=(\CC^*)^d=\Spec(\CC[t_1^{\pm1},\dots,t_d^{\pm1}])
\]
with coordinates $\boldt=t_1,\dots,t_d$ on $\widehat{\CC}^n$ by
\begin{equation}\label{16}
(\eta,\xi)\mapsto \eta\cdot\xi:=(\eta^{\bolda_1}\xi_1,\dots,\eta^{\bolda_n}\xi_n).
\end{equation}
This action induces a grading
\[
R_A=\bigoplus_{\bolda\in\ZZ A} (R_A)_\bolda
\]
on $R_A$, where
\[
\deg(\del_j)=\bolda_j;
\] we refer to this as the  \emph{$A$-grading}. There is a natural
extension to $D_A$ if one sets
\[
\deg(x_j)=-\bolda_j
\]
that makes every Euler operator $A$-graded of degree zero.

The coordinate ring of the orbit closure
through $(1,\dots,1)$ is the \emph{toric ring}
\[
S_A:=\CC[t^{\bolda_1},\cdots,t^{\bolda_n}]=\CC[\NN A]=R_A/I_A.
\]
\begin{rmk}
  The semigroup ring $S_A$ is
  normal (and hence Cohen--Macaulay by Hochster's Theorem 1 in \cite{Hochster-tori}) if and only
  if $\NN A$ is saturated in the sense of Definition \ref{dfn-normal}.
\schluss\end{rmk}
We shall identify subsets of columns of $A$ with subsets of column
indices or submatrices. For such a subset $\tau \subset A$, set
\[
(\boldone_\tau)_j:=
\begin{cases}
1&\text{if}\quad\bolda_j\in\tau,\\
0&\text{if}\quad\bolda_j\notin\tau,
\end{cases}
\]
denote by $O_A^\tau$ the orbit of $\boldone_\tau$, and its Zariski closure by $\bar O_A^\tau$. Moreover, we write $S_A^\tau$ for the coordinate ring of $\bar O_A^\tau$.

Let $I_A^\tau$ be the $R_A$-ideal generated by $I_A$ and all $\bolddel^\boldu$ with
$A\cdot\boldu\not\in\tau$. It is $A$-graded and prime and we have $S_\tau=R_A/I_A^\tau$. Note that
\[
O_A^\tau=\Var(I_A^\tau)\minus
\bigcup_{\tau'\subsetneq \tau}\Var(I_A^{\tau'}),
\]
with $\dim(\tau)=\dim(\Var(I_A^\tau))=\dim(O_A^\tau)$.

The following sets are then in one-to-one correspondence:
\[
\left\{\text{faces $\tau$ of $\RR_{\geq 0}\cdot
  A$}\right\}\leftrightarrow \left\{\text{$A$-graded primes
  $I_A^\tau\supseteq I_A$ of $R_A$}\right\}\leftrightarrow
\left\{\text{$\TT$-orbits $O^\tau_A$}\right\}.
\]

\subsection{Toric category and Euler--Koszul technology}

The following set of constructions and results is taken from \cite{MMW05}.

Note that $E_i-\beta_i\in D_A$ can be viewed as a
left $D$-linear endomorphism on $A$-graded $D_A$-modules $M$ by sending
a $\ZZ A$-homogeneous $y\in M$ to
\begin{equation}\label{86}
(E_i-\beta_i)\circ y:=(E_i-\beta_i-\deg_i(y))y,
\end{equation}
and that these morphisms commute with one another.

\begin{dfn}[Degrees and Euler--Koszul complex]\label{19}
  Let
  \[
  N=\bigoplus_{\bolda\in\ZZ A}N_\bolda
  \]
  be an $A$-graded $R_A$-module and pick $\beta\in\CC^d$.  Let
  $\tdeg_A(M)$ be the \emph{true $A$-degrees of $N$}, given as the set
  of points $A\cdot\boldu$ in $\ZZ A$ for which the graded component
  $N_\boldu$ is nonzero,
  \[
  \tdeg_A(N):=\{\bolda\in\ZZ^d\mid N_\bolda\neq 0\}.
  \]
  Write $\qdeg_A(N)$ for the Zariski closure of $\tdeg_A(N)\subseteq
  \ZZ A$ inside $\CC^d$.

  The
\emph{Euler--Koszul complex} $K_{A,\bullet}(N;\beta)$ is the Koszul
complex of the endomorphisms $E-\beta$ on the left $D_A$-module
$D_A\otimes_RN$ equipped with the natural $A$-grading.  Its $i$-th homology
$$
H_{A,i}(N;\beta):=H_i(K_{A,\bullet}(N;\beta))
$$
is the $i$-th
\emph{Euler--Koszul homology} of $N$.  Note that
$H_{A,0}(S_A;\beta)=M_A(\beta)$.
\schluss\end{dfn}

\begin{rmk}
A (commutative graded) precursor of the Euler--Koszul complex when $N=S_A$
appears already in \cite{GKZ89} for proving holonomicity of $M_A(\beta)$ when
$S_A$ is a Cohen--Macaulay ring, and in Adolphson
\cite{Adolphson-duke94,Adolphson99} a modified version of the complex is
discussed.
\schluss\end{rmk}

The properties of the Euler--Koszul complex are most pleasant when
$N$ is in the category of {\em toric modules}. These are $A$-graded
$R_A$-modules that have a finite composition series whose successive
quotients are $\ZZ A$-shifted quotients of $S_A$.

\begin{rmk}
  There is a generalization in \cite{SchulzeWalther-ekdi} to
  \emph{quasi-toric} (\emph{i.e.}, certain
  non-Noetherian $A$-graded) modules that is useful for the
  interplay of Euler--Koszul complexes on local cohomology modules or
  on localizations such as $\CC[\ZZ A]$.

  A different generalization (\emph{toral} modules) is given and used
  in \cite{DMM-duke}.
\schluss\end{rmk}

By \cite{MMW05}, short exact sequences $0\to N'\to N\to N''\to 0$ of
toric modules give rise to long exact sequences of Euler--Koszul
homology modules that are all holonomic (see Definition \ref{dfn-holonomic}).
Moreover, vanishing of $H_{A,0}(N;\beta)$ implies vanishing
of all $H_{A,i}(N;\beta)$ and this vanishing is equivalent to $-\beta$
not being in the quasi-degrees of $N$.

\begin{rmk}
  While Euler--Koszul complexes were initially defined for the study of the
  size of the solution space of $A$-hypergeometric systems
  \cite{MMW05}, they have turned out to be remarkably successful when
  investigating other issues such as irregularity
  (see Section \ref{sec:Irreg} and \cite{SchulzeWalther-duke}), reducibility of the monodromy
  \cite{Walther-taka,MariCruz-RMI19}, comparisons with direct image
  functors (see the next subsection as well as \cite{SchulzeWalther-ekdi,Avi-JA19,Avi-JPAA19}), more
  general classes of binomial $D$-modules
  \cite{DMM-duke,BMW-torusInv,BMW-holSing}, the study of Horn hypergeometric
  systems \cite{DMM-duke,BMW-normalHorn}, resonance
  \cite{SchulzeWalther-resgkz}, or Hodge theoretic
  aspects (see sections \ref{sec:Hodge} and \ref{sec:MirrorSym} as well as \cite{Reich2,ReiSe,ReiSe2,ReiSe-Hodge,RW-weight}).
\schluss\end{rmk}

\subsection{Fourier--Laplace transformed GKZ-systems}\label{subsec-FGKZ}

 We noted in section \ref{subsec-torusaction} that the torus $\TT$
 acts on the Fourier--Laplace dual space $\widehat{\CC}^n$. The orbit closure through $(1,\ldots,1)$ is an affine toric variety $X_A := \Spec(S_A)$. We identify its dense open orbit $O_A$ with the torus $\TT$. This gives rise to inclusions
 \[
  \TT \overset{j_A}\longrightarrow  X_A \overset{i_A}\longrightarrow \widehat{\CC}^n
 \]
 where $j_A$ is an open embedding and $i_A$ is a closed embedding. We set
 \begin{gather}\label{eqn-h_A}
 h_A := i_A \circ  j_A.
 \end{gather}

We denote the Fourier--Laplace transform of $M_A(\beta)$ by
$\widehat{M}_A(\beta)$, and the corresponding quasi-coherent sheaves by
$\calM_A(\beta)$ and $\widehat{\calM}_A(\beta)$ respectively. Using
the definition of the Fourier--Laplace transform \eqref{eq:FLconcrete} one easily sees that
$\widehat{\calM}_A(\beta)$ has support on the toric variety $X_A$. In
\cite{SchulzeWalther-ekdi} the parameters $\beta$ were
identified for which  there is an isomorphism
$\widehat{\calM}_A(\beta)\simeq (h_{A})_+ \calO_{\TT}^\beta$ between
the Fourier--Laplace transform of $\calM_A(\beta)$ and
the direct image under $h_A$ of the twisted structure sheaf
\[
\calO_{\TT}^\beta = \calD_{\TT} / \calD_{\TT} ( \del_{t_1}t_1+\beta_1, \ldots , \del_{t_d}t_d + \beta_d).
\]
The relevant definition is the following one.
 \begin{dfn}\label{def-sRes}\cite{SchulzeWalther-ekdi}
  The elements of
  \[
   \sRes(A) := \bigcup_{j=1}^n \sRes_j(A)
  \]
  where
  \[
  \sRes_j(A) := \{\beta \in \CC^d \mid \beta \in -(\NN+1)\bolda_j + \qdeg_A(S_A / (t^{\bolda_j}))\}
  \]
  are the \emph{strongly resonant parameters} of $A$.
\schluss\end{dfn}
Strong resonance, as the language suggests, is a strengthening of
resonance, defined next.
\begin{dfn}\label{dfn-resonance}
  The parameter $\beta$ is \emph{resonant} for $A$ if
  $\beta+\ZZ^d$ meets the complexified boundary hyperplanes of the cone
  $\RR_{\geq 0}A$.
\schluss\end{dfn}
\begin{rmk}

Strongly resonant parameters are resonant.
The resonant parameters contain $\NN A$, but the strongly
resonant ones usually do not. For example, if the semigroup $\NN A$ is saturated, then
$\NN A \cap \sRes(A) = \emptyset$. In particular, 0 is not an element of $\sRes(A)$ in this case, a fact
that will become useful later.
\schluss\end{rmk}

\begin{exa}
Consider the matrix
\[
A = \left(\begin{matrix}
-1 & 0 & 1 & 2 \\
1 & 1 & 1 & 1
\end{matrix} \right)
\]
the sets $\tdeg_A(S_A)$ and $\sRes(A)$ and the cone $\RR_{\geq 0}A$
are sketched below. Since $d=2$, fullness of $A$ implies that we have
$\qdeg_A(S_A)=\CC^2$.
\begin{figure}[h]
  \caption{Cone, true, and strongly resonant degrees.}
\begin{center}
\newdimen\scale
\scale=0.9cm
\begin{tikzpicture}
 \filldraw[cyan,opacity=.4] (0,0) -- (\scale*-2.05,\scale*2.05) -- (\scale*4.05,\scale*2.05)  -- (0,0);

 \foreach \x in {-3,...,4}{
   \foreach \y in {-2,...,2}{
     \node[draw,circle,inner sep=0.5pt,fill] at (\scale*\x,\scale*\y) {};
   }
 }
\draw (\scale*-3,\scale* -2) -- (\scale*4, \scale *1.5);
\draw (\scale*-2,\scale* -2) -- (\scale*4, \scale *1);
\draw (\scale*-1,\scale* -2) -- (\scale*4, \scale *0.5);
\draw (\scale*0,\scale* -2) -- (\scale*4, \scale *0);
\draw (\scale*1,\scale* -2) -- (\scale*4, \scale *-0.5);
\draw (\scale*2,\scale* -2) -- (\scale*4, \scale *-1);
\draw (\scale*3,\scale* -2) -- (\scale*4, \scale *-1.5);
\draw (\scale*-3,\scale* 2) -- (\scale*1, \scale *-2);
\draw (\scale*-3,\scale* 1) -- (\scale*0, \scale *-2);
\draw (\scale*-3,\scale* 0) -- (\scale*-1, \scale *-2);
\draw (\scale*-3,\scale* -1) -- (\scale*-2, \scale *-2);

 \foreach \y in {0,...,2}{
   \foreach \x in {-\y,...,\y}{
     \node[draw,circle,inner sep=1.5pt,fill] at (\scale*\x,\scale*\y) {};
    }
  }

 \node[draw,circle,inner sep=1.5pt,fill] at (\scale*2,\scale*1) {};
 \node[draw,circle,inner sep=1.5pt,fill] at (\scale*3,\scale*2) {};
 \node[draw,circle,inner sep=1.5pt,fill] at (\scale*4,\scale*2) {};

\node[draw,circle,inner sep=1.5pt,fill] at (-0.6,-2.35) {};
\node[] at (0.4,-2.39) {$\tdeg_A(S_A)$};
\draw(2.0,-2.35)--(2.5,-2.35);
\node[] at (3.5,-2.39) {$\sRes_A(S_A)$};
\filldraw[cyan,opacity=.4] (-3.3,-2.45) rectangle (-2.7,-2.25);
\node[] at (-2.1,-2.41) {$\RR_{\geq 0} A$};

\end{tikzpicture}
\end{center}
\end{figure}
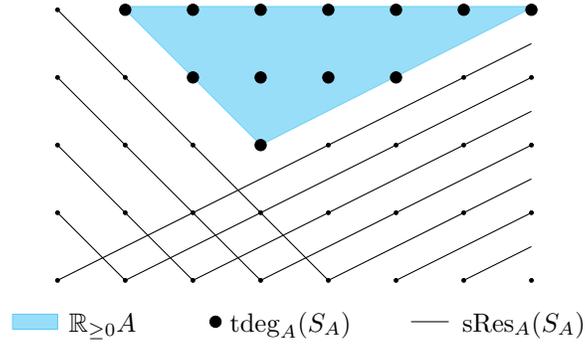
\schluss\end{exa}

\begin{thm}\label{thm:FL-GKZ-DirectImage}
Let $A\in \ZZ^{d \times n}$ be as above, then the following statements are equivalent
\begin{enumerate}
\item $\beta \not \in \sRes(A)$
\item $\widehat{\calM}_A(\beta) \simeq (h_A)_+ \calO_{\TT}^\beta$
\item Left multiplication with $\xi_i$ is invertible on $\widehat{M}_A(\beta)$.\qed
\end{enumerate}
\end{thm}
\begin{rmk}
  The idea of linking $\widehat{\calM}_A(\beta)$ to the direct
  image $(h_A)_+ \calO_{\TT}^\beta$ originates with
  \cite{GGZ87} where it was shown that $\beta$ non-resonant gives
  the desired isomorphism.  The precise computation in Theorem
  \ref{thm:FL-GKZ-DirectImage} comes from
  \cite{SchulzeWalther-ekdi}.  These results were refined and
  extended to the strongly resonant case in
  \cite{Avi-JA19,Avi-JPAA19} where Steiner uses a combination of
  direct and proper direct image functors.
  \schluss\end{rmk}

\subsection{Holonomicity, Rank, and Singular Locus}

Suppose $M=D_A/I$ is some left $D_A$-module, and $\calM=\calD_{\CC^n}/\calI$ the associated sheaf of $\calD_{\CC^n}$-modules. Then its analytification
$\calM^{\an}=\calD_{\CC^n}^{\an}/\calD_{\CC^n}^{\an}\calI$ is obtained by replacing $\calD_{\CC^n}$
by the sheaf $\calD_{\CC^n}^{\an}$ of analytic linear differential operators  on $\CC^n$ where now $\calI\subset \calD_{\CC^n}\subset\calD^{\an}_{\CC^n}$ generates a left ideal of analytic linear differential operators.

Choose $\pointx\in\CC^{n}$ and denote stalks by subscripts.  Consider the functor
\[
\Sol_\pointx(-)=\Hom_{\calD^{\an}_{\CC^n,\pointx}}(-,\calO_{\CC^n,\pointx}^{\an})
\]
from germs of left $\calD_{A,\pointx}^{\an}$-modules to vector spaces.\footnote{Notice that we do not consider derived solutions here;
  so our use of the symbol $\Sol$ differs from many other texts on $\calD$-modules.} If
$\calM^{\an}=\calD_{\CC^n}^{\an}/\calD_{\CC^n}^{\an}\calI$ then $\eta\in
\Sol_\pointx(\calM^{\an})$ corresponds to the analytic solution
$\eta(1+\calD_{\CC^n}^{\an}\calI)$ near $\pointx$. The dimension of the vector
space of solutions to $\calM$ at $\pointx$ is \emph{the rank of $M$ at
  $\pointx$}. When we mean the rank at a generic point $\pointx$ we
speak of just \emph{the rank of $M$}.

Typically, $\Sol_\pointx(\calM^{\an})$ is infinitely generated. But for the select
class of \emph{holonomic modules} it is always finite.
\begin{dfn}\label{dfn-holonomic}
  Any principal $D_A$-module (resp.\ $\calD_{\CC^n}^{\an}$-module) $M$
  (resp.\ $\calM$) with generator $m$ has a natural \emph{order
    filtration $F^\ord_\bullet$} by $R_A$-modules (resp.\ $\calO_{\CC^n}$-modules)
  where $F^\ord_k(M)$ (or, on the stalk, $F^\ord_k(\calM_\pointx)$) is generated by
  the cosets of $\bolddel^\boldu$ with $|\boldu|\le k$. The notion readily
  extends to any module with chosen set of generators and
  behaves well under analytification.

  If $\calM=\calD_{\CC^n}^{\an}$ is the sheaf of differential operators itself,
  the associated graded object is on the stalk isomorphic to the
  regular ring
  $\calO_\pointx[\boldy]$ where $\boldy=y_1,\ldots,y_n$ is the set of
  symbols to $\del_1,\ldots,\del_n$. For any $M$ (resp.\ $\calM$), the
  associated graded object $\gr^F(-)$ becomes a module over
  $\gr^F(D_A)$ (resp.\ $\gr^F(\calD_{\CC^n}^{\an})$).

  The module is \emph{holonomic} if the associated graded module has
  Krull dimension $n$.  \schluss\end{dfn} It was shown in
  \cite{GGZ87,GKZ89} that many, and then in \cite{Adolphson-duke94}
  that in fact all $A$-hypergeometric systems are holonomic; an
  elementary proof is given in \cite{BGM-pacific}. Holonomicity was
  then extended in \cite{MMW05,SchulzeWalther-ekdi} to all
  Euler--Koszul homology modules derived from quasi-toric input.

\medskip

By \cite{SKK-LNM73,Gabber81}, the characteristic variety is always involutive
and has all components of dimension $n$ or larger. This implies that
holonomic modules have finite length and satisfy a
Krull--Remak--Schmidt theorem (have well-defined sets of simple composition
factors with multiplicity taken into account). Moreover, the quantity
\[
\rk(M):=\dim_\CC(\CC(\boldx)\otimes_{\CC[\boldx]}M)
\]
agrees with the rank of $M$ in a generic point $\pointx\in\CC^n$ by the
Cauchy--Kovalevskaya--Kashiwara Theorem \cite[p.~37]{SST00}.

For many important $A$-hypergeometric systems, a search of explicit
natural power series solutions leads to rank many independent
solutions, compare \cite{GGZ87,SST00}. It was claimed in \cite{GKZ89}
that the rank of $M_A(\beta)$ is
\[
\rk(M_A(\beta))=\vol(A),\]
where $\vol(A)$  is the (simplicial) \emph{volume} of $A$, a
purely  combinatorial quantity given by the
quotient of the measure of the convex hull of the origin and the columns of $A$, divided by the measure of the standard $n$-simplex.
Adolphson \cite{Adolphson-duke94} pointed at a possible
flaw in the argument, and \cite{ST98} eventually provided a
counter-example that is worth looking at.

\begin{exa}[The 0134-curve, \cite{ST98}]\label{exa-0134}
Let $A=\begin{pmatrix}1&1&1&1\\0&1&3&4\end{pmatrix}$. The volume of
$A$ is 4, equal to the volume of the interval $(0,4)$ inside
$\RR$. (Since the interval is $1$-dimensional, usual
volume---length---and simplicial volume agree).

The toric ideal $I_A$ is homogeneous here, defining the pinched
rational normal space curve. In \cite{SST00} it is shown that series
solution methods based on weight vectors and the computation of
certain initial ideals of $H_A(\beta)$ always lead to volume many
independent series solutions, as long as $A$ is homogeneous.  This
generalized the na\"ive series written out in \cite{GGZ87,GKZ89} to
the case where logarithmic terms can appear in the series solutions.

For almost all $\beta$, the rank of $M_A(\beta)$ in a generic
point is $4$, spanned by functions
\begin{eqnarray*}
x_1^{(4\beta_1-\beta_2)/4}x_4^{\beta_2/4}+\ldots,\qquad
x_1^{(4\beta_1-\beta_2-3)/4}x_2x_4^{(\beta_2-1)/4}+\ldots,\\
x_1^{(4\beta_1-\beta_2-1)/4}x_3x_4^{(\beta_2-3)/4}+\ldots,\qquad
x_1^{(4\beta_1-\beta_2-6)/4}x_2^2x_4^{(\beta_2-2)/4}+\ldots,
\end{eqnarray*}
where the dots indicate a (usually infinite) series of terms ordered
by the weight vector $(0,1,2,0)$. (The particular weight is
immaterial, but it needs to be sufficiently generic; this one is so for
this example).
If one now deforms $\beta$ into $(1,2)$ then the four independent
solutions above degenerate into a linearly dependent set of rank
three. On the other hand, the functions
\[
\frac{x_2^2}{x_1},\qquad\qquad \frac{x_3^2}{x_4}
\]
are new, not-deforming (in $\beta$) solutions to $M_A((1,2))$. It
follows that the ``rank jumps at $\beta=(1,2)$'', from 4 to $5=4-1+2$.
\schluss\end{exa}

Shortly after the discovery of rank jumps, the case of homogeneous
monomial curves was completely discussed in \cite{CDA-duke}: the
``holes'' of $\NN A$ (the finitely many elements of $(\RR_{\geq
  0}A\cap\ZZ A)\minus \NN A$) are exactly the rank-jumping parameters,
and each rank jump is by $1$.  It was then shown in \cite{MMW05} that
as $\beta$ varies, the rank of $M_A(\beta)$ is upper-semicontinuous,
so that it can only go up under specialization (formation of a limit)
of $\beta$. In fact, \cite[Cor.~9.3]{MMW05} shows that the
\emph{exceptional set} $\calE_A$ of points where rank exceeds volume
is Zariski closed and equals a certain subspace arrangement. To
understand the origins of $\calE_A$ one must view the local cohomology
modules $H^i_\bolddel(S_A)$ with $i<d$ as quasi-toric modules; their
elements are then witnesses to the failure of $S_A$ to be
Cohen--Macaulay, while the union of their quasi-degrees forms the
exceptional arrangement.  The fact, also observed in \cite{MMW05},
that this arrangement has codimension at least two explains why
finding rank-jumps at all turned out to be very hard and involved
extensive computer experiments in \cite{ST98}.

\begin{exa}[Continuation of Example \ref{exa-0134}]
  In Example \ref{exa-0134}, $d=2$ and
so $\calE_A$ can be at most a finite set of isolated points. The local
cohomology $H^0_\bolddel(S_A)$ is zero and $H^1_\bolddel(S_A)$ is a
1-dimensional vector space generated by the \v Cech cocycle
$(\del_2^2/\del_1,\del_3^2/\del_4)$. To see this,  note that
$(\del_1,\del_4)$ is primary to $\bolddel$ in $S_A$. Thus,
$H^1_\bolddel(S_A)$ can be computed $A$-degree by $A$-degree from
the \v Cech complex on $S_A$ induced by $\del_1,\del_4$. Each degree
component in $S_A$ and its monomial localizations are $1$-dimensional
$\CC$-spaces; we use this to depict these localizations in the \v Cech
complex
by dots as follows:
\begin{figure}[H]
  \centering
  \caption{The \v Cech complex to the 0134-curve}
  \psfrag{o}{$\oplus$}
  \psfrag{l}{$\CC[\NN A]$}
  \psfrag{m}{$\CC[\NN (A,-\bolda_1)]\oplus\CC[\NN (A,-\bolda_4)]$}
  \psfrag{r}{$\CC[\ZZ A]$}
  \includegraphics[width=0.95\textwidth]{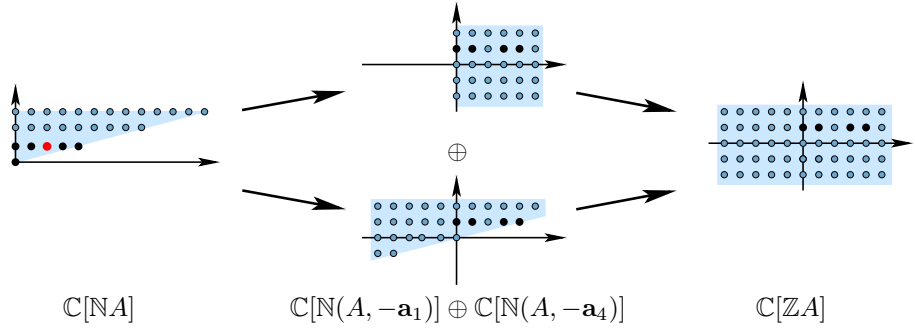}
\end{figure}
In this picture, the blue area indicates the directions in which the
semigroup in question extends, black dots are the elements of $A$ and
the red dot indicates a ``missing'' element in the semigroup.
Taking cohomology ``dot-by-dot'' one identifies the local
cohomology groups $H^1_\frakm(S_A)$, $H^0_\frakm(S_A)$ as claimed.

It is remarkable
that the components of the $H^1_\frakm(S_A)$-cocycle are precisely the ``new''
solutions that appear at $\beta=(1,2)$ that do not deform to other
$\beta$. While this is not always literally true, a weaker form is
typical and an explanation of this phenomenon involving Laurent
polynomials is given in
\cite{
BerkeschForsgardMatusevich-TAMS,
BerkeschForsgardMatusevich-Adv16}, especially for $d=2$. Compare also
Remark \ref{rem:SolsIrreg}.
\schluss\end{exa}

\begin{rmk}In \cite{Berkesch-rank} it is
proved that there is a purely combinatorial recipe (involving the
relative positioning of $\beta$ to the degrees of $\NN A$)
that determines the rank of $M_A(\beta)$. The procedure to arrive at
the exact rank is very involved.

The only known closed rank formula is for non-jumping
parameters, where the rank is just the volume.\footnote{This is not
  entirely true: if $d=2$ and the columns of $A$ lie in a hyperplane
  not containing the origin, then all rank jumps are by 1,
  as shown in  \cite{CDA-duke}.} The best known general
bound is exponential \cite{SST00}, in the sense that the rank of
$M_A(\beta)$ is bounded above by $2^{2d}\vol(A)$.
It was shown in
\cite{MatusevichWalther} that for every $d$ there are rank jump examples with
$\rk(M_A(\beta))=\vol(A)+d-1$.
This is improved in
\cite{Fernandez-jumps} to the existence of $a\in \RR$ greater than $1$ and
families of matrices $A_{(d)}$
of size  $d\times n_d$ and with
parameters $\beta_{(d)}$ such
that the rank of $M_{A_{(d)}}(\beta_{(d)})$ exceeds $a^d\vol(A)$.
It would be interesting to know how far the bound from \cite{SST00} is
from the the worst examples that exist.
\schluss\end{rmk}

\medskip

There is an open subset of $\CC^n$ on which the solutions for
$M_A(\beta)$ form a vector bundle of rank $\rk(M_A(\beta))$. The
complement (the \emph{singular locus of the module}) of
this set is algebraic, cut out by the \emph{$A$-discriminant}, a
product of individual discriminants to polynomial systems, one for
each face of the cone over $A$. For a very detailed discussion on
this, see the books \cite{GKZ-book}, and \cite{SST00}. If one moves
from general to special $\pointx$, rank can go down due to singularities in
the solutions.  In contrast to rank in generic points, rank at special
$\pointx$ is not known to be upper-semicontinuous. For the case of $A$ as
in Example \ref{exa-0134}, this is worked out in \cite{W-exp}, which
discusses the more general question of stratifying $\CC^n$ by the
\emph{restriction diagrams}, which encode the behavior of the
$D$-module theoretic (derived) pull-back to $\pointx\in\CC^n$; the
elementary pull-back just counts rank at $\pointx$.

\subsection{Better behaved systems and contiguity}\label{subsec:Contiguity}

For each $\beta'=\bolda_j+\beta$ there is a  natural \emph{contiguity morphism}
\[
  c_{\beta,\beta+\bolda_j}\colon M_A(\beta)\stackrel{\del_j}{\to} M_A(\beta')
\]
of degree $\bolda_j$, induced by right multiplication with $\del_j$ on
$S_A$ through the Euler--Koszul functor. The existence of these
morphisms is a consequence of the fact that
$(E_i-\beta_i)\cdot\del_j=\del_j(E_i-\beta_i-a_{i,j})$; this is a
special case of Equation \eqref{86} when $y=\del_j$.  Since elements
in $I_A$ act as zero on $S_A$, any composition of contiguity morphisms
of fixed total degree $\gamma\in\NN A$ acts the same way as morphism
$c_{\beta,\beta+\gamma}$ from $M_A(\beta)$ to $M_A(\beta+\gamma)$.

Contiguity morphisms
have turned out to be a very useful tool in the study of
$A$-hypergeometric systems since for $k\gg 0$,
$c_{\beta+k\bolda_j,\beta+(k+1)\bolda_j}$ and $c_{\beta-(k+1)\bolda_j,\beta-k\bolda_j}$ are
isomorphisms (and one can determine explicit bounds in terms of $A,\beta$
for $k$ being sufficiently big). Contiguity maps have been used in
\cite{Saito-compositio01} to identify combinatorially the isomorphism
classes of $A$-hypergeometric systems, in \cite{Walther-taka} to study
irreducibility and holonomic duality of $M_A(\beta)$ as a $D_A$-module,
and in \cite{Reich2, ReiSe-Hodge} for investigating
the Hodge module structure on certain $\calM_A(\beta)$. For a study of
Gau\ss\ hypergeometric functions via contiguity operators see \cite{Beukers07}.

On the level of solutions, a map in the reverse direction is induced
that literally takes the derivative by $x_j$. For certain applications
in mirror symmetry it is desirable to know that every contiguity
operator induces an isomorphism on (the solutions of) $M_A(\beta)$. In
case one has a generic $\beta$, this is automatic. But in practical
situations it is more likely that $\beta$ is integer, or at least resonant.
In the present context, resonance
encapsulates the lack of genericity of a parameter
$\beta$ to admit contiguity isomorphisms (in both
directions). Resonance and contiguity operators were refined and used
in \cite{Adolphson-duke94,Saito-compositio01,Saito-compositio11,
  Okuyama-Tohoku06,CatDickRod-imrn11,SchulzeWalther-resgkz,Beukers-resonance,
  Beukers-crelle16} to study reducibility and general structure of
$M_A(\beta)$.

Now consider the quasi-toric module $F_A$ equal to the
ring $\CC[\ZZ A]$. It arises  as the localization of $S_A$ at all
$\del_j$, or alternatively at one monomial whose degree is in the
interior of $\RR_{\geq 0}A$. By definition, multiplication by $\del_j$
on $F_A$ is an isomorphism, and therefore the same applies to the
generalized $A$-hypergeometric system that arises as the Euler--Koszul
homology $H_{A,0}(F_A;\beta)$, for every $\beta$. Since $F_A$ is a
maximal Cohen--Macaulay $S_A$-module, there is no other Euler--Koszul
homology, \cite{MMW05,SchulzeWalther-ekdi}.

This module $H_{A,0}(F_A;\beta)$ was studied in
\cite{BorisovHorja-MathAnn13,BorisovHorja-Adv15} and termed
\emph{better behaved GKZ-system}.  A variant of these systems,
considered in \cite{Mo15}, can be described as the Euler-Koszul
homology $H_{A,0}(\CC[\RR_{\geq 0}A\cap\ZZ^d];\beta)$ of the
normalization of $S_A$. In Section \ref{sec:Hodge} below we will
discuss Hodge theoretic ramifications of the main result of
\cite{Mo15}.

\section{Irregularity}
\label{sec:Irreg}

In this section we discuss regularity issues of hypergeometric
$D$-modules; this is a multi-variate form of essential
singularities. We start with discussing more general filtrations than
the one by order. A combinatorial object can be derived from this
process that governs the convergence behavior of solutions to
$A$-hypergeometric systems near coordinate hyperplanes. Via results of
Laurent and Mebkhout we discuss a generalized classical Fuchs criterion
this gives information on the irregular solutions.

\subsection{The Fuchs criterion and regularity}\label{subsec-Fuchs}

A univariate function $f(t)$, analytic on a small open disk around
$t=0$ but singular at $t=0$, can behave in two
essentially different ways: the growth of $f(t)$ as $t\rightarrow 0$
could be bounded by a polynomial, or not. In the former case, $f$ has
a pole, in the latter an essential singularity. If $f$ arises as
solution to a differential equation we say $0$ is a \emph{regular
  singular point} of the equation in the first, and an irregular singular
point in the second case.

For linear differential equations $P\bullet f(z)=0$ in the local
parameter $z$, Fuchs gave the following practical procedure for
determining regularity of the origin. If $\calO_0:=\CC\{z\}$ is
the ring of convergent power series near $z=0$, write $P$ as a linear
combination
\[
P=\sum_{k=0}^m p_k(z)\cdot \frac{\del^k}{\del z^k},
\]
$m$ being the order of $P$, and $p_k=\sum_{i=n_k}^\infty c_{k,i}z^i\in
\calO_0$ with $c_{k,n_k}\neq 0$ indicating the lowest order term of
$p_k(z)$.  Writing $\del_z$ for differentiation by $z$, for a monomial
$z^r\del_z^s$ we use the two weights
\begin{align*}
  &&V(z^r\del_z^s)&:=s-r&&&&{\text{$V$-filtration at $0$};}\\
  &&F(z^r\del_z^s)&:=s&&&&{\text{order filtration}.}
  \end{align*}
Then plot for each $k$ the weights of $c_{k,n_k}\del_z^k$ in the $(F,V)$-plane:

\begin{centering}
  \begin{figure}[H]
    \caption{Two Fuchs polygons}
\psfrag{F}{$F$}
\psfrag{V}{$V$}
\psfrag{P1}{$P=z\del_z+1$}
\psfrag{P2}{$P=z^3\del_z+2$}
\includegraphics[width=0.2\textwidth]{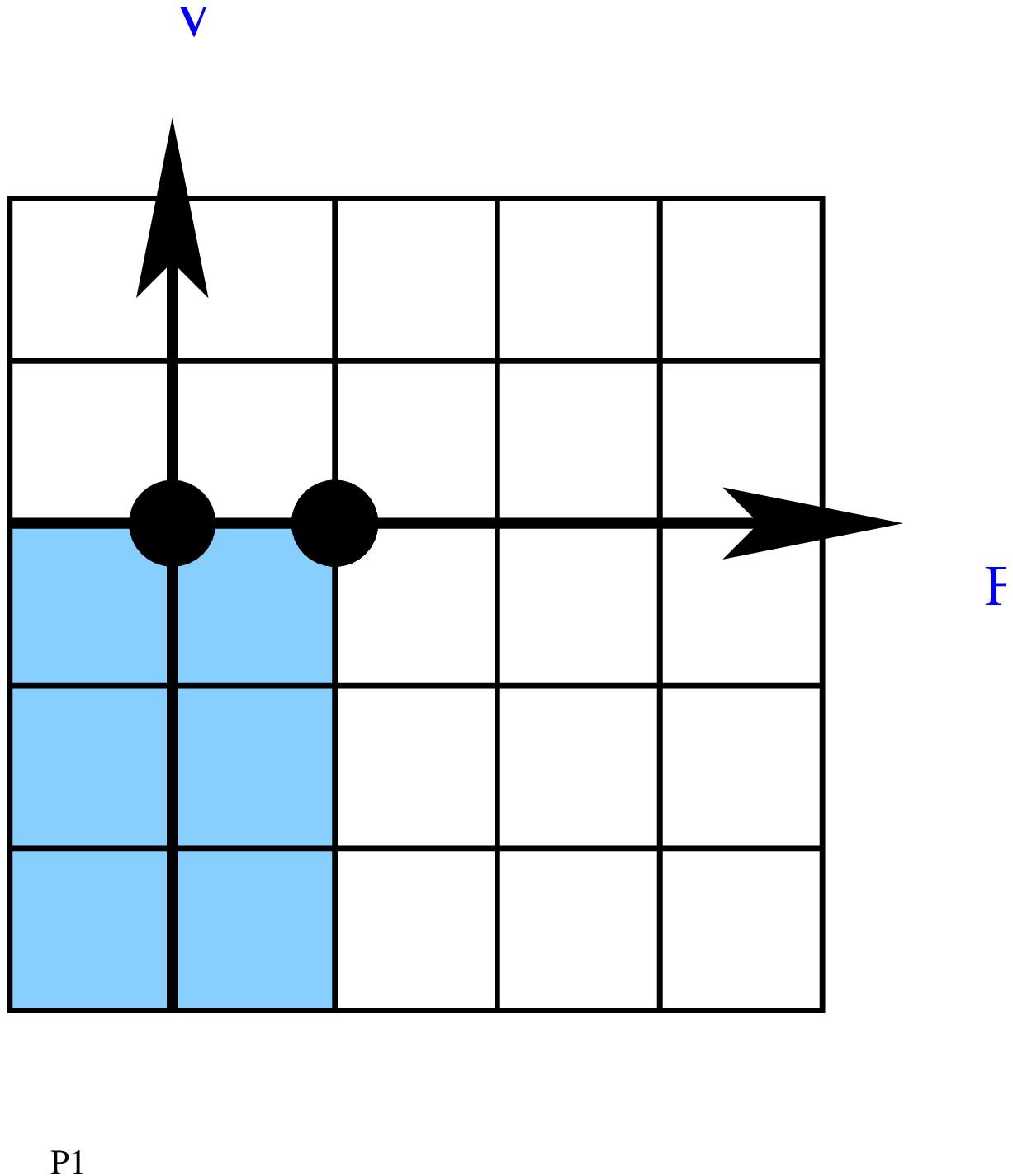}
\hskip0.2\textwidth
\includegraphics[width=0.2\textwidth]{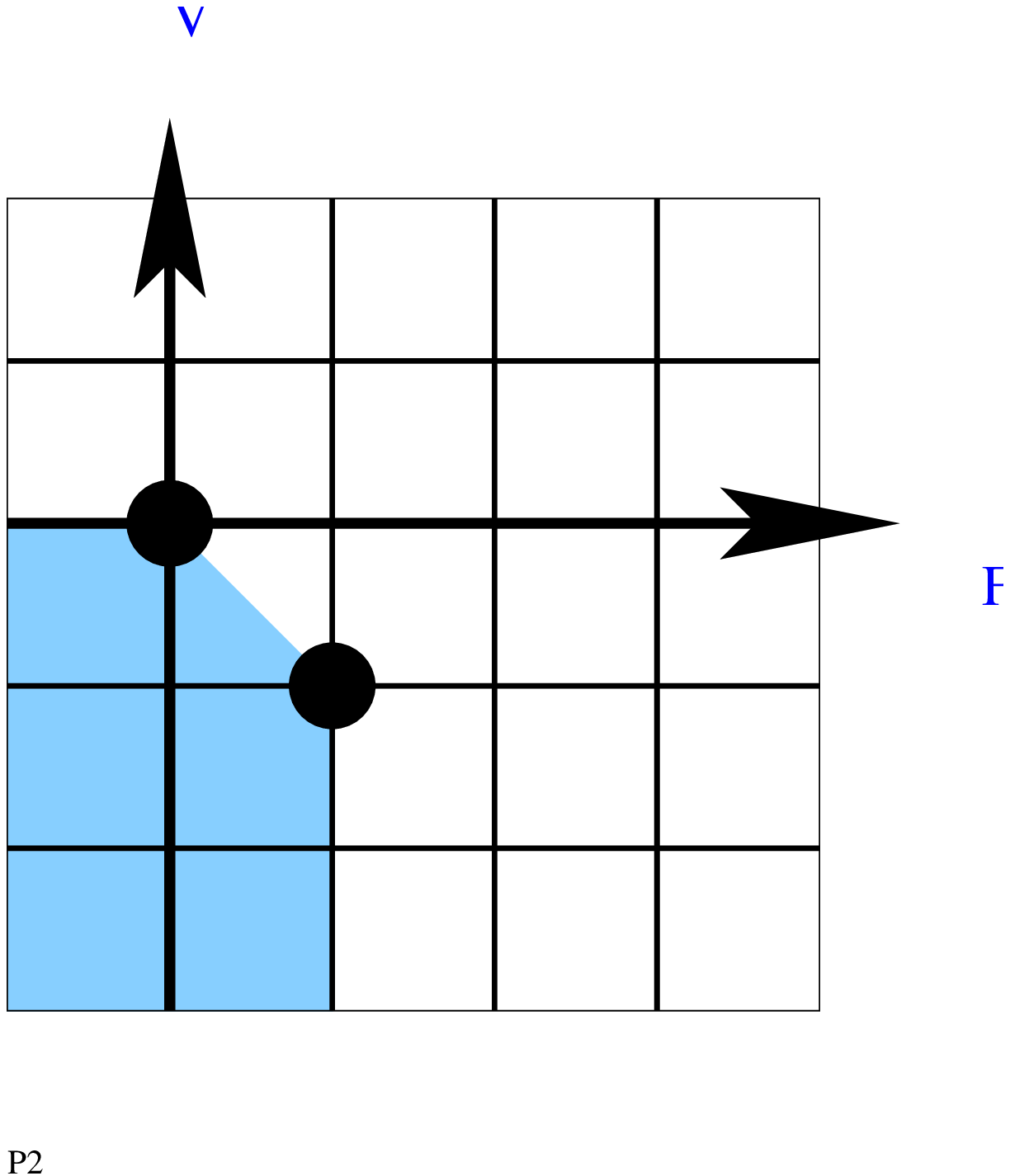}
\end{figure}
\end{centering}
The shaded region (the \emph{Fuchs polygon} of the operator) is the
lower left convex hull of the (finitely many) points so obtained. It
is, by definition, stable under shifts in negative $F$- and
$V$-direction, and hence
unchanged under analytic automorphisms that keep the origin fixed
(this is a consequence of taking the lower left hull).

Two cases arise, indicated in the picture:
\begin{enumerate}
  \item The  Fuchs polygon has one vertex, in the upper right
    corner (left).
    \item There are two or more corners. This is tantamount to the
      boundary of the shaded region having one or more finite boundary
      segments
      with \emph{slopes} different from $0$ and $-\infty$ (right).
\end{enumerate}
Fuchs' criterion (see \cite{Gray-BAMS,Ince44} for a detailed account)
states that $P$ has a regular singularity at  the origin if and only
if the Fuchs polygon of $P$ has no slopes.

Regular differential equations are much better behaved than irregular
ones, both theoretically and practically. On the theoretic side, they
form an ingredient of the Riemann--Hilbert correspondence that links
regular holonomic $D$-modules to perverse sheaves,
which for irreducible modules restricts to a bijection with
intersection cohomology complexes; on the practical side regular
differential equations are amenable to the Frobenius method since
their solutions come from the Nilsson ring
\cite{Kashiwara-RIMS84,Mebkhout-LNP80,Mebkhout-Compositio84,SST00}.

In higher dimensions, the concept of regularity is more difficult. One
way of defining it proceeds via pullbacks: the $\calD$-module $\calM$
on the analytic space $\CC^n$ is regular if and only if the pullback
of $\calM$ along any analytic morphism $\iota\colon \Delta^*\to \CC^n$, where
$\Delta^*$ is a punctured disk, leads to a module with regular
singularities at the origin on $\Delta^*$. The problem is that there
are many such morphisms to be tested.

Laurent \cite{Laurent87} and later with Mebkhout
\cite{LaurentMebkhout} found a way to translate regularity in more
than one variable into a condition that resembles the Fuchs
criterion. For that, we need to discuss filtrations and initial ideals
on $D$-modules in more detail.

\subsection{Initial ideals and triangulations}\label{subsec-triang}

A general technique to understand (non-commutative) algebraic structures is the reduction to a simpler (commutative) situation by applying a grading with respect to a filtration.
For $D$-modules, the filtration by the order of differential operators
leads to the characteristic variety which carries various bits of
information on the $D$-module.  The process of grading is rather
cumbersome but can be performed algorithmically in various situations
using Gr\"obner basis methods.  The simplest case is that of a generic
weight vector because the resulting graded ideal will be monomial;
this invites the use of techniques developed in \cite{SST00} and
\cite{Sturmfels96}.

So, let $L=(L_1,\dots,L_n)\in\QQ^n$ be a generic weight vector on
$R_A$; genericity is needed to assure that  $\gr^L(I_A)$ is a
monomial ideal.  (In $\RR^n$ there are weights $L$ that are generic
for all ideals of $R_A$ simultaneously. There is no rational weight
with this property, but for a finite number of
ideals a Zariski open set of the rational weight space consists of  generic weights.)

\begin{exa}\label{exa-Kummer-slopes}
  For the matrix $A=\begin{pmatrix}1&0&1\\0&1&1\end{pmatrix}$, with
  columns indicated with solid bullets, the following picture sketches the
  possible initial ideals that arise from the weights in the family
  $L^t=\begin{pmatrix}1&1&t\end{pmatrix}$, $t>0$. Note that
  $\bolda_1=\bolda_1/L^t_1$ and $\bolda_2=\bolda_2/L^t_2$ for all
  $t$.
  Plotted with
  hollow bullets  are the points $\bolda_3/L_3^t$ for the indicated
  choices of $t$.

  \noindent\begin{minipage}[h]{0.5\textwidth}
    \begin{figure}[H]
                                                \psfrag{2}{{\scriptsize $t=6/2$}}
      \psfrag{3}{{\scriptsize $t=6/3$}}
      \psfrag{4}{{\scriptsize $t=6/4$}}
      \psfrag{5}{{\scriptsize $t=6/5$}}
      \psfrag{6}{{\scriptsize $t=6/6$}}
      \psfrag{7}{{\scriptsize $t=6/7$}}
      \includegraphics[height=3cm]{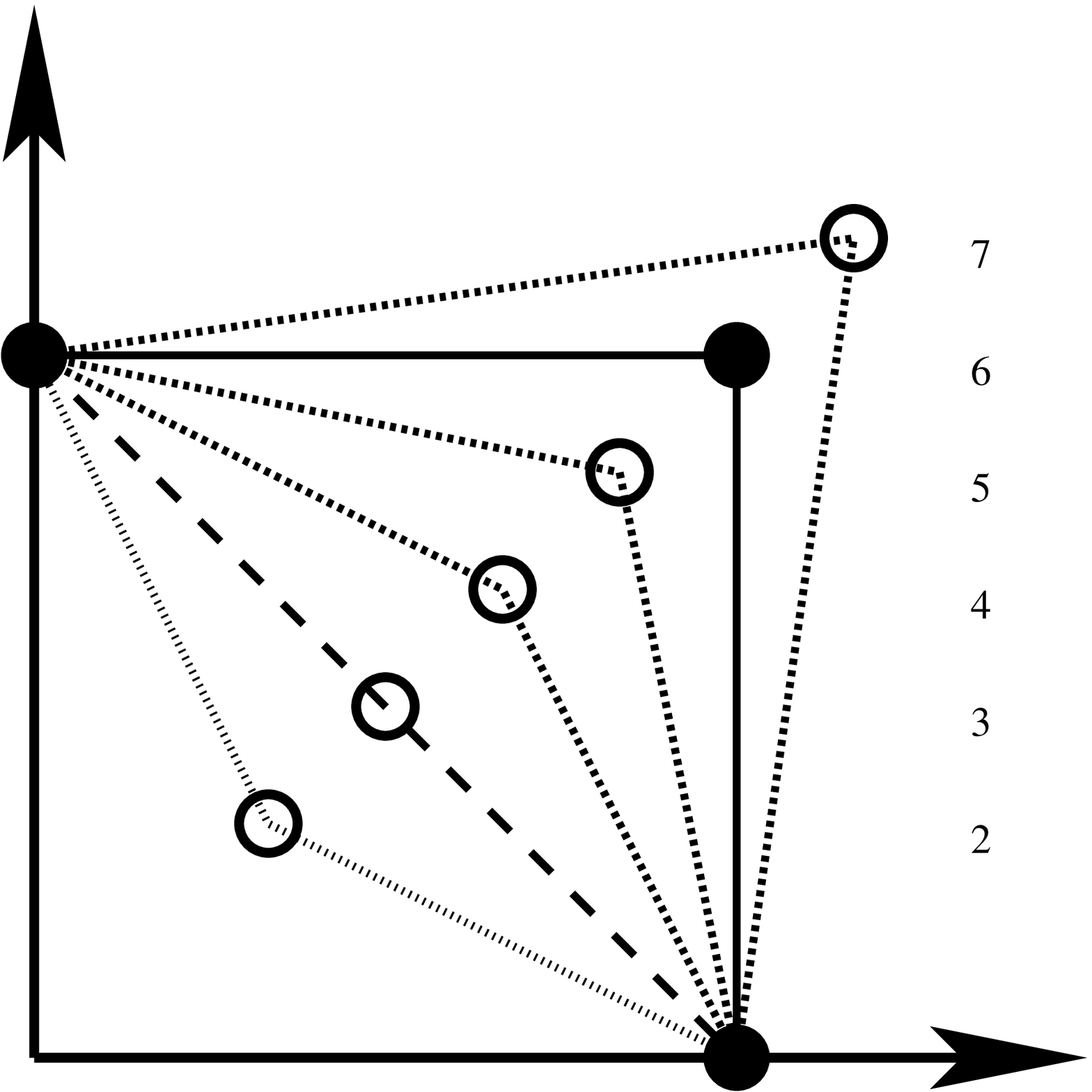}
    \end{figure}
  \end{minipage}
    \begin{minipage}[h]{0.5\textwidth}
     \begin{eqnarray*}I_A=\ideal{\del_1\del_2-\del_3}=\gr^{L^2}(I_A)\\\\
   \gr^{L^t}(I_A)=\left\{\begin{array}{ccc}\ideal{\del_1\del_2} & \text{ if }
   & 0<t<2;\\\ideal{\del_3} & \text{ if } & t>2.\end{array}\right.
   \end{eqnarray*}
\end{minipage}
Collinearity of $\{\bolda_1/L^t_1,\bolda_2/L^t_2,\bolda_3/L^t_3\}$
is equivalent to $L^t$-homogeneity of $I_A$.  \schluss\end{exa}

\begin{dfn}\label{dfn-Sigma}
  Associated to the generic weight $L$ and the $R_A$-ideal $I$ is an
  \emph{initial simplicial complex} $\Sigma^L_I$ that arises as
  follows. A collection $\tau$ of indices contained in $[n]$ forms a
  face of $\Sigma^L_I$ if and only if there is no monomial in
  $\gr^L(I)$ whose support is precisely $\tau$. Put another way,
  $\Sigma^L_I$ is the simplicial complex whose Stanley--Reisner ideal
  is the radical of $\gr^L(I)$.

If $I=I_A$ we write $\Sigma^L_A$ for $\Sigma^L_{I_A}$.
\schluss\end{dfn}

For example, suppose $I_A$ is the principal ideal
generated by $\del_1\del_2\del_3-\del_4\del_5^2$. Then $I_A$ admits
two distinct monomial initial ideals whose corresponding simplicial
complexes are:

\begin{figure}[H]
\caption{The initial simplicial complexes $\Sigma^L_A$  for $I_A=(\del_1\del_2\del_3-\del_4\del_5^2)$.}
\centering
\hfill
\subfigure[The join of a line segment with a 3-cycle, $\gr^L(I_A)=\del_1\del_2\del_3$.]{
 \psfrag{1}{$1$}
 \psfrag{2}{$2$}
 \psfrag{3}{$3$}
 \psfrag{4}{$4$}
 \psfrag{5}{$5$}
 \qquad\includegraphics[width=0.25\textwidth]{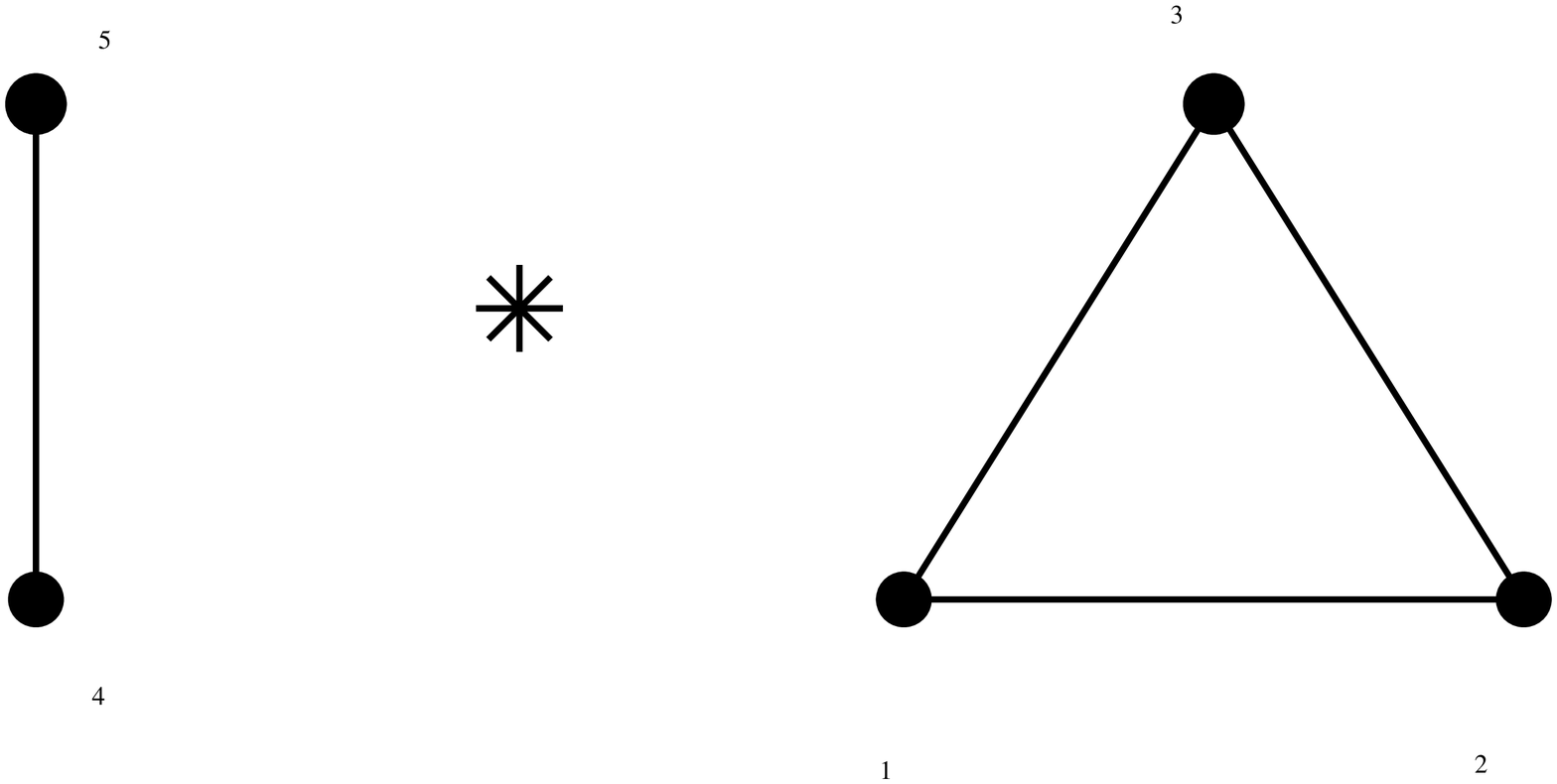}\qquad
}
\hfill\hfill
\subfigure[The join of two points with a triangle, $\gr^L(I_A)=\del_4\del_5$.]{
 \psfrag{1}{$1$}
 \psfrag{2}{$2$}
 \psfrag{3}{$3$}
 \psfrag{4}{$4$}
 \psfrag{5}{$5$}
 \qquad\includegraphics[width=0.25\textwidth]{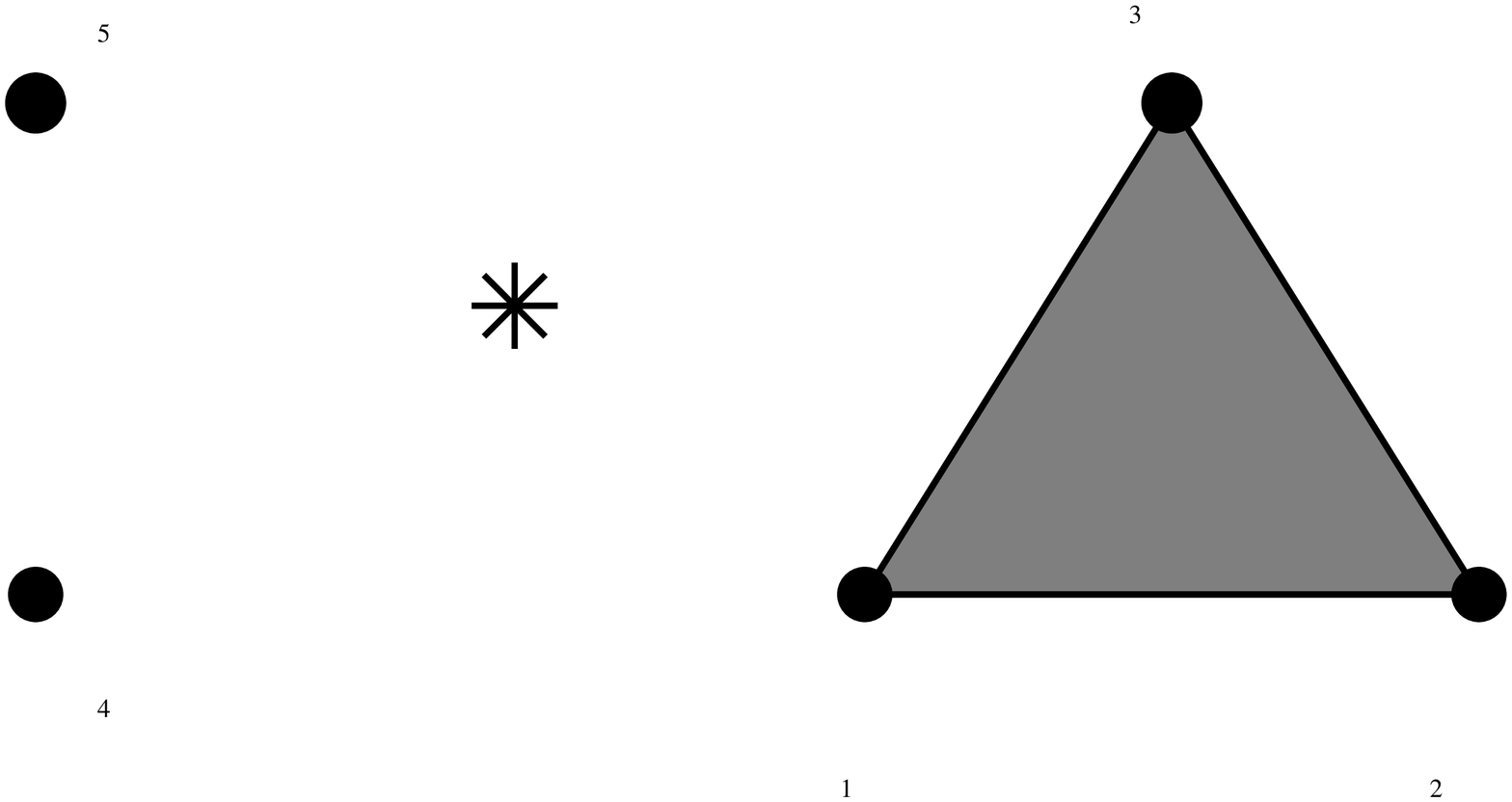}\qquad
}
\hfill
\label{fig-two-complexes}
\end{figure}

The generic weight $L$ also induces a triangulation of $[n]$ as
follows. Consider the points $\hat A=\{(\bolda_j,L_j)\in\RR^d\times
\RR\}_{1\le j\le n}$. The faces of the triangulation are those faces
of the cone $\RR_{\geq 0}\hat A$ of $\hat A$ that are visible from the
point $(\boldzero,-\infty)$; these are exactly those faces whose outer
normal vectors have negative last component. A triangulation of $[n]$
is \emph{regular} (or \emph{coherent}) if it arises this way for some
$L$. This property is strongly tied to $A$, and not all triangulations
of $A$ have to be regular.
\begin{figure}[H]
\caption{A non-regular triangulation of a triangle.}
\centering
\includegraphics[width=0.15\textwidth]{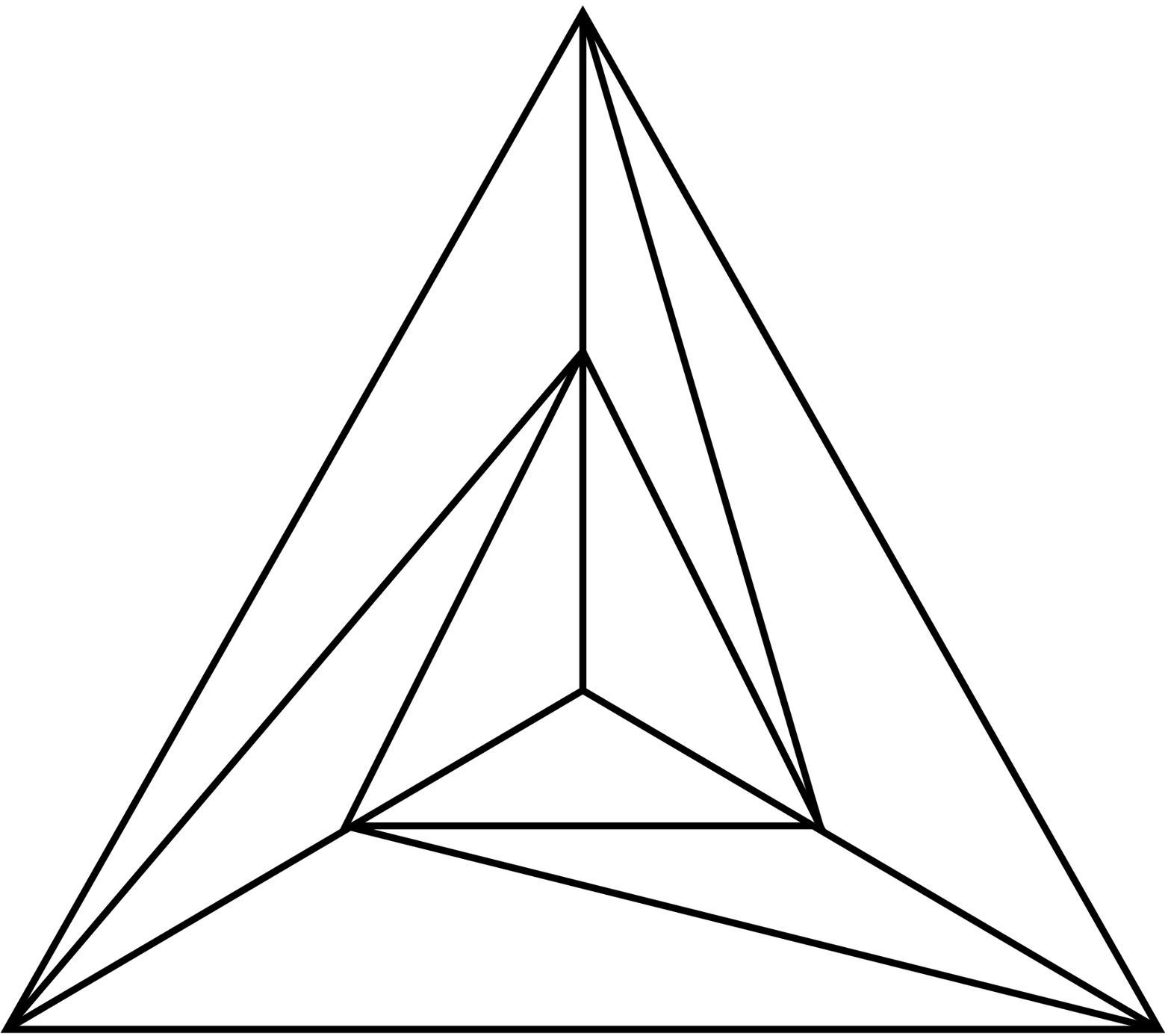}
\label{fig-nonreg}
\end{figure}

The collection of regular triangulations of $A$ turns out to be in (the
obvious) bijection with the initial complexes of $A$.
There is a third combinatorial object associated to $L$ and $A$,
namely the collection $\calS(\gr^L(I_A))$ of \emph{standard pairs} of
$\gr^L(I_A)$, introduced in \cite{SturmfelsTrungVogel-MathAnn95}.
A standard pair $(\bolddel^\boldb,\sigma)$ of the monomial ideal $\bar
I$
is a monomial and a subset of $[n]$ such that
\begin{itemize}
\item $\supp(\boldb)\cap\sigma=\emptyset$,
\item $\bolddel^\boldb\mod \bar I$ is not $(\prod_{j\in\sigma}\del_j)$-torsion, but
\item $\bolddel^\boldb\mod \bar I$ is $\del_k(\prod_{j\in\sigma}\del_j)$-torsion for all $k\not\in\sigma$.
\end{itemize}
For example, if the monomial ideal is $(\del_4\del_5^2)$ the standard
pairs are $(1,\{1,2,3,4\})$, $(\del_5,\{1,2,3,4\})$, and
$(1,\{1,2,3,5\})$.  The standard pairs yield immediately a
decomposition into irreducible ideals by
\[
\bar I=\bigcap_{(\bolddel^\boldb,\sigma)\in \calS(\bar I)}(\{\del_j^{b_j+1}\mid
j\not\in\sigma\}).
\]
For $\bar I$ as above we obtain $\bar I=(\del_5)\cap (\del_5^2)\cap (\del_4^1)$.

The standard pairs hence contain all information needed to recover
$\bar I$
and its triangulations. In particular, the facets of $\Sigma^L_A$ are
precisely the subsets $\sigma$ that are listed in the standard pairs.

\begin{exa}\label{exa-Kummer-more-slopes}
  We consider Example \ref{exa-Kummer-slopes} from this new angle.
  We fix the weights $L_1=L_2=1$ and
vary the weight $t=L_3$.  For $L_3<2$, $\gr^LI_A=\ideal{\del_1\del_2}$
and the facets of $\Sigma_A^L$ are
$\{1,3\},\{2,3\}$.  We could interpret this as the complex of faces,
not containing $\boldzero$,
of the convex hull of $\boldzero$ and the columns of $A$.
Similarly we obtain $\Sigma_A^L=\{1,2\}$ for $L_3>2$, which can be
read as a convex hull as before, but with $\bolda_3$ not in the
picture.
For $L_3=2$, $\gr^LI_A=I_A$ is prime and
$\Sigma_A^L$ should now equal $\{1,2,3\}$: we would like to view $\bolda_3$ as
``collinear with $\bolda_1,\bolda_2$'' in this case.
This is the topic
of the next section; the following is a teaser:
in order to view the three cases from a unifying angle, note that
scaling a weight component $L_i$ by $\lambda$ and ``scaling the degree
$\bolda_i$ of $\del_i$'' by $1/\lambda$ have the same effect on the
initial terms (and also on the face complex of $\Sigma^L_A$).
One is thus led to replace $\bolda_3$ by $\bolda_3/L_3$;
then  the
resulting convex hull yields the face complex generated by
$\{1,2,3\}$ if $L_3=2$, by $\{1,2\}$ is $L_3>2$, and by $\{1,3\}$ and
$\{2,3\}$ if $L_3<2$.
\schluss\end{exa}

\subsection{Slopes and the $(A,L)$-umbrella}

In case of a $D_A$-module $M=D_A/J$, $J$ an ideal in $D_A$, we will want
to grade with respect to a filtration on $D_A$ defined by (and
identified with) a weight vector $L\in\QQ^d\times\QQ^d$ for the
variables $x_1,\dots,x_n,\del_1,\dots,\del_n$.  We denote the
$L$-leading term of $P\in D_A$ by $\sigma^L(P)$ and call it the
\emph{$L$-symbol}.
\begin{cnv}
  We assume that there is a positive real constant $c$ such that
  \[
  L_{x_j}+L_{\del_j}=c>0
  \]
  for all $j$ simultaneously.
\schluss\end{cnv}

This  hypothesis has the effect that
\[
W_A:=\gr^L(D_A)\cong \CC[\boldx,\bolddel]
\]
is a (commutative) polynomial ring whose spectrum is naturally
identified with the total space of the cotangent bundle $T^* \CC^n$ of
$\CC^n$. Moreover, each $E_i$ is $L$-homogeneous of positive degree.

The $W_A$-ideal $\gr^L(J)$ defines the \emph
{$L$-characteristic variety} $\ch^L(M)$ of the module $M$; for a
holonomic module $M$ it is purely $n$-dimensional by a result of
G.G.~Smith \cite{Smith01}.

We record
the special case
\[
\ch^L(M_A(\beta))=\Var(\gr^L(H_A(\beta)))\subseteq T^* \CC^n
\]
when $M=M_A(\beta)$. Our plan is to connect this construction to
analytic information as follows.

Suppose $X'\subseteq X=\CC^{n,\an}$ is an analytic subspace with a smooth
point $\pointx\in X'$. Then in suitable local coordinates at $\pointx$ one can
write $X'$ as the zero set of the first $n-\dim X'$ coordinates on
$X$. In the stalk at $\pointx$ consider the grading of the $D$-module $M$
by the filtrations induced by the weights $L^{p/q}:=pF+qV$ where as always $F$
is the order filtration and $V$ is the $V$-filtration along $X'$
(compare Subsection \ref{subsec-Fuchs}):
\[
V(x_i)=V(\del_i)=0\text{ if } i>n-\dim(X');\qquad -V(x_i)=V(\del_i)=1\text{ if
}i\le n-\dim(X').
\]
(There is an obvious identification of graded objects for $L^{p/q}$
and $L^{p'/q'}$ when $p/q=p'/q'$).
\begin{dfn}
  With notation as just introduced,
  $p/q\in\QQ$ is a \emph{slope of $M$ along $X'$} if
$\ch^L(M)=\supp(\gr^L(M))$ jumps at $p/q$. This means that
$\ch^{L^\eps}(M)$ is for small $\eps\in\RR_+$  constant on $(-\eps+\frac{p}{q},\frac{p}{q})$
  and $(\frac{p}{q},\frac{p}{q}+\eps)$ but not on $(-\eps+\frac{p}{q},\frac{p}{q}+\eps)$.
\schluss\end{dfn}
This definition is taken from \cite{Laurent87}. By
\cite{LaurentMebkhout}, Laurent's algebraic slopes constructed from
filtrations agree with
Mebkhout's transcendental slopes given as jumps of the Gevrey
filtration on the irregularity sheaf and hence provide a measure of
growth for the solutions of $M$. The central question in this section
is to study the behavior of $\ch^L(M_A(\beta))$ under changes of $L$
and $\beta$.

\medskip

We illustrate the link of slopes of $M_A(\beta)$
with Fuchs' criterion in an example.

\begin{exa}
  It is clear from the series expansion \eqref{eqn-pFq} that
  the Kummer confluent series ${}_1F_1(a;b;z)$ is analytic at every
  finite $z$ for all $a,b$. On the other hand, it follows from the
  integral definition of the error function that at $z=\infty$ there
  is an essential singularity (and algebraic changes of coordinates do
  not eradicate essential singularities). If we denote $-1/z$ by $u$, then the
  differential operator $\theta_z(\theta_z+1/2)-z(\theta_z-1/2)$ turns
  into $u\theta_u(\theta_u-1/2)-(\theta_u+1/2)$ for the resulting
  inverse Kummer confluent series.

  The Fuchs polygons are:
  \[
  \begin{minipage}{\textwidth}
    \begin{figure}[H]
      \caption{Fuchs polygon for Kummer (left) and inverse
        Kummer (right)}
      \psfrag{V}{$V$}
      \psfrag{F}{$F$}
      \includegraphics[width=0.2\textwidth]{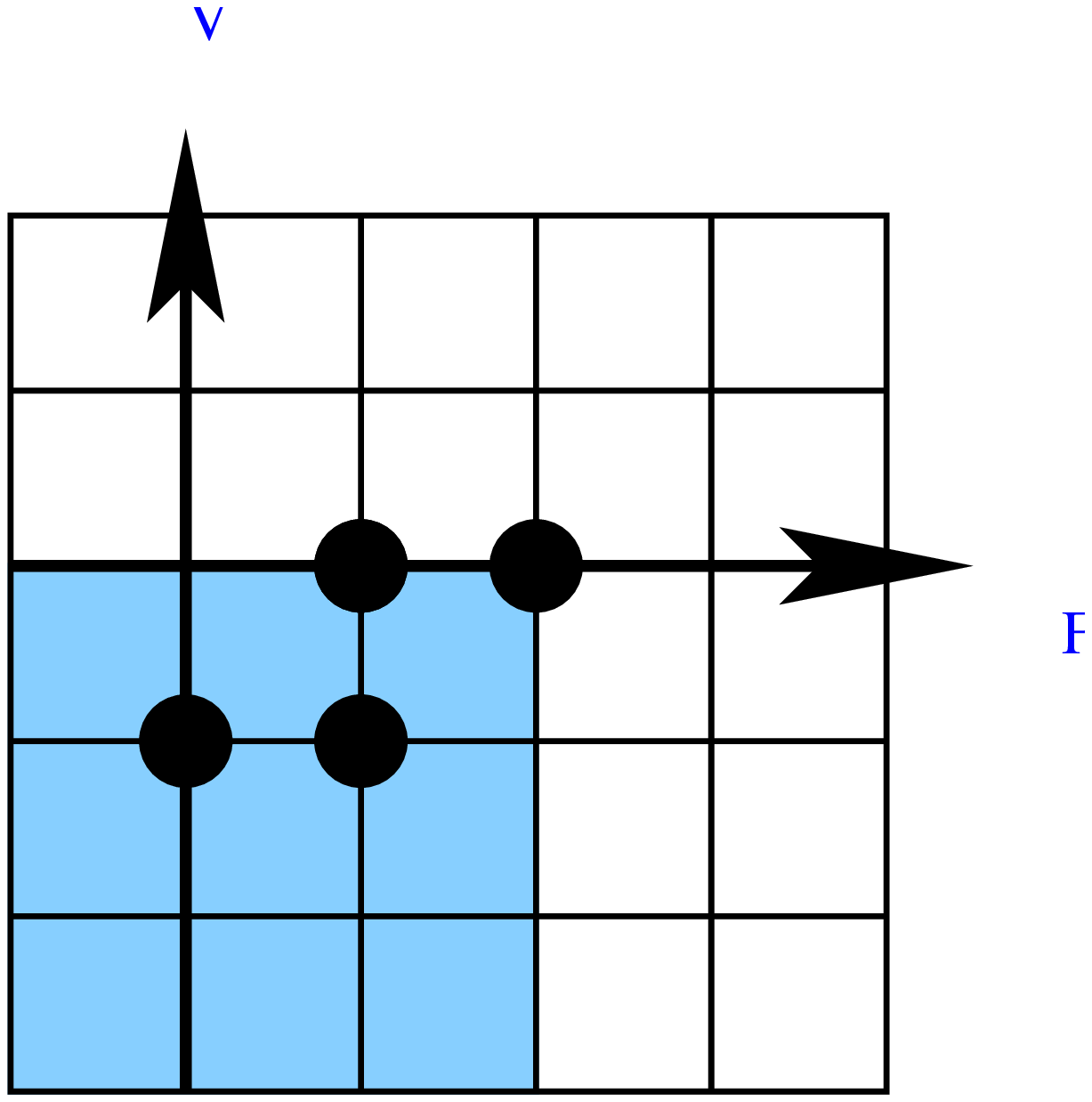}
      \hspace{0.2\textwidth}
      \psfrag{V}{$V$}
      \psfrag{F}{$F$}
      \includegraphics[width=0.2\textwidth]{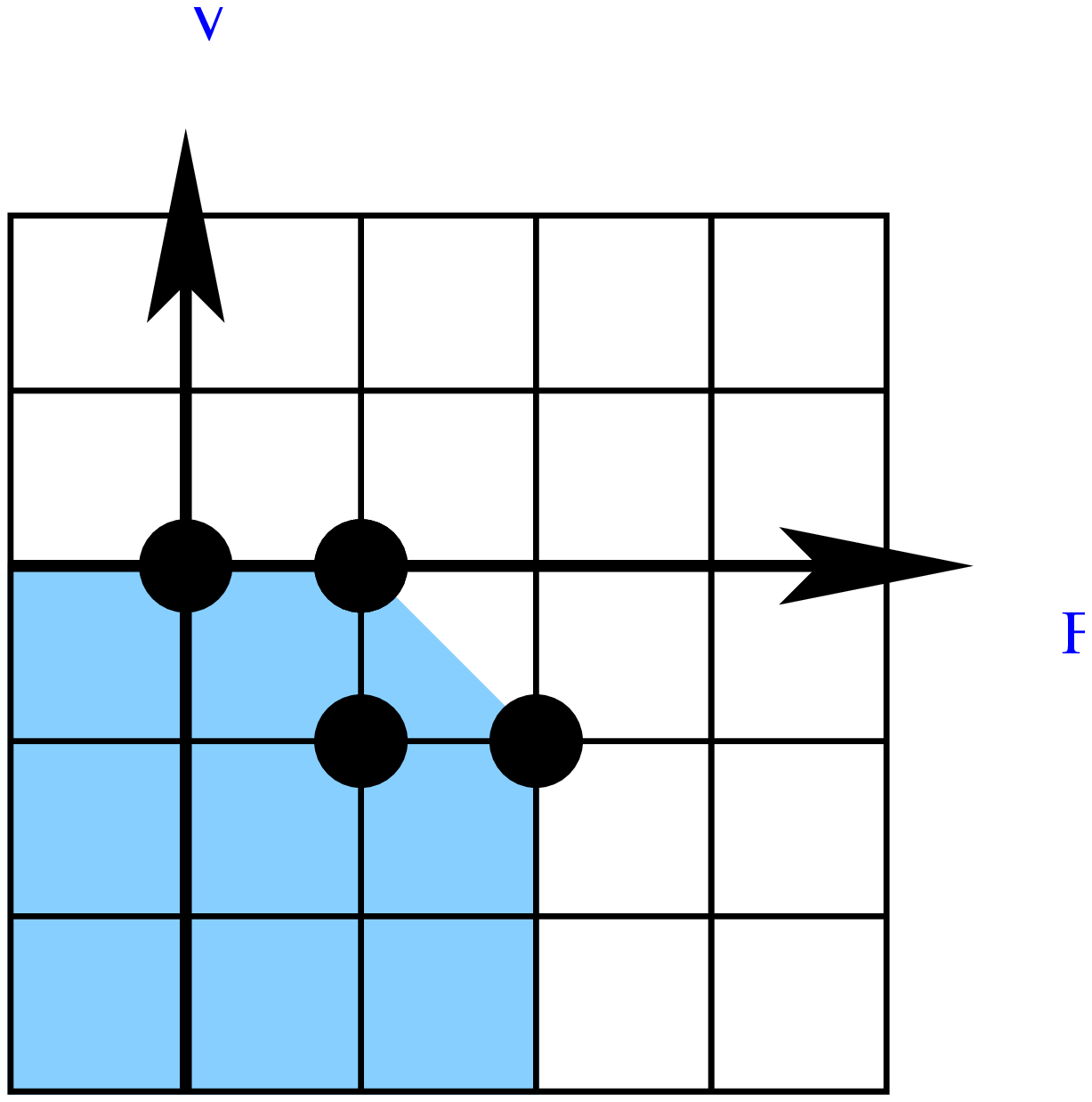}
    \end{figure}
  \end{minipage}
  \]

  So, the Kummer series has (of course) regular ``singularities'' at
  the origin, while the inverse Kummer series has a slope of
  $-1$. This reflects the fact that, up to multiplication by a
  function bounded by a polynomial, the Kummer series at $0$ behaves
  like $\exp(z^0)$, while the inverse Kummer series behaves like
  $\exp(z^{-1})$: the Kummer series grows (up to polynomially bounded
  factors) near $\infty$ like $\exp(z)$.

For the translation to the $A$-hypergeometric setting we can use in
both cases $A=\begin{pmatrix}1&0&1\\0&1&1\end{pmatrix}$, with $\boldv$
being $(1,1,-1)$ or $(-1,-1,1)$.  The toric ideal is then
$I_A=\ideal{\del_1\del_2-\del_3}$.

We know from Example \ref{exa-Kummer-slopes} that for the family
$L^t=(1,1,t)$ there is a jump at $t=2$ in the $L^t$-graded ideal of
$I_A$ since at that moment $\Box_\boldv$ becomes $L$-homogeneous. It
turns out that the $L^t$-characteristic variety of $H_A(\beta)$ for
any $\beta$ also changes at $t=2$, so that $M_A(\beta)$ has a slope of
$2$ along the hyperplane $x_3=0$.

The correspondence between these numbers is
encapsulated by the equation $\frac{1}{s_F}=\frac{1/s_L}{1/s_L-1}$, where
$s_F$ is the slope of the Fuchs polygon (and indicates exponential growth
behavior with exponent $s_F$), and $s_L$ is the slope at which
Laurent's filtrations jump.
\schluss\end{exa}

\medskip

We now discuss  ``regular triangulations to  non-monomial
graded toric ideals'' coming from non-generic weight vectors in
greater generality,  the details
being taken from \cite{SchulzeWalther-duke}.
For the transition, suppose $J$ is generated by
elements inside $R_A\subseteq D_A$. Then one can restrict the weight
to $L_\bolddel$ on $R_A$ and compute $\gr^{L_\bolddel}(J\cap R_A)$ in
the commutative situation of Subsection \ref{subsec-triang}. Note that then
$\gr^L(J)=\gr^L(D_A)\cdot \gr^{L_\bolddel}(J\cap R_A)$.
Specifically, we write
\[
I_A^L:=\gr^L(I_A)\cap R_A,\quad S_A^L:=\gr^L(S_A)\cong R_A/I^L_A.
\]

Let $L=(L_1,\dots,L_n)\in\QQ^n$ be any weight vector on $R_A$. As
$L$ may have zero components, possible division (as suggested in
Example \ref{exa-Kummer-more-slopes}) by $L_i=0$ forces us
into work in a projective space:
\[
\bolda_1,\dots,\bolda_d\in\ZZ A\subseteq \QQ^d\subseteq\PP^d_\QQ.
\]
In $\PP^d_\QQ$, any two distinct points $\bolda,\boldb\in\PP^d_\QQ$
are joined by two line segments.  If the hyperplane $H$ in $\PP^d_\QQ$
contains neither $\bolda$ nor $\boldb$, one may define the
convex hull of $\bolda,\boldb$ as the line segment not
intersecting $H$.  Similarly one can define the convex hull
$\conv_H(S)$ of a subset $S\subseteq\PP^d_\QQ$ disjoint from $H$ as
the \emph{convex hull} of $S$ in the affine space $\PP^d_\QQ\smallsetminus
H$.

\begin{dfn}[The $(A,L)$-umbrella $\Phi^L_A$]\label{96}
We set $\bolda_j^L:=\bolda_j/L_j\in\PP^d_\QQ$.
Choose a linear functional $f\colon \ZZ A\to\ZZ$ for which $f(\bolda_j)>0$
for all $j$ and $\eps>0$ such that $|f(\bolda_j)|>\eps\cdot|L_j|$;
such form exists since $A$ is pointed.
Let $H_\eps:=f^{-1}(-\eps)$ and call
\[
\Delta^L_A:=\conv_{H_\eps}(\{\boldzero,\bolda_1^L,\dots,\bolda_n^L\})\subseteq\PP^d_\QQ
\]
the {\em $(A,L)$-polyhedron}.
Let the {\em $(A,L)$-umbrella} be the set $\Phi_A^L$ of faces of
$\Delta^L_A$ which do not contain $\boldzero$; write
$
\Phi^{L,k}_A
$
for its $k$-skeleton.

The matrix $A$ is called {\em
  $L$-homogeneous} if all $\bolda^L_j$ lie on a common hyperplane of
$\PP_\QQ^d$.  Every $A$ is $\boldzero$-homogeneous and we call
$\Phi_A:=\Phi_A^0$\index{PA@$\Phi_A$} the \emph{$A$-umbrella}.  Note
that $\Phi_A$ can be identified with the face lattice of the
polyhedral cone $\RR_{\geq 0}A$.
\schluss\end{dfn}

Parts of this definition, taken from \cite{SchulzeWalther-duke} are
foreshadowed by \cite[Prop.~4]{GKZ89}.

\begin{figure}[H]
\caption{$(A,L)$-umbrellas for Example \ref{exa-4slopes}. (Blue $\Delta^L_A$ with  boundary $\Phi^L_A$.)}\label{70}

       \psfrag{a1}{$\bolda_1$}
       \psfrag{a2}{$\bolda_2$}
       \psfrag{a3}{$\bolda_3$}
       \psfrag{a4}{$\bolda_4$}
       \psfrag{a4L}{$\bolda_4^L$}
       \psfrag{0}{$0$}
       \psfrag{1}{$1$}
       \psfrag{2}{$2$}
       \psfrag{3}{$3$}
       \psfrag{Ra4}{$\RR\bolda_4$}

       \includegraphics[width=0.24\textwidth]{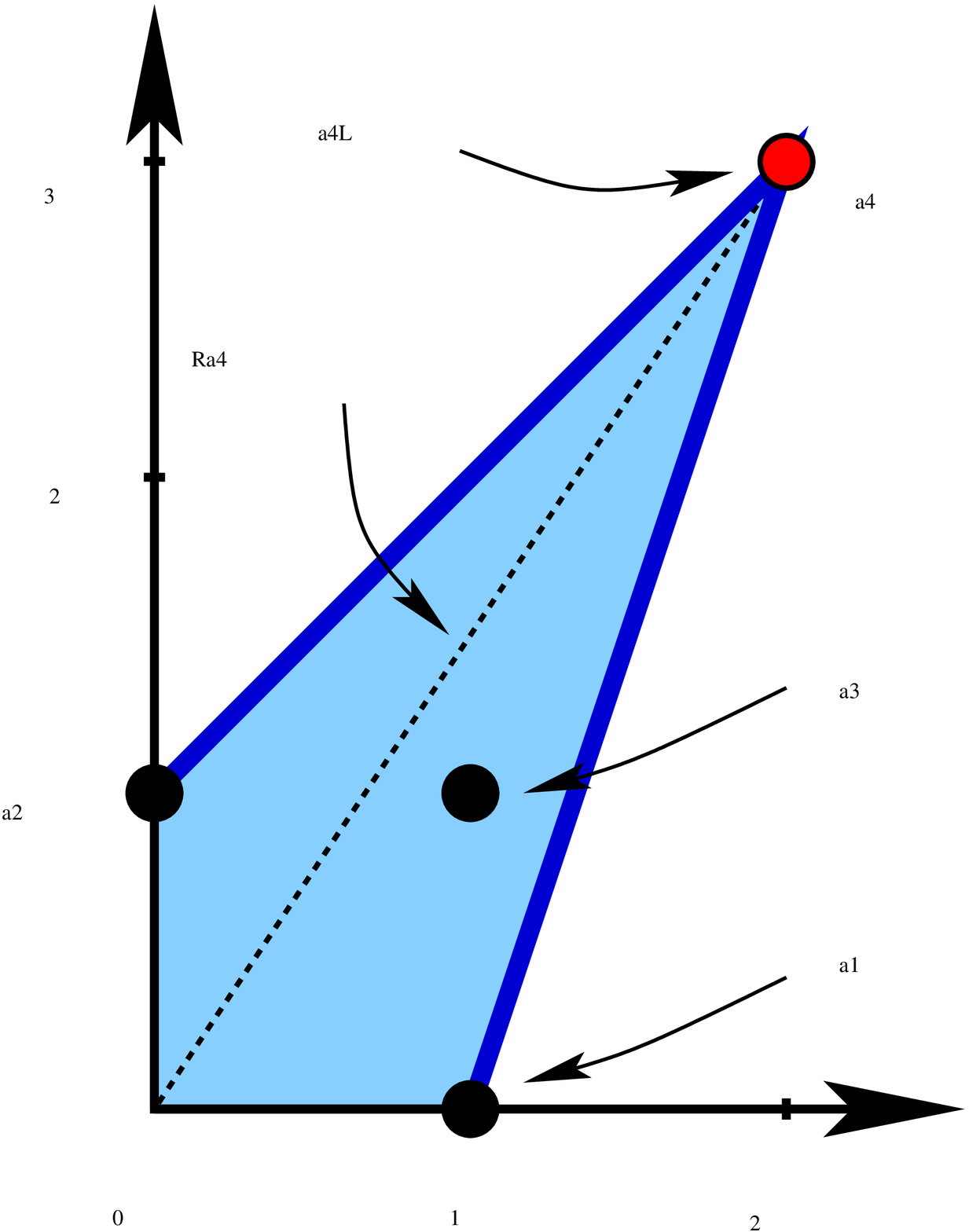}
       \includegraphics[width=0.24\textwidth]{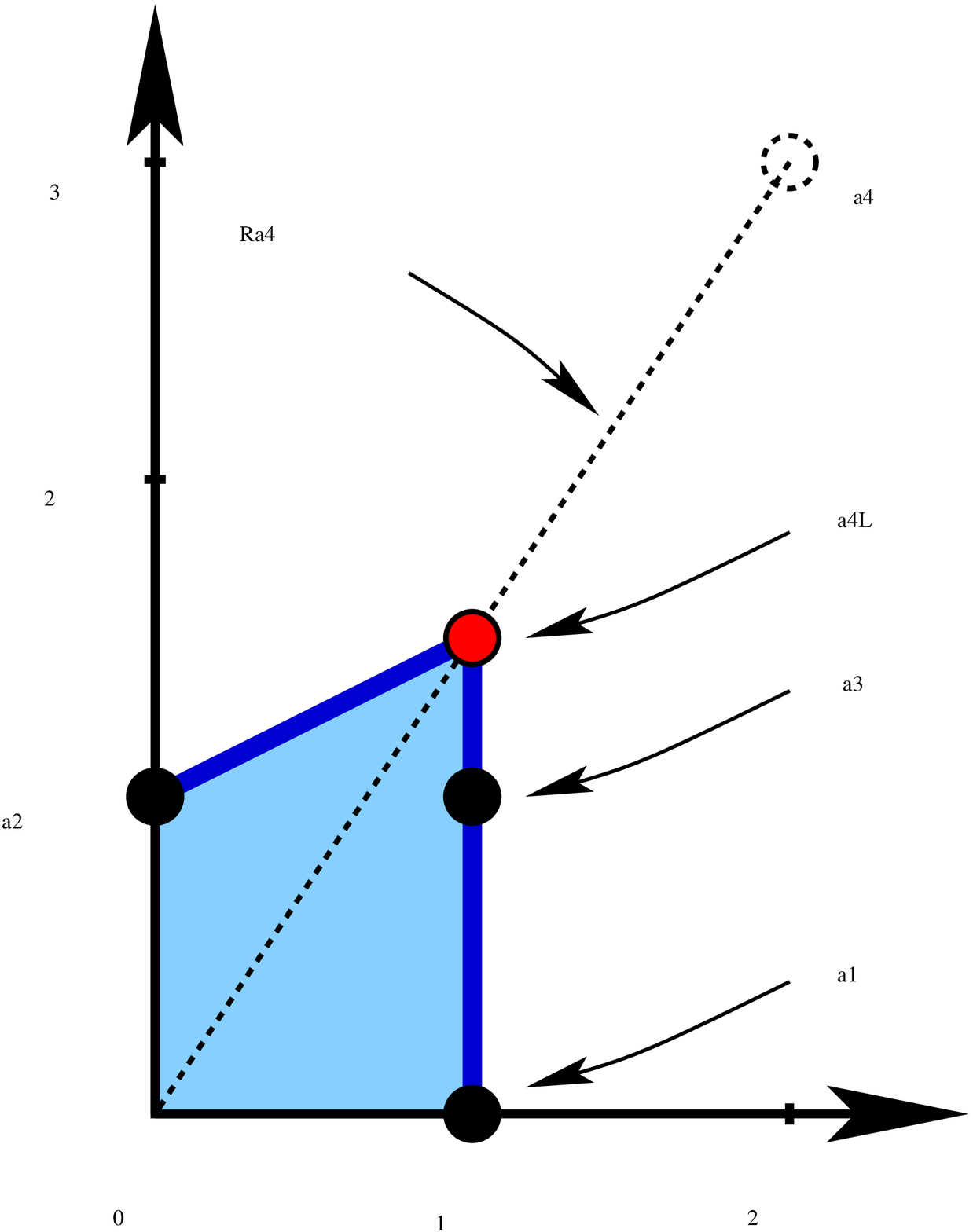}
       \includegraphics[width=0.24\textwidth]{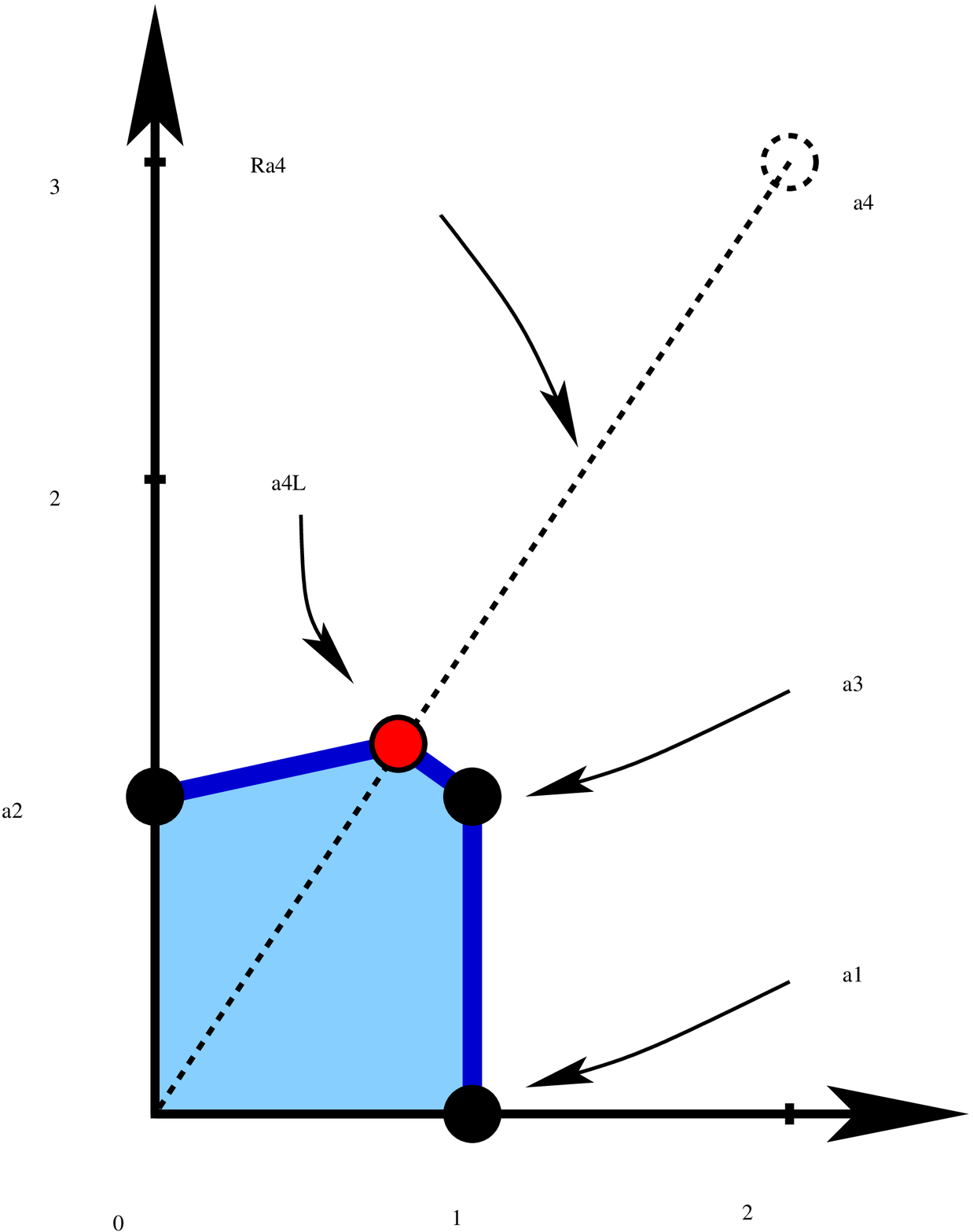}
       \includegraphics[width=0.24\textwidth]{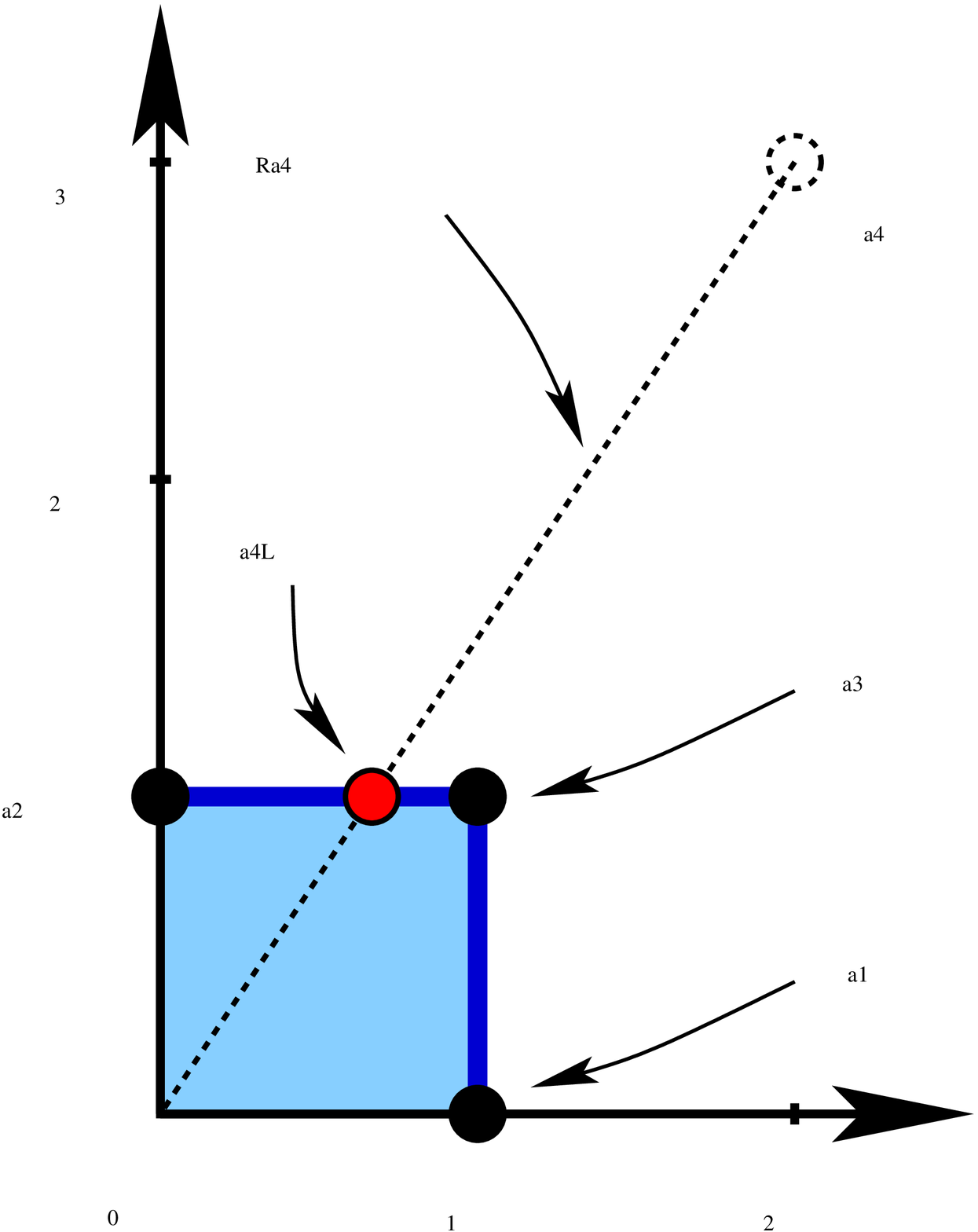}
\end{figure}

\begin{exa}\label{exa-4slopes}
Figure \ref{70} shows the $(A,L)$-umbrella for the matrix
$A=\begin{pmatrix}1&0&1&2\\0&1&1&3\end{pmatrix}$ for various
filtrations in the family $L^t=(1,1,1,t)$.  While moving the parameter,
$\Phi_A^L$ jumps exactly at $t=2$ and $t=3$.  For the intervals $t<2$, $t=2$,
$2<t<3$, $t=3$, $t>3$, the corresponding complexes $\Phi_A^L$ are
generated by $\{\{1,4\},\{2,4\}\}$, $\{\{1,3,4\},\{2,4\}\}$,
$\{\{1,4\},\{2,4\},\{3,4\}\}$, $\{\{1,3\},\{2,3,4\}\}$,
$\{\{1,3\},\{2,3\}\}$.
\schluss\end{exa}

\begin{rmk}
  In order to see how $\Phi_A^L$ generalizes $\Sigma_A^L$ for positive
weights, embed $\PP^d_\QQ\subseteq \PP^{d+1}_\QQ$
as the hyperplane $\{a_{d+1}=a_0\}$, and assume that $L$ is positive
and generic. A subset of
$\{\bolda_1^L,\ldots,\bolda_n^L\}\subseteq \AA^d_\QQ\subseteq
\PP^d_\QQ$ maximizes a linear functional $q(t_1/t_0,\ldots,t_d/t_0)$
with value $c$ if and only if the corresponding subcollection of
$\{(\bolda_j,L(\bolda_j)\}_1^n\subseteq \AA^{d+1}_\QQ \subseteq
\PP^{d+1}_\QQ$ maximizes \emph{with value zero} the linear functional
$q(t_1/t_0,\ldots,t_d/t_0) -t_{d+1}/t_0$. So, the faces of
$\Delta^L_A\times\{1\}\subseteq \AA^{d+1}_\QQ$ are in bijection with those of
the cone spanned by it from the origin in $\AA^{d+1}_\QQ$ that have
outer normal vector ``pointing down'', and this is the same cone as
the one spanned by the appropriate collection inside
$\{(\bolda_j,L(\bolda_j)\}_1^n$.
\schluss\end{rmk}

Just like $\Sigma_A^L$ in the monomial case, $\Phi_A^L$ corresponds to
minimal prime ideals of $\gr^L(I_A)$.  More precisely the following holds.

\begin{thm}[{\cite[Thm.~2.14]{SchulzeWalther-duke}}]\label{13}
The set of $A$-graded prime ideals containing $I^L_A$ equals
$\{I_A^\tau\mid\tau\in\Phi_A^L\}$ and so
\[
\Spec(S^L_A)=\Var(I_A^L)=\bigcup_{\tau\in\Phi_A^{L,d-1}}\bar O_A^\tau=\bigsqcup_{\tau\in\Phi_A^L}O_A^\tau
\subset \widehat{\CC}^n.
\]
In particular, the $(A,L)$-umbrella
encodes the geometry of $S^L_A$.\qed
\end{thm}

\subsection{$L$-characteristic varieties}\label{subsec-Lchar}

Equipped with the knowledge from the previous section, we can return
to the question of describing
\[
\Upsilon_A^L:=\ch^L(M_A(\beta)).
\]
For a weight $L\in\QQ^n\times\QQ^n$,
the $L$-symbols $\sigma^L(E_i)$ span the tangent spaces of every
torus orbit and hence impose the conormal condition to $O_A^\tau$ for
all $\tau\in\Phi_A^L$ (compare \cite{GKZ89, SchulzeWalther-duke}).
The inclusion
\begin{equation}\label{95}
\gr^L(H_A(\beta))\supseteq\ideal{\sigma^L(E)}+\gr^L(D_A\cdot I_A^L)
\end{equation} appears
already in \cite{GKZ89,Adolphson-duke94} and shows that $\ch^L(M_A(\beta))$
must be contained in the union of the closures of all these conormals.

One might hope that \eqref{95} is always an equality; this would
simplify the problem of describing $\ch^L(M_A(\beta))$.  The right
hand side is the \emph{fake initial ideal} and equality holds if
$I_A^L$ is Cohen--Macaulay, \cite[Thm.~4.3.8]{SST00}.  Unfortunately,
this inclusion can be strict in general as the following example shows.
\begin{exa}
For $A=\begin{pmatrix}1&1&1&1\\0&1&3&4\end{pmatrix}$ and
$L=(\boldzero,\boldone)$ inducing the order filtration one has
$\gr^L(H_A(\beta))=\gr^L(D_A\cdot I_A)+\ideal{\sigma^L(E)}$ for
$\beta=(1,2)$, but in fact for all parameters
\[
\gr^L(H_A(\beta))=\gr^L(D_A\cdot I_A)+\ideal{\sigma^L(E)}+\ideal{P}
\]
where
\[
P=
(\beta_2-2)x_1\del_1^2+
(\beta_2-\beta_1-1)x_2\del_1\del_3+
(\beta_2-3\beta_1+1)x_3\del_2\del_4+
(\beta_2-4\beta_1+2)x_4\del_3^2.
\]
\schluss\end{exa}

Notwithstanding this example, the following is true.

\begin{thm}\label{3}
The $L$-characteristic variety of the $A$-hypergeometric system is
\[
\Upsilon_A^L=\ch^L(M_A(\beta))=\bigcup_{\tau\in\Phi_A^L}\bar\Upsilon_A^\tau=\bigsqcup_{\tau\in\Phi_A^L}\Upsilon_A^\tau,
\]
where for $\tau\in\Phi_A^L$, we denote by $\Upsilon_A^\tau\subseteq T^* \widehat{\CC}^n$
the conormal to the orbit $O_A^\tau\subseteq \CC^n$, and where we use the identification $T^* \CC^n \cong T^* \widehat{\CC}^n$.
\end{thm}

By Theorem \ref{3} the two ideals in (\ref{95}) differ along minimal
components only by their multiplicities.  Taking into account this
information turns the $L$-characteristic variety $\ch^L(M_A(\beta))$
into the \emph{$L$-characteristic cycle} $\cc^L(M_A(\beta))$ of
$M_A(\beta)$.  Let $\mu^{L,\tau}_{A,0}(\beta)$ be the multiplicity of
$\Upsilon^\tau_A$ in $\cc^L(H_A(\beta))$.  This number is bounded from
below by the intersection multiplicity $\mu^{L,\tau}_A$ between the \emph{Euler
  variety}
\[
\Var(\gr^L(E_1,\ldots,E_d))\subseteq \CC^n
\]
and the component of $\gr^L(I_A)$ along $\Upsilon^\tau_A$. Moreover,
$\mu^{L,\tau}_{A,0}(\beta)$ agrees with this lower bound for a
Zariski-open set of parameters $\beta$, but may exceed it for special
values of $\beta$; see \cite{SchulzeWalther-duke}.

For $\tau\subseteq \tau'\in\Phi^{L,d-1}_A$, denote
\[
\pi_{\tau,\tau'}:\ZZ\tau'\to\ZZ\tau'/(\ZZ\tau'\cap\QQ\tau)\index{pt@$\pi_{\tau,\tau'}$}
\]
the natural projections, and define the polyhedra
\[
P_{\tau,\tau'}:=\conv(\pi_{\tau,\tau'}(\tau'\cup\{0\}))\index{Ptt@$P_{\tau,\tau'}$},\quad Q_{\tau,\tau'}:=\conv(\pi_{\tau,\tau'}(\tau'\smallsetminus\tau))\index{Qtt@$Q_{\tau,\tau'}$}.
\]
Using this notation, with volume functions normalized such that they return
unity on the standard simplex,
\[
\mu^{L,\tau}_A=
\sum_{\tau\subseteq\tau'\in\Phi^{L,d-1}_A}[\ZZ A:\ZZ\tau']\cdot[(\ZZ\tau'\cap\QQ\tau):\ZZ\tau]\cdot\vol_{\tau,\tau'}(P_{\tau,\tau'}\smallsetminus
Q_{\tau,\tau'})\geq 1.
\]
In particular, this formula proves that the slopes of the $D$-module
$M_A(\beta)$ are determined entirely by combinatorics of $A^L$, since
this is true for their $L$-characteristic varieties.
(For the empty face $\tau$, if $\NN A$ is saturated, this simplifies to the
formula already in \cite{GKZ89} that rank is then equal to the volume of $A$).

\begin{rmk}\label{rmk-slopes}
  If an $A$-hypergeometric system is homogeneous, it can have no
  slopes since it is regular holonomic \cite{HottaEq}. On the other
  hand, an inhomogeneous $H_A(\beta)$ has at least one slope along the
  subspace cut out by the variables corresponding to any of the faces
  of the umbrella of $A$ that do not touch the boundary of the
  umbrella, as moving it will eventually change the shape of the
  umbrella (compare \cite{SchulzeWalther-duke}).  By Laurent's
  results, regularity of $M_A(\beta)$ is hence equivalent to
  homogeneity and independent of $\beta$.  \schluss\end{rmk}

\begin{rmk}
\label{rem:SolsIrreg}
A natural question is whether
one can
find a stratification of the parameter space such that rank is
constant on each stratum and whether one can give a family of parametric
solutions that deform analytically to rank many solutions on the
chosen stratum. This is indeed so; the details are worked out in
\cite{BerkeschForsgardPassare-Mich14,BerkeschForsgardMatusevich-Adv16,BerkeschForsgardMatusevich-TAMS}.

For confluent systems, when the Nilsson ring does not contain all
solutions, the approach of Gevrey series can be used.
Early focus was on the irregularity sheaves of Mebkhout introduced in
\cite{Mebkhout-GrothendieckIII}. In a series of papers, Castro-Jimenez
and Fernandez-Fernandez
\cite{Fernandez-Adv10,CastroFernandez-TAMS11,CastroFernandez-JA11,CastroFernandez-MathZ12},
study theory and construction of solutions. Another point of interest
is asymptotics.
In
\cite{CastroGranger-IMRN15} it is worked out how this plays out in the
$d=1$ case ($A$ is a single row matrix): Gevrey series solution
along the singular locus of the system appear as asymptotes
of holomorphic solutions along suitable paths of integration. A
similar result for modified systems is proved in \cite{CFTT-TAMS15}.

A related problem is that of determining the monodromy of
$A$-hypergeometric systems. This turns out to be an extraordinarily
difficult problem, and only limited information is available at this
point. We mention the work of Ando, Esterov and Takeuchi
\cite{AndoEsterovTakeuchi-Adv15} that determines the monodromy at
infinity for confluent (inhomogeneous) systems, building on
\cite{Takeuchi-MathAnn10} for the homogeneous case. Hien's rapid decay
cycles (\cite{Hien-Inventiones09}) make an entry here via
\cite{EsterovTakeuchi-AJM15}, replacing the classical integral
representations of Gel{$'$}fand et al.
\schluss\end{rmk}

\section{Hodge theory of GKZ-systems}\label{sec:Hodge}

In this section we show that certain GKZ-systems carry a mixed Hodge
module structure in the sense of \cite{SaitoMHM} and investigate some
consequences of this fact. Since the definition of mixed Hodge modules
(MHM) is rather involved, we give here a simplified version which is
enough for our purpose. Assuming the reader to be at least somewhat
acquainted with the Riemann--Hilbert correspondence, we start with a
brief outline of the cornerstones of the theory of mixed Hodge
modules. We then give (certain) $A$-hypergeometric systems an
interpretation as Gau\ss--Manin systems and use it to define an MHM
structureon these $A$-hypergeometric systems. We then discuss two
induced filtrations on these GKZ-systems.

\subsection{Section setup, and basics on mixed Hodge
  modules}\label{subsec-setupHodge}

An algebraic mixed Hodge module on a smooth algebraic variety $X$ is
an algebraic, regular holonomic $\calD_X$-module $\calM$ together with
an increasing filtration by coherent $\calO_X$-modules $F^\Hodge_\bullet
\calM$ called the \emph{Hodge filtration} and an increasing
$\calD_X$-module filtration $W_\bullet \calM$ called the \emph{weight
  filtration}.  The $\calD_X$-module $\calM$ and the filtrations
$F^\Hodge_\bullet \calM$ and $W_\bullet \calM$ are required to satisfy rather
subtle compatibility conditions; in particular there are strong
conditions concerning the boundary behavior along every divisor of
$X$.  The category $\MHM(X)$ of algebraic mixed Hodge modules on $X$
is Abelian. Given a mixed Hodge $\calM$, its graded parts
\[
\Gr^W_k
(\calM) := W_k \calM / W_{k-1} \calM
\]
are pure Hodge modules. The category $\HM(X)$ of \emph{pure
  Hodge modules} is semi-simple; \emph{i.e.}, each graded part is a sum
a simple objects.  The simple $\HM(X)$-objects correspond via the de
Rham functor to intersection complexes $\IC_Y(\calL)$ supported on an
irreducible subvariety $Y$ of $X$, where $\calL$ is an irreducible
local system on an open, smooth subset of $Y$.  In particular, the
restriction of a pure Hodge module to the Zariski open set on which
the underlying $\calD$-module is smooth turns it to a variation of
pure Hodge structures on that smooth locus.

The standard example of a (mixed) Hodge module
on a smooth
variety $X$ is the structure sheaf $\calO_X$: it carries a canonical mixed
Hodge module structure, which satisfies
\begin{align*}
\Gr^{F^\Hodge}_p \calO_X :=& F^\Hodge_p \calO_X / F^\Hodge_{p-1} \calO_X = 0&\text{ if
}p \neq 0,\\
\Gr^W_p \calO_X =&\,0&\text{ for }p \neq  \dim X.
\end{align*}

\begin{ntn}\label{ntn-f}
  If $f\colon X\to Y$ is a morphism of smooth complex algebraic
  varieties, four basic functors on $\calD$-modules are induced. The
  most immediate one is the (left exact) na\"ive inverse image functor
  that arises from the chain rule, \cite[\S1.3]{HTT}. Its left derived
  functor, shifted by $\dim(X)-\dim(Y)$, is the \emph{inverse image
    functor} $f^+$ that is denoted by $f^\dagger$ in \cite[Rmk.~1.5.10]{HTT}. Conjugating
  $f^+$ by the holonomic duality functor from \cite[\S2.6]{HTT}
  leads to the \emph{exceptional inverse image} $f^\dagger$ that is denoted
  $f^\star$ in \cite[Dfn.~3.2.13]{HTT}.

  There is a \emph{direct image} functor as well, but its definition is
  more technical because the chain rule cannot be reversed in
  general. Again, one proceeds by defining a na\"ive version (neither
  left nor right exact) as in
  \cite[\S1.3]{HTT}, from which a derived functor $f_+$ can be
  defined; this functor is denoted $\int_f$ in
  \cite[p.~40]{HTT}. Conjugation by the duality functor leads to the
  \emph{exceptional direct image} functor $f^\dagger$, which is denoted
    $\int_{f!}$ in \cite[\S3.2]{HTT}.\schluss
\end{ntn}

Due to the groundbreaking work of
M.~Saito (\cite{Saito1, SaitoMHM}), for each morphism $f:X \to Y$ there are lifts of the functors
$f_+, f_\dag, f^+, f^\dag$ to the category of mixed Hodge modules
which we denote by
\begin{align*}
f_*, f_! : D^b \MHM(X) &\to D^b \MHM(Y)\\
f^{!}, f^*: D^b \MHM(Y) &\to D^b \MHM(X)\, .
\end{align*}
The proof of the existence of these functors on MHM require various
rather deep results from Hodge theory (such as the existence of a
Hodge structure on the cohomology of a degenerating VHS on a curve
which was established by S.~Zucker using $L^2$-cohomology), the theory
of filtered $\calD$-modules, compatibility properties of $V$- and
$F$-filtration (also known as strict specializability), as well as a
tricky formalism of induced modules.

\medskip

Our starting point is Section \ref{subsec-FGKZ}, where we have seen
that if $\beta \not \in \sRes(A)$ then $\widehat{\calM}_A(\beta)
\simeq (h_A)_+ \calO_{\TT}^\beta$.  So, in particular, if
$\calO^\beta_\TT$ is in $\MHM(X)$ then so is
$\widehat{\calM}_A(\beta)$ whenever $\beta\not\in \sRes(A)$.
Now in order for its (inverse) Fourier--Laplace transform
to be a mixed Hodge module,
the GKZ-system ${\calM}_A(\beta)$
should of
course in particular be regular holonomic. By Remark \ref{rmk-slopes}
and Definition \ref{dfn-homogeneous}, this property is equivalent to
$I_A$ being homogeneous. In other words, for the GKZ-system to have
any hope of being an MHM module we must require that the vector
$(1,1,\ldots ,1)$ is in the row span of $A$.  Fortuitously, this
requirement on $A$ provides also the solution to the translation of
MHM structures from $\widehat{\calM}_A(\beta)$ to
${\calM}_A(\beta)$. Indeed, while the (inverse) Fourier--Laplace
transform does in general not preserve mixed Hodge modules, we shall
employ a Radon transform (which makes only sense in the homogeneous case)
in order to construct a mixed Hodge module structure on the GKZ-system via
$\widehat{\calM}_A(\beta)$.

In order
to simplify the statement of some formulas in the remainder of the
article, we make now the following convention on $A$.
\begin{cnv}\label{cnv-newA}
   From now on, $A$ is in $ \ZZ^{(d+1) \times (n+1)}$ and we assume
   that $A$ is homogeneous, full, pointed, and generates a saturated
   semigroup.
\schluss\end{cnv}

Since a
GKZ-system derived from a pair $(A,\beta)$ is unchanged under an
invertible $\ZZ$-linear transformation of the rows we can moreover
assume that the matrix $A$ has the following shape
\begin{gather}\label{eqn-A,B}
A=\left(\begin{array}{c|ccc}1 & 1 & \ldots & 1 \\\hline
0 & \\
\vdots & & B & \\
0 &  \end{array}\right)
\end{gather}
where $B \in \ZZ^{d \times n}$ is full but is not necessarily  pointed or
homogeneous. Notice also that if $\NN A$ is saturated,
then so is $\NN B$; however, the converse implication is not true in general.

\subsection{Geometric interpretation of GKZ-systems}\label{subsec:geoHodge}

The aim of this section is to express certain GKZ systems as objects
which are built from consecutive applications of (possibly proper)
direct image and (possibly exceptional) inverse image functors applied
to a structure sheaf. From the discussion above it follows then that
these GKZ systems carry a mixed Hodge module structure. In order to
achieve this we have to introduce various integral transformations and
their relations.

Define a pairing
\begin{eqnarray}\label{eqn-exp}
\langle -, - \rangle\colon \widehat{\CC}^{n+1} \times \CC^{n+1}
&\to& \CC\\
(\pointy,\pointx) &\mapsto& \sum_{j=0}^n \pointy_j \pointx_j,\nonumber
\end{eqnarray}
and a free rank one $\calO_{\widehat{\CC}^{n+1} \times \CC^{n+1}}$-module
\[
\calL := \calO_{\widehat{\CC}^{n+1} \times \CC^{n+1}}
\cdot \exp\left({(-1)\cdot  \langle -, - \rangle}\right)
\]
which acquires a  $\calD_{\widehat{\CC}^{n+1} \times \CC^{n+1}}$-module structure
via the product rule. We denote by $p_1$ and $p_2$ the
projections from $\widehat{\CC}^{n+1} \times \CC^{n+1}$ to the first and second factor respectively. The
sheafified version of the Fourier-Laplace transform is given by
\begin{equation}\label{eq:FL}
 \FL (\calN) := p_{2+}(p_1^\dag \calN \overset{L}\otimes_\calO \calL)[n+1]
\end{equation}
 and one has $\FL\circ\FL = - \id$. Although defined at the level of derived
 categories, $\FL$ is an exact functor, and an instructive
 exercise shows that on the level of global sections it is given by
 formula \eqref{eq:FLconcrete}.
 Theorem \ref{thm:FL-GKZ-DirectImage} now implies that, whenever
 $\beta\not\in\sRes(A)$, we have
\[
\FL((h_A)_+ \calO_{\TT}^\beta) \simeq \FL^2(\calM_A(\beta)) \simeq \calM_A(\beta).
\]
Here, the final identification holds due to the homogeneity of $I_A$ even though
$\FL^2$ is not the identity.

\medskip

The second type of transformation we will need is the
\emph{Radon transformation} of $\calD$-modules introduced by Brylinski
\cite{Brylinski}; some variations were later discussed by d'Agnolo and
Eastwood \cite{AE}.

Let
\[
U:=\{\sum_{j=1}^n \pointy_j \pointx_j\neq 0\}\subseteq \PP(\widehat{\CC^{n+1}})\times\CC^{n+1}
\]
be the complement of the universal hypersurface
\[
Z:=\{\sum_{j=1}^n \pointy_j \pointx_j =0\}\subseteq \PP(\widehat{\CC^{n+1}})\times\CC^{n+1}
\]
defined by the vanishing of the pairing
$\langle-,-\rangle$. For the sake of readability,
we denote $\PP(\widehat{\CC^{n+1}})$ form now on simply by $\widehat{\PP}^n$. Consider
the following commutative diagram
\[
\xymatrix{& U \ar[dl]_{\pi_1^U}  \ar[d]^{j_U} \ar[dr]^{\pi_2^U} & \\ \widehat{\PP}^{n}  & \widehat{\PP}^{n} \times \CC^{n+1} \ar[l]_{\pi_1} \ar[r]^{\pi_2}& \CC^{n+1} \\ & Z \ar[ul]^{\pi_1^Z} \ar[u]_{i_Z} \ar[ur]_{\pi_2^Z} &}
\]
The Radon transformation is the functor
$\Radon\colon
D^b_{\rh}(\calD_{\widehat{\PP}^{n}}) \to D^b_{\rh}(\calD_{\CC^{n+1}})$ given by
\[
\Radon(\calN) := (\pi_2^Z)_+ (\pi_1^Z)^\dag \calN \simeq (\pi_2)_+ (i_{Z})_+ i_Z^\dag \pi_1^\dag \calN,
\]
and it permits variations $\Radon^\circ_c, \Radon_{\cst}:
D^b_{\rh}(\calD_{\widehat{\PP}^{n}}) \to D^b_{\rh}(\calD_{\CC^{n+1}})$  given
by
\begin{align*}
\Radon^\circ_{\mathrm c}(\calN) &:= (\pi_2^U)_\dag (\pi_1^U)^\dag \calN \simeq (\pi_2)_+ (i_Z)_+ i_Z^\dag \pi_1^\dag \calN \\
\Radon_{\cst}(\calN) &:= (\pi_2)_+ \pi_1^\dag \calN
\end{align*}
The adjunction triangle $(j_U)_\dag j_U^\dag \to \id \to (i_Z)_+ i_Z^\dag \overset{+1}\to$ gives rise to a triangle
\begin{equation}\label{eq-adjtriangle}
\Radon^\circ_c \to \Radon_{\cst} \to  \Radon \overset{+1}\to
\end{equation}

Let
\[
\pi: \widehat{\CC}^{n+1} \setminus \{0\} \to \widehat{\PP}^{n}
\]
be the canonical projection and denote by
\[
\pi_\VV\colon \VV \to \widehat{\PP}^{n}
\]
the total space of the tautological bundle
$\calO_{\widehat{\PP}^{n}}(-1)$. Recall that $\VV$ can be identified with the
blow-up of the point $\{0\}$ of $\widehat{\CC}^{n+1}$ and $\widehat{\PP}^{n}$ with the
exceptional divisor $E$. We denote by $\pi'_{\VV,E}\colon E \to
\{0\}\to \widehat{\CC}^{n+1}$ the restriction of the blow up map $\pi'_\VV: \VV
\to \widehat{\CC}^{n+1}$. The following proposition relates the Fourier-Laplace and
Radon transformations.
\begin{prp}\cite[Proposition 1]{AE}\label{prop-DE}
Let $\calN \in D^b_{\rh}(\calD_{\widehat{\PP}^{n}})$. There are the following isomorphisms
\begin{eqnarray*}
\Radon(\calN)&\simeq&\FL((\pi'_\VV)_+(\pi_\VV)^+\calN),\\
\Radon^\circ_{\mathrm c}(\calN)&\simeq&\FL(j_+ \pi^+ \calN),\\
\Radon_{\cst}(\calN)&\simeq&\FL((\pi'_{\VV,E})_+ \calN),
\end{eqnarray*}
where $j: \widehat{\CC}^{n+1} \setminus \{0\} \into \widehat{\CC}^{n+1}$ is the canonical
inclusion. \qed
\end{prp}
In particular, if $\calN$ is a mixed Hodge module, then
the above isomorphisms allow us to equip the right hand sides with induced
MHM structures.
\medskip

To simplify the presentation, we will focus  now (and this until Definition \ref{dfn-admissible} below) primarily on the case $\beta = 0$. For $\beta \neq 0$ a twisted variant of the Radon transformation is
needed: see \cite{ReiSe-Hodge} for details.
We start with the following commutative diagram
  \begin{gather}\label{eqn-h-diagram}
\begin{tikzcd}[ampersand replacement=\&,column sep=8ex]                                    \TT \arrow{d}{\pi_0} \arrow{r}[swap]{h'_A}
    \arrow[bend left]{rr}{h_A}
    \& \widehat{\CC}^{n+1} \setminus \{0\} \arrow{r}[swap]{j}
    \arrow{d}{\pi}
    \&  \widehat{\CC}^{n+1} \\
    \overline{\TT} \arrow{r}{g_B} \& \widehat{\PP}^n
    \end{tikzcd}
\end{gather}
where
\[
\pi_0\colon (\CC^*)^{d+1}=\TT \longrightarrow (\CC^*)^d=:\overline{\TT}
\]
is the projection to the last $d$ variables
and where
\begin{eqnarray}\label{eq:Embedding_g}
  g_B\colon \overline{\TT}&\hookrightarrow&\widehat{\PP}^n\\
(\pointt_1,\ldots,\pointt_d)=\pointt&\mapsto&
(1:\pointt^{\boldb_1}:\ldots:\pointt^{\boldb_n}).\nonumber
\end{eqnarray}
In particular,
\[
h_A\colon \TT\to \widehat{\CC}^{n+1}
\]
is as in \eqref{eqn-h_A} earlier (with the caveat that now $A$ is as
in Convention \ref{cnv-newA}).
We then observe that
$(h_A)_+ \calO_\TT \simeq (h_A)_+ \pi_T^+ \calO_{\overline{\TT}}
\simeq j_+ \pi^+ (g_B)_+ \calO_{\overline{\TT}}$, and
with Proposition \ref{prop-DE}, the isomorphisms
\begin{equation}\label{eq-GKZisHodge}
\calM_A(0) \simeq \FL((h_A)_+ \calO_{\TT}) \simeq
\Radon^\circ_{\mathrm c}((g_B)_+ \calO_{\overline{\TT}})
\end{equation}
endow the GKZ-system $\calM_A(0)$ with the structure of a mixed Hodge
module.

\medskip

We now consider a part of the long exact sequence of the adjunction
triangle \eqref{eq-adjtriangle} applied to $(g_B)_+
\calO_{\overline{\TT}}$.  In order to identify the individuals terms
we introduce a family of Laurent polynomials defined on $(\CC^*)^d
\times \CC^{n}=\overline\TT \times \CC^{n}$ using the columns
$\boldb_1,\ldots, \boldb_n$ of the matrix $B$ from \eqref{eqn-A,B}. We
define
\begin{eqnarray}\label{eq:DefVarphi}
\varphi: \overline{\TT} \times \CC^{n} &\to& \CC^{n+1} \\
(\pointt,\pointx) &\mapsto& (-\sum_{j=1}^{n} \pointx_j \pointt^{\boldb_i}, \pointx_1,\ldots,\pointx_n)
\end{eqnarray}

\begin{thm}[{\cite[Cor.~2.3]{Reich2}}]\label{theo:4TermSeq}
There is the following commutative diagram with exact rows where all
vertical maps are all isomorphisms; just for this statement we abbreviate for typesetting
reasons $g_B$ by $g$ and denote the Radon transform by just $\Radonshort$.
\[
\adjustbox{scale=0.89,center}{
\begin{tikzcd}[column sep=0.2in]
 \calH^{n}(\Radonshort_{\cst}(g_+\calO_{\overline{\TT}}))
 \ar[hookrightarrow]{r} \ar{d} &
 \calH^{n}(\Radonshort (g_+\calO_{\overline{\TT}})) \ar{r} \ar{d} &
 \calH^{n+1}(\Radonshort^\circ_{\mathrm c}(g_+\calO_{\overline{\TT}})) \ar[twoheadrightarrow]{r} \ar{d} & \calH^{n+1}(\Radonshort_{\cst}(g_+\calO_{\overline{\TT}}))  \ar{d}\\
 H^{d-1}(\overline{\TT};\CC) \otimes_\CC \calO_{\CC^{n+1}} \ar[hookrightarrow]{r} & \calH^0(\varphi_+ \calO_{\overline{\TT} \times \CC^{n-1}}) \ar{r} & \calM_A(0) \ar[twoheadrightarrow]{r}& H^{d}(\overline{\TT};\CC) \otimes_\CC \calO_{\CC^{n+1}}
\end{tikzcd}
}
\]
As a consequence, the lower exact sequence underlies a sequence of mixed Hodge modules.\qed
\end{thm}

\subsection{Hodge-filtration on GKZ-systems}
\label{sec:HodgeWeight}

Although the isomorphism \eqref{eq-GKZisHodge} equips the GKZ system
$\calM_A(0)$ with the structure of a mixed Hodge module, it is far
from clear what the Hodge and weight filtrations look like. The first
step in this direction was carried out by Stienstra \cite{Sti},
relying heavily on work of Batyrev \cite{Bat4}, who computed the Hodge
and weight filtration on the smooth part of the GKZ system.

Denote
\[
\Delta := \conv(\bolda_0,\ldots,\bolda_n)
\]
the convex hull of the points $\bolda_0,\ldots,\bolda_n$, and note that this is
the decone of the $A$-polyhedron from Definition \ref{96}.
Let $\tau \subseteq \Delta$ be a face of $\Delta$, let  $\pointx \in \CC^n$, and set
\[
F^\tau_{A,\pointx} := \sum_{j : a_j \in \tau} \pointx_j \boldt^{a_j}.
\]
The Laurent polynomial
$F_{A,\pointx}:=F^A_{A,\pointx}$ is called
\emph{non-degenerate} (see, \emph{e.g.}, \cite[Definition 3.3]{Bat4}) if for
every face $\tau$ of $\Delta$ the equations
\[
F^\tau_{A,x} = t_0 \frac{\del}{\del t_0}( F^\tau_{A,x}) = \ldots = t_d
  \frac{\del}{\del t_d}( F^\tau_{A,x}) = 0
\]
have no common solutions in $\TT$. Then, for $0\le i\le d$,
define
the differential operators
\[
P_i := \sum_{j=0}^n \left(a_{i,j} \pointx_j \boldt^{\bolda_j} + t_j
\del_{t_j} \right)
\]
which are elements of the Weyl algebra $D_{\CC[\boldt^\pm]}$ on $t_0,\ldots,t_d$
localized at $t_0\cdots t_d$.
One checks that these operate on
the semigroup ring $S_A\subseteq \CC[t_0^{\pm 1}, \ldots, t_d^{\pm 1}]$,
$P_i(S_A)\subseteq S_A$, so  they are
differential operators on the affine toric variety $X_A=\Spec(S_A)$.

\medskip

Before we can state Stienstra's result mentioned in the introduction
to this section, we need some more terminology. Let
\[
I_\Delta^{(0)} \subseteq I^{(1)}_\Delta \subseteq \ldots \subseteq I_\Delta^{(d+1)} \subseteq I_\Delta^{(d+2)} = S_A
\]
be the ascending sequence of homogeneous ideals in $S_A$ where
$I_\Delta^{(k)}$ is generated by all elements $\boldt^\bolda$ with
$\bolda \in \NN A$ that are not contained in any codimension $k$ face
of $\RR_{\geq 0} A$. Define a decreasing sequence of $\CC$-vector
spaces in $S_A$
\[
\cdots \supseteq \calE^{-k} \supseteq \calE^{-k+1} \supseteq \cdots \supseteq \calE^{-1} \supseteq \calE^0 \supseteq \calE^1 = 0
\]
where $\calE^{-k}$ is spanned by monomials $\boldt^{\bold c}$ such that
$\bold c = (c_0,\ldots,c_d) \in \NN A$ satisfies $c_0 \leq k$.

Stienstra proved the following result
\begin{thm}\cite{Bat4,Sti,ReiSe-Hodge}\label{theo:Batyrev}
Let $\pointx \in \CC^{n+1}$ be such that the Laurent polynomial
$F_{A,\pointx}$ is non-degenerate and consider the canonical inclusion
$i_\pointx : \{\pointx\} \into \CC^{n+1}$. Then, with $\varphi$ denoting
the family from \eqref{eq:DefVarphi},
\[
 H^{d}(\overline{\TT},\varphi^{-1}(\pointx);\CC) \simeq
 i_\pointx^+ \calM_A(0) \simeq S_A / \sum_{i=0}^d P_i S_A.
\]
Under this isomorphism, the Hodge filtration  is given by
\begin{align}\label{eq:BatyrevHodge}
  F^{d-k} H^d(\overline{\TT},\varphi^{-1}(\pointx);\CC)
  &&\simeq  &&
  \image\left( \calE^{-k}\to  S_A / \sum_{i=0}^d P_i S_A\right).
\end{align}
If the matrix $B \in \ZZ^{d\times n}$ is homogeneous, then the weight filtration on $H^{d}(\overline{\TT},\varphi^{-1}(x);\CC)$ is given by
\begin{align}\label{eq:BatyrevWeight}
  W_{k+d-1} H^{d}(\overline{\TT},\varphi^{-1}(x);\CC)
  && \simeq&  &
  \image\left(I^{(k)}_\Delta\to  S_B / \sum_{i=1}^{d} P_i S_B\right),
\end{align}
where the semigroup ring $S_B$, the ideals $I^{(k)}_\Delta$ and the
differential operators $P_i$ are now derived from $B$.\qed
\end{thm}
Equation \eqref{eq:BatyrevHodge} is shown in \cite{Sti} for
homogeneous $A$; the general case is treated in \cite{ReiSe-Hodge}.

The surjection
$D_A \to M_A(\beta)$  induces from the order filtration
$F^\ord_\bullet$ on $D_A$ a  filtration
on $M_A(\beta)$ which we denote by $F^{\ord}_\bullet M_A(\beta)$; we
proceed
similarly to define a filtration $F^{\ord}_\bullet$ on the sheaf
$\calM_A(\beta)$. The following theorem gives a comparison between
this order filtration and the Hodge filtration $F^\Hodge_\bullet
\calM_A(\beta)$ (in the sense of mixed Hodge modules), this extends
the first part of the above Theorem \ref{theo:Batyrev}. Since we will
formulate the result for certain parameter vectors $\beta$ different
from $0$, we first need to introduce the following definition.

\begin{dfn}\label{dfn-admissible}
  The set of \emph{admissible parameters} $\beta \in \RR^{d+1} \subseteq \CC^{d+1}$ is defined by
\[
\mathfrak{A}_A:= \bigcap_{\tau : \tau \text{ facet}} \{\RR \cdot \tau - [0,\frac{1}{e_\tau}) \cdot \boldeps_A \}
\]
where $\boldeps_A := \bolda_0 + \ldots +\bolda_n$, $e_\tau:= \langle
n_\tau,\boldeps_A\rangle \in \ZZ_{>0}$ and $n_\tau$ is the unit, inward
pointing, normal vector of $\tau$.
\schluss\end{dfn}

\begin{exa}\label{exa-twisted-cubic}
For the matrix
\[
A = \left(\begin{matrix}1 & \phantom{-}1 & 1 & 1 \\0 & -1 & 1 & 2  \end{matrix}\right),
\]
the following picture
\[
\newdimen\scale
\scale=0.8cm
\begin{tikzpicture}
 \filldraw[blue,opacity=.2] (0,0) -- (\scale*-2.05,\scale*2.05) -- (\scale*4.05,\scale*2.05)  -- (0,0);
 \filldraw[cyan,opacity=.8] (0,0) -- (\scale*-.66,\scale*-.33) -- (\scale*-.33,\scale*-.66) -- (\scale*.33,\scale*-.33)  -- (0,0);

 \foreach \x in {-3,...,4}{
   \foreach \y in {-2,...,2}{
     \node[draw,circle,inner sep=0.5pt,fill] at (\scale*\x,\scale*\y) {};
   }
 }
\draw (\scale*-3,\scale* -2) -- (\scale*4, \scale *1.5);
\draw (\scale*-2,\scale* -2) -- (\scale*4, \scale *1);
\draw (\scale*-1,\scale* -2) -- (\scale*4, \scale *0.5);
\draw (\scale*0,\scale* -2) -- (\scale*4, \scale *0);
\draw (\scale*1,\scale* -2) -- (\scale*4, \scale *-0.5);
\draw (\scale*2,\scale* -2) -- (\scale*4, \scale *-1);
\draw (\scale*3,\scale* -2) -- (\scale*4, \scale *-1.5);
\draw (\scale*-3,\scale* 2) -- (\scale*1, \scale *-2);
\draw (\scale*-3,\scale* 1) -- (\scale*0, \scale *-2);
\draw (\scale*-3,\scale* 0) -- (\scale*-1, \scale *-2);
\draw (\scale*-3,\scale* -1) -- (\scale*-2, \scale *-2);

\draw(-2.5,-2.35)--(-1.9,-2.35);
\node[] at (-1.0,-2.39) {$\sRes(A)$};
\filldraw[blue,opacity=.2] (0,-2.2) -- (0.6,-2.2) -- (0.6,-2.5) -- (0,-2.5) -- (0,-2.2);
\node[] at (1.2,-2.41) {$\RR_{\geq 0} A$};
\filldraw[cyan,opacity=.8] (2,-2.2) -- (2.6,-2.2) -- (2.6,-2.5) -- (2,-2.5) -- (2,-2.2);
\node[] at (3.05,-2.41) {$\mathfrak{A}_A$};
\end{tikzpicture}
\]
shows the sets $\sRes(A)$ (see Definition \ref{def-sRes} above) and
$\mathfrak{A}_A$.
\schluss\end{exa}

We can now state a result, taken from \cite[Theorem 5.35]{ReiSe-Hodge} which describes the Hodge filtration on the GKZ-systems in a rather precise way.
\begin{thm}\label{thm:HodgeOnGKZ}
Let $A \in \ZZ^{(d+1) \times (n+1)}$ be as in Convention \ref{cnv-newA}, $\beta \in
\mathfrak{A}_A$ and $\beta_0 \in (-1,0]$. Then the Hodge filtration on
$\calM_A(\beta)$ is given by the shifted order filtration, so that
we have the following equality of filtered $\calD_{\CC^{n+1}}$-modules
\[
(\calM_A(\beta) ,F^\Hodge_\bullet) = (\calM_A(\beta), F^{\ord}_{\bullet+d})
\]
\end{thm}
It has been shown in \cite[Theorem 5.43]{ReiSe-Hodge} that the first
part of the above Theorem \ref{theo:Batyrev}, and so Formula
\eqref{eq:BatyrevHodge} is a rather direct consequence of the
comparison between the Hodge and the order filtration on $\calM_A(0)$.
\begin{rmk}
As already noted in Section \ref{sec:EK-Rank} above, a variant of Borisov--Horja's better behaved GKZ-systems has been considered in \cite{Mo15}.
If we suppose that $A$ is normal (as we do throughout this section),
then the definition in \cite{Mo15} coincides with the one for ordinary
GKZ-systems as given in \ref{def:GKZ} above.
However, the matrix $A$ is not supposed to be homogeneous in \cite{Mo15}.
The module $\calM_A(\beta)$ will have irregular singularities then, as discussed in Section \ref{sec:Irreg} above. One may ask what kind of Hodge theoretic information can be derived from  $\calM_A(\beta)$ in this case. This is similar to the statements on the ordinary versus irregular Hodge filtration on univariate hypergeometric systems that we will discuss below.

In \cite[Prop.~1.4]{Mo15}, Mochizuki proves the the following statement which can be considered as an irregular variant of Theorem \ref{thm:HodgeOnGKZ} above.
Let $B\in\ZZ^{d\times n}$ be such that $\ZZ B=\ZZ^d$. Suppose for the simplicity of the exposition that $\NN B=\RR_{\geq 0} B\cap \ZZ^d$.
Consider the non-commutative ``Rees ring''
\begin{eqnarray}\label{eq-Rees}
\Rees_{\CC \times \CC^n} &=& \CC[z,x_1,\ldots,x_n]\langle z^2 \del_z, z\del_{x_1},\ldots, z \del_{x_n}\rangle
\end{eqnarray}
and
the corresponding sheaf $\shRees_{\CC\times\CC^n}$.
Let $\calH_A^z(0)$ be the left $\shRees_{\CC\times\CC^n}$-ideal generated by
\begin{eqnarray}
  \label{eq:zGKZ-ideal}
  \hspace*{.1\textwidth}
\begin{minipage}{.8\textwidth}
  \begin{eqnarray*}
    \widehat{E}_0 &:= & z^2 \del_z + \sum_{j=1}^nz x_j \del_{x_j}; \\
\widehat{E}_i &:= & \sum_{j=1}^n a_{i,j}\, z x_j \del_{x_j} \qquad\qquad \text{for} \quad k=1,\ldots, d;\\
\widehat{\Box}_{\bold u} &:= &\prod_{j: u_j > 0} (z \del_{x_j})^{u_j}
- \prod_{j: u_j < 0} (z \del_{x_j})^{-u_j} \qquad \text{for all} \; \boldu
\in \ker(B).
\end{eqnarray*}
\end{minipage}
\end{eqnarray}
Then the left $\shRees_{\CC\times\CC^n}$-module
$\shRees_{\CC\times\CC^n}/\calH_A^z(0)$ underlies
a \emph{mixed twistor module} on $\CC^n$, a notion that in many
respects is the correct replacement of a mixed Hodge module in the irregular setup. In particular, any mixed Hodge module can be considered as a special mixed twistor module, and therefore the case $\beta=0$ of Theorem
\ref{thm:HodgeOnGKZ} can be deduced from Mochizuki's result. Using a
filtered variant of the Fourier--Laplace transformation (compare the
discussion in Section \ref{sec:MirrorSym} below), one can also obtain the latter from Theorem \ref{thm:HodgeOnGKZ}, as has been demonstrated in
\cite[Corollary 4.8]{SevCastReich}.
\schluss\end{rmk}

As another application of Theorem \ref{thm:HodgeOnGKZ}, we will describe some results about the Hodge structure of univariate hypergeometric equations (see the discussion in Subsection \ref{subsec-dim-red} above). Consider again the operator
\begin{equation}
\label{eq:OneDimHyperOp}
P=\prod_{i=1}^{m'}(\theta_z-\lambda_i)-z\cdot\prod_{j=1}^m(\theta_z-\mu_j) \in \CC[z]\langle \partial_z \rangle
\end{equation}
(compare with equation \eqref{eqn-gauss-diffeq}, where $m'=q+1$, $m=p$ and where $\lambda_1=0,\lambda_i=1-\beta_{i+1}$, $\mu_j=-\alpha_j$)
for some \emph{real} numbers $\lambda_i,\mu_j$.  The corresponding cyclic module
\[
\calH(\lambda;\mu):=\calD_{\AA^1}/\calD_{\AA^1}\cdot P,
\]
is irreducible if and only if for all $i,j$ we have $\lambda_i-\mu_j\notin \ZZ$.
The modules $\calH(\lambda;\mu)$  are the most basic examples of \emph{rigid} $\calD$-modules (see \cite{KatzN1,Arinkin}). A first consequence of this property is that if $\calH(\lambda;\mu)$ is irreducible, then it is isomorphic to some
$\calH(\lambda';\mu')$ whenever $\mu-\mu'$ and
$\lambda-\lambda'$ are integer vectors.  We can thus assume that
$0\leq \lambda_1 \leq \ldots, \lambda_{m'} <1$, $0\leq \mu_1 \leq \ldots
\leq \mu_m <1$ and that $\lambda_i\neq \mu_j$ for all $i,j$. It is
obvious that $\calH(\lambda;\mu)$ is regular exactly when $m'=m$ and in
that case it has the three singular points $\{0,1,\infty\}$. On the
other hand, if $m'\neq m$ then $\Sing(\calH(\lambda;\mu) =
\{0,\infty\}$.

In the regular case, that is, if $m'=m$, the rigidity property can be stated at the level of the
the local system $\calL$ on
$\PP^1\backslash \{0,1,\infty\}$ of solutions of $P$: it simply says that the local monodromies around the singular points determine the (global) monodromy representation defined by $\calL$. From there it follows
by \label{page:SimpsonFedorov}\cite[Cor.~8.1]{Si2} and also
\cite[Prop. 1.13]{DeligneRigid} that $\calL$
underlies a
\emph{complex variation of Hodge structures}. Then the following
formula for its Hodge numbers has been shown in \cite[Thm.~1]{Fedorov}
\begin{equation}\label{eq:Fedorov}
\begin{array}{rcl}
\dim\gr_k^{F^\Hodge}\calL&:=&\dim(F^\Hodge_k \calL/F^\Hodge_{k+1}\calL) \\ \\
&=& \# \left\{s:1\leq s \leq m', k=\# \{i:\lambda_i < \mu_s\}-s\right\}.
\end{array}
\end{equation}
The Picard-Fuchs equation of the family of elliptic curves in Example
\ref{exa-elliptic} corresponds, as we computed there, to the
hypergeometric differential equation given by the module
$\calH(0,0;1/2,1/2)$. Applying Fedorov's formula yields $\dim(\gr_0^F
\calL)=\dim(\gr_1^F \calL)=1$, confirming our computation in Example
\ref{exa-elliptic}.
Notice also that in this case the local system
$\calL$ underlies a real (and even rational) variation of Hodge
structures, which is consistent with \cite[Theorem 2]{Fedorov}.

If $m'\neq m$ (and, up to a change of the
coordinate $z\mapsto 1/z$ we can assume that $m'>m$),
then $\calH(\lambda;\mu)$ is irregular and can no longer support a
variation of Hodge structures. In \cite{Sa15}, a category of
\emph{irregular Hodge modules} is developed, which can roughly be seen
as lying between the category of mixed Hodge modules and the category
of mixed twistor modules. A possibly irregular $\calD_X$-module
$\calM$ on a complex manifold $X$ underlying an irregular Hodge module
comes equipped with an \emph{irregular Hodge filtration}, an
increasing filtration $F^{\irr}_\alpha \calM$ by coherent
$\calO_X$-modules indexed by the real numbers (contrarily to the regular case); we write
$F^{\irr}_{<\alpha} \calM:=\bigcup_{\beta<\alpha}
F^{\irr}_\beta\calM$.  However, the
indexing set is determined by a finite set $I\subseteq [0,1)$
  having the property that
  \[
  \gr^{F^{\irr}}_\alpha \calM := F^{\irr}_\alpha
  \calM / F^{\irr}_{<\alpha} \calM=0 \qquad\text{ if }\alpha\notin
  I+\ZZ.
  \]
  In \cite{SaYu18}, the following formula for the irregular Hodge
  numbers has been found (see also \cite{CDS} and \cite{SevCastReich},
  where the Hodge filtration itself is determined in some cases, using
  Theorem \ref{thm:HodgeOnGKZ} from above):
  \begin{align}\label{eqn:SaYu}
\dim \gr_\alpha^{F^{\irr}} \calH(\lambda;\mu) &= \# \left\{s:1\leq s\leq
m', \alpha =
\# \{i:\mu_i < \lambda_s\}
+
(m'-m)\alpha_s-s
\right\}.
\end{align}
For $m'=m$, this gives back the formula
\eqref{eq:Fedorov} up to the fact that the local system $\calL$ is in the
regular case in \cite{Fedorov} the one of the solutions of
$\calH(\lambda;\mu)$, whereas  formula \eqref{eqn:SaYu} gives (for $m'=m$) Hodge
numbers of a filtration defined on the dual local system of flat
sections.

\subsection{Weight filtration on GKZ systems}
\label{subsec-weight}

In the remainder of this section, we discuss results concerning the
weight filtration on GKZ-systems. Recall that we equipped the GKZ-system
$\calM_A(0)$ in subsection \ref{subsec:geoHodge} with a mixed Hodge module
structure by rewriting it as certain Radon transform of a direct image
of a structure sheaf (cf. \eqref{eq-GKZisHodge}). In this subsection
we endow the GKZ systems with an \emph{a priori} different mixed Hodge module structure.
If the matrix $A$ is chosen to be homogeneous then the GKZ-system $\calM_A(0)$ is a monodromic $\calD$-module. In this case the Fourier--Laplace transformation can be replaced by the Fourier--Sato transformation (or monodromic Fourier--Laplace transformation) (cf. \cite[Th\'{e}or\`{e}me 7.24]{Brylinski}) which happens to be a functor of mixed Hodge modules.

Denote by
\[
\theta: \CC^*  \times \widehat{\CC}^{n+1} \longrightarrow \widehat{\CC}^{n+1}
\]
the standard $\CC^*$-action on $\widehat{\CC}^{n+1}$. We refer to the push-forward $\theta_*(z \del_z)$ as the Euler vector field $\mathfrak{E}$, where $z$ is a coordinate on $\CC^*$. A regular holonomic $\calD$-module $\calM$ is called \textit{monodromic}, if the Euler field $\mathfrak{E}$ acts finitely on the global sections of $\calM$.

Consider the diagram
\[
\xymatrix{ & \widehat{\CC}^{n+1} \times \CC^{n+1} \ar[dl]_{p_1} \ar[dr]^{\omega} & & \\
\widehat{\CC}^{n+1} & & \CC_z \times \CC^{n+1} & \{0\} \times \CC^{n+1} \ar[l]_{i_0}}
\]
where $p_1$ is the projections to the first factor, $i_0$ is the canonical inclusion and the map $\omega$ is given by
\begin{align*}
\omega: \widehat{\CC}^{n+1} \times \CC^{n+1} &\longrightarrow \CC_z \times \CC^{n+1}\\
(\pointy, \pointx) &\mapsto (\pointz= \sum_i \pointy_i \pointx_i,\; \pointx)
\end{align*}
The Fourier--Sato transformation (or monodromic Fourier transformation) is defined by
\begin{align*}
\FS: \MHM(\widehat{\CC}^{n+1}) &\longrightarrow \MHM(\CC^{n+1}) \\
\calM &\mapsto \phi_z \omega_* p_1^* \calM [n+1]
\end{align*}
where $\phi_z$ is the vanishing cycle functor along $z=0$.

It was shown in \cite[Proposition 4.12]{RW-weight} that the Fourier--Sato transformation respects the weight filtration of monodromic $\calD$-modules which are localized along $\{0\} \in \widehat{\CC}^{n+1}$ (up to a shift). Hence, a weight filtration on the GKZ-system is induced by the following isomorphisms:
\[
W_{k+n+1} \calM_A(0) := W_{k+n+1} \FS((h_A)_+ \calO_{\TT}) \simeq \FS(W_k (h_A)_+ \calO_{\TT})
\]

Since the Fourier--Sato transform is an equivalence of categories it is therefore enough to compute the weight filtration on $\widehat{\calM}_A(0)= (h_A)_+ \calO_{\TT}$ which will be done below.

Recall that the graded parts $\Gr^W_k \calM$ of a mixed Hodge module
are pure Hodge modules and as such are semi-simple, splitting as
direct sums of intersection complexes (which are simple
$\calD$-modules). Because the number of simple objects (counted with
multiplicity) is independent on the chosen (weight) filtration this
also gives us the simple objects occurring in the weight filtration
induced by the Radon transform (but possibly in another
order). However, we conjecture that the Fourier--Sato transformation
and the Radon transformation are actually isomorphic on the level of mixed
Hodge modules.
\begin{cnj}
For $\calN \in \MHM(\PP^n)$:
\[
\FS(j_*\pi^!N) \simeq \Radon^\circ_c(\calN)
\]
\schluss\end{cnj}

We will now proceed to state the result on the weight filtration of $\widehat{\calM}_A(0) = (h_A)_+ \calO_{\TT}$:

Let $ \tau \subseteq \gamma \subseteq \sigma$ be faces of a cone $\sigma \subset \RR^{d+1}$.  The quotient face of $\gamma$ by $\tau$ is defined as:
\[
\gamma/\tau := (\gamma + \tau_\RR)/\tau_\RR \subseteq \RR^{d+1} / \tau_\RR
\]
where $\tau_\RR$ is the linear span of the cone $\tau$. Define
\[
\gamma^{\mho} := \{ f \in \Hom_\RR(\RR^{d+1},\RR) / \gamma^\bot \mid f(\pointx) \geq 0\;\; \forall \pointx \in \gamma \}
\]
The cone $\gamma^\mho$ is the dual of $\gamma$ in its own span, hence independent of $\sigma$. For cones $\tau \subseteq \gamma$ denote by $X_{\gamma/\tau}$ the spectrum of the semigroup ring induced by the cone $\gamma /\tau$ in its natural lattice. Set $Y_{\gamma/ \tau}:= X_{(\gamma/\tau)^\mho}$.

In the following, we denote the cone $\RR_{\geq 0} A$ by $\sigma$.
The Fourier--Laplace transformed GKZ system $\widehat{\calM}_A(0)$ is
isomorphic to $(h_A)_+ \calO_{\TT}$ and has support on the affine
toric variety $X_A = X_\sigma$. For a face $\tau$ of $\sigma$ write
$d_\tau$ for its dimension.  We have seen in Subsection
\ref{subsec-torusaction} that the $d_\tau$-dimensional $\TT$-orbits
$O^\tau_A$ in $X_\sigma$ are in one-to-one correspondence with the
faces $\tau$ of $\sigma$. The closure of an orbit $O^\tau_A$ is
$X_\tau$.

It turns out that the varieties $X_\tau$ are exactly those which occur
as support varieties of the summands in the semisimple decompositions of the graded parts  $\gr^W\widehat{\calM}_A(0)$.

Let $\calL_{(\tau,d+e)}$ be the constant local system of rank $\dim
\IH^{d+1-d_\tau-e}(Y_{\sigma/\tau})$ on $O^\tau_A$.  In order to
simplify the notation, we use the symbol $\IC_Y(\calL)$ for the
intersection cohomology $\calD$-module on some smooth variety $X$ with
support on the closed subset $Y\subseteq X$, and where $\calL$ is a
local system on a Zariski open subset of $Y$.
\begin{thm}
Let $A \in \ZZ^{(d+1) \times (n+1)}$ be full, pointed, saturated, but
not necessarily homogeneous.
The weight graded parts of the mixed Hodge module $\widehat{\calM}_A(0)$ are given by
\[
\gr^W_{d+1+e} \widehat{\calM}_A(0) \simeq \bigoplus_{\tau \subseteq \sigma} \IC_{X_\tau}(\calL_{(\tau,d+1+e)}).
\]
\end{thm}

\begin{cor}
Let $A \in \ZZ^{(d+1) \times (n+1)}$ be as above. The length of the GKZ system $\calM_A(0)$ is
\[
\sum_{\tau \subseteq \sigma} \sum_{e=0}^{d+1-d_\tau} \dim \IH^{e}(Y_{\sigma/\tau}) = \sum_{\tau \subseteq \sigma} \dim \IH^*(Y_{\sigma / \tau}).
\]
\end{cor}

\section{Application to toric mirror symmetry}
\label{sec:MirrorSym}

The aim of this final section is to discuss some
results concerning the so-called mirror symmetry phenomenon,
which links enumerative geometry of projective algebraic, and more
generally symplectic varieties (called \sideA-model) to complex
geometry, in particular, Hodge theory of their so-called
\sideB-models. The \sideB-model is usually given by a family of algebraic
varieties which may have singularities and which need not be
projective (which forces one
to consider compactifications, see below). Often these families on the \sideB-side
are referred to as \emph{Landau--Ginzburg models}.

 The first example of mirror symmetry was given by Candelas, de la
 Ossa, Green and Parkes \cite{Candelas} who predicted a virtual number
 of rational curves on a quintic threefold (later referred to as the
 genus $0$ Gromov--Witten invariants) by period computations for the
 mirror partner (the \sideB-model). These predictions were
 verified and also generalized to numerically effective smooth
 complete intersections in toric varieties by Givental \cite{Giv6},
 \cite{Giv7}. His celebrated mirror theorem shows that the
 $J$-function, a generating function for the genus $0$ GW-invariants
 of such varieties, is computable in terms of a cohomology-valued
 hypergeometric function. Givental also conjectured that the
 components of this function are given as oscillating integrals. This
 was much later proved by Iritani in \cite{Ir2} (even treating the
 case where the toric variety in question is an orbifold), some
 details of the construction described below are parallel to his
 paper. However, an algebraic construction of the correct Hodge
 theoretic \sideB-model was still missing. Our purpose in this section
 is to give an overview of techniques and results (mainly referring to
 \cite{ReiSe, ReiSe2, ReiSe-Hodge} as well as to \cite{Mo15}), where
 the machinery of GKZ-systems as discussed in the previous sections is
 used to obtain a purely algebraic Hodge theoretic (and $\calD$-module
 based) mirror correspondence for certain smooth toric varieties
 resp. subvarieties of them.

\subsection{Gromov--Witten invariants and Dubrovin connection}
\label{subsec-GW+Dub}

Let $X$ be a toric smooth projective variety.  For the purpose of this
exposition, we assume further that $X$ is Fano, so the anticanonical
class $[-K_X]$ is ample. A good part of the results discussed below
also applies if one considers \emph{weak Fano} manifolds, meaning that
$[-K_X]$ is a numerically effective (nef) class. There are however a
few technical modifications needed in the nef case, which is why we
refrain from discussing it here. Developing the mirror symmetry
picture described below in the absence of any positivity assumption on
$X$ remains a subject of active current research (see, \emph{e.g.},
\cite{Ir4}, \cite{GKR}, \cite{IritaniBigQCohom}).

Let $\beta \in H_2(X,\ZZ)$ and choose $\gamma_1, \gamma_2, \gamma_3 \in H^*(X,\QQ)$. The genus zero, three point Gromov--Witten invariants
\[\label{eq:GWinv}
\langle I_{0,3,\beta}\rangle(-,-,-): H^*(X,\QQ)^{\otimes 3} \to \QQ
\]
intuitively count the number of stable maps $f$ from rational curves
$C$ with---in this case---three marked points, satisfying $f_*([C]) =
\beta$ and $f(C) \cap \PD(\gamma_i) \neq \emptyset$ for $i=1,2,3$. (Here and
elsewhere, $\PD(-)$ denotes the Poincar\'e dual). Technically, they
are obtained as follows: pull back the (three) arguments of $\langle
I_{0,3,\beta}\rangle$ to the moduli space of such maps (along the three induced
evaluation maps to $X$), take their cup product and evaluate against this product by integration over a certain
\emph{virtual fundamental class} on the moduli space. Constructing
this latter class is a major issue in Gromov--Witten theory (see,
\emph{e.g.} \cite{FP97} and \cite{BF97}).

\medskip

We choose a homogeneous basis $T_0,T_1,\ldots,T_r,T_{r+1},\ldots, T_s$
of $H^*(X;\ZZ)$ such that $T_0  \in H^0(X;\ZZ)$, the classes
$T_1,\ldots, T_r \in H^2(X;\ZZ)$ lie in the nef cone of $X$ and
$T_{r+1},\ldots, T_s \in H^{>2}(X;\ZZ)$. Let $g_{ij} := (T_i,T_j)$ be
the Poincar\'{e} pairing between the elements $T_i$ and $T_j$ and
define
\[
T^i := \sum_j g_{ij}T_j.
\]

With $\delta \in H^2(X;\CC)$, the three point Gromov--Witten invariants can be used as structure
constants for a family of multiplications
\begin{equation}\label{eq-smallquantum}
\gamma_1 \ast \gamma_2 := \sum_{\beta \in H_2(X,\ZZ)}\sum_{i=0}^s
\exp({\delta(\beta)})\cdot \langle I_{0,3,\beta}\rangle(\gamma_1,\gamma_2,T_i)T^i
\end{equation}
on $H^*(X;\CC)$. This product structure is  the
\emph{small quantum product} of $X$ and
 parameterized
by the cosets of $\delta$ in the \emph{complexified K\"ahler moduli space}
\[
\calK:= H^2(X;\CC)/2\pi \sqrt{-1}\cdot H^2(X,\ZZ).
\]

\emph{A priori} it is far from clear that the sum in
\eqref{eq-smallquantum} is convergent. However, the Gromov--Witten
invariants satisfy (among others) the following properties:
\begin{align*}
&\textbf{Effectivity}: &&\langle I_{0,3,\beta}\rangle = 0 \qquad
  \text{if $\beta$ does not lie in the Mori cone}\\
&\textbf{Degree}: && \langle I_{0,3,\beta} \rangle (T_i,T_j,T_k) = 0 \qquad \text{unless} \quad \sum_{i=1}^3 \deg(T_i) = 2 \dim X - 2 c_1(X)(\beta) \\
&\textbf{Point Mapping}: && \langle I_{0,3,0}\rangle(T_i,T_j,T_k) = (T_i \cup T_j \cup T_k)([X])
\end{align*}
where we recall that the Mori cone is the cone in $H_2(X;\RR)$ of effective classes of curves. It is dual to the cone of nef divisors in $H^2(X;\RR)$.
The effectivity axiom together with our assumption that $X$ be Fano---
so that the class $c_1(X)$ be ample---show that $\langle
I_{0,3,\beta}\rangle$ is zero unless $c_1(X)(\beta) \geq 0$. The
degree axiom now tells us that for fixed $T_i,T_j,T_k$ there are only
finitely many $\beta$ in the Mori cone such that $\langle
I_{0,3,\beta}\rangle(T_i,T_j,T_k)$ is non-zero. Hence the product
defined in \eqref{eq-smallquantum} is finite and therefore defined on
the whole space $\calK$.

It can be seen from other axioms that the small quantum product is
commutative, associative and that $T_0$ acts as identity. Let
$\eta_1, \ldots, \eta_r \in H_2(X,\ZZ)$
such that $T_i(\eta_j)$ is the Kronecker $\delta_{i,j}$
for $1\le i,j\le r$. If we write
\begin{eqnarray*}
\delta &=& t_1 T_1 + \ldots+t_rT_r \in H^2(X;\CC),\\
\beta&=&\beta_1\eta_1+\ldots+\beta_r\eta_r\in H_2(X;\CC),
\end{eqnarray*}
and set $q_i := \exp({t_i})$ for $i=1,\ldots ,r$,
we get
\[
\exp({\delta(\beta)}) = q_1^{\beta_1}\ldots q_r^{\beta_r}.
\]
Then, under the exponential map
from $H^2(X;\CC)$ to $\calK$,
$\mathboldq=\{q_i\}_{i=1,\ldots,r}$ become  coordinates on $\calK$
corresponding to $\boldt=\{t_i\}_{i=1,\ldots,r}$ on $H^2(X;\CC)$ and induce an explicit
isomorphism $\calK \simeq (\CC^*)^r$. Since $T_1,\ldots,T_r$ lie in
the nef cone, the cone generated by the dual basis
$(\eta_j)_{j=1,\ldots,r}$ contains the Mori cone and therefore all
monomials $q_1^{\beta_1}\ldots q_r^{\beta_r}$ have non-negative
exponents. Hence the quantum product extends to the partial
compactification
\begin{equation}\label{eq:CompactK}
\overline{\calK}:= \CC^r \hookleftarrow (\CC^*)^r =
\calK.
\end{equation}
The point mapping property of the Gromov--Witten invariants
shows that the small quantum product degenerates to the ordinary cup
product at $q=0$.

\begin{exa}\label{exa-Hirzebruch}
Consider the first Hirzebruch surface $F_1$ which is induced by the
following fan (left); on the right is shown the space
$H^2(F_1;\RR)$ using the coordinate system given by the classes of
$D_1$ and $D_2$. (See the start of Subsection \ref{subsec-LGmodels}
for information on how to view $H^2(X;\ZZ)$).

\begin{center}
\newdimen\scale
\newdimen\offsetx
\newdimen\offsety
\scale=0.9cm
\offsetx=6cm
\offsety=-1.0cm
\begin{tikzpicture}

\node[] at (0,-2.89) {Fan of Hirzebruch $F_1$};

\filldraw[color=blue!60] (\scale*0,\scale*0) -- (\scale*0,\scale*2) -- (\scale*2,\scale*2) -- (\scale*2,\scale*0) -- (\scale*0,\scale*0);
\filldraw[color=blue!30] (\scale*0,\scale*0) -- (\scale*0,\scale*2) -- (\scale*-2,\scale*2) -- (\scale*-2,\scale*-2) -- (\scale*0,\scale*0);
\filldraw[color=blue!50] (\scale*0,\scale*0) -- (\scale*-2,\scale*-2) -- (\scale*0,\scale*-2)  -- (\scale*0,\scale*0);
\filldraw[color=blue!40] (\scale*0,\scale*0) -- (\scale*0,\scale*-2) -- (\scale*2,\scale*-2)  -- (\scale*2,\scale*0) -- (\scale*0,\scale*0);
\draw(\scale*0,\scale*0)--(\scale*2,\scale*0);
\draw(\scale*0,\scale*0)--(\scale*0,\scale*2);
\draw(\scale*0,\scale*0)--(\scale*-2,\scale*-2);
\draw(\scale*0,\scale*0)--(\scale*0,\scale*-2);

\filldraw[cyan,opacity=.2] (\offsetx+\scale*0,\offsety+\scale*0) -- (\offsetx+\scale*0,\offsety+\scale*2) -- (\offsetx+\scale*2,\offsety+\scale*2) -- (\offsetx+\scale*2,\offsety+\scale*0) -- (\offsetx+\scale*0,\offsety+\scale*0);

\draw(\offsetx+\scale*0,\offsety+\scale*0)--(\offsetx+\scale*2,\offsety+\scale*0);
\draw[->,line width=1.8pt](\offsetx+\scale*0,\offsety+\scale*0) -- (\offsetx+\scale*1.5,\offsety+\scale*0);
\draw(\offsetx+\scale*0,\offsety+\scale*0)--(\offsetx+\scale*0,\offsety+\scale*2);
\draw[->,line width=1.8pt](\offsetx+\scale*0,\offsety+\scale*0) -- (\offsetx+\scale*0,\offsety+\scale*1.5);
\draw(\offsetx+\scale*0,\offsety+\scale*0)--(\offsetx+\scale*-2,\offsety+\scale*2);
\draw[->,line width=1.8pt](\offsetx+\scale*0,\offsety+\scale*0) -- (\offsetx+\scale*-1.5,\offsety+\scale*1.5);

\node[] at (2.2,0) {$D_1$};
\node[] at (0,2.2) {$D_2$};
\node[] at (-2.1,-2.1) {$D_3$};
\node[] at (0,-2.2) {$D_4$};
\node[] at (\offsetx+ -1.5cm,\offsety + 0.8cm) {$[D_4]$};
\node[] at (\offsetx+ 0.5cm,\offsety + 1.0cm) {$[D_2]$};
\node[] at (\offsetx+ 1.3cm,\offsety + -0.5cm) {$[D_1],[D_3]$};
\filldraw[cyan,opacity=.2] (4,-2.7) -- (4.6,-2.7) -- (4.6,-3) -- (4,-3) -- (4,-2.7);
\node[] at (7.15,-2.84) {nef cone of $F_1$ inside $H^2(F_1;\RR)$};
\end{tikzpicture}
\end{center}
We choose the homogeneous basis $T_0=1,\; T_1=[D_1],\; T_2=[D_2],\; T_3 =\PD(\{pt\})$.
The small quantum cohomology product of $F_1$ is determined by
\[
\begin{array}{llll}
  T_1 * T_0 = T_1,&  T_1 * T_1 = -q_1 T_1 + q_1 T_2, & T_1 * T_2 = T_3,                    & T_1 * T_3 = q_1q_2T_0 \\\\
  T_2 * T_0 = T_2,&  T_2 * T_1 =  T_3,              & T_2 * T_2 = q_2T_0 + T_3,&  T_2 * T_3= q_2 T_1+ q_1 q_2 T_0
\end{array}
\]
since one can conclude that
\[
T_3 * T_3 = T_3*(T_1 *T_2) = (T_3 * T_1) *T_2 = q_1q_2 T_0 *T_2 = q_1q_2 T_2.
\]
The small quantum cohomology ring of $F_1$ is therefore given by
\[
\CC[q_1,q_2,T_1,T_2]/\left(T_1^2 + q_1 T_1 -q_2 T_2 , T_2^2- T_1 T_2 -q_2, T_1 T_2^2 -q_2 T_1 - q_1q_2 \right).
\]
Restricting this ring to $q_1 =q_2 = 0$ gives $\CC[T_1,T_2]/ (T_1^2,
T_2^2- T_1T_2, T_1 T_2^2)$ which is isomorphic to the cohomology ring
(cf. \cite[Section 5.2]{Fulton}),
\[
H^*(F_1;\CC) \equiv \CC[D_1,D_2,D_3, D_4]/(D_1 D_3, D_2D_4, D_1D_2 D_4, D_1 - D_3, D_2 -D_3 -D_4)
\]
under the map $T_1 \mapsto D_1, T_2 \mapsto D_2$.
\schluss\end{exa}

We are going to give a reformulation of the quantum cohomology algebra in terms of certain differential systems. The intrinsic reason of the appearance of differential equations in this context is best understood when studying the \emph{big quantum product} instead of the small one as we have done above. It basically means to have a product on $H^*(X;\CC)$ which is parameterized by any class $\delta\in H^*(X;\CC)$
instead of a class in $H^2(X;\CC)$ (more precisely, instead of a representative of a coset in $\calK$). One can show that the structure constants of the big quantum product can be obtained as third derivatives of a generating function, referred to as the Gromov-Witten potential. This fact reveals an intrinsic integrability property of the (big) quantum product. Moreover, the associativity then boils  down to a famous third order non-linear partial differential equation satisfied by the GW-potential, abbreviated as WDVV-equation (after Witten, Dijkgraaf, Verlinde, Verlinde, see, \emph{e.g.} \cite{ManinFrob}). It turns out that using the next definition, this equation can be rewritten as a flatness property of a system of \emph{linear} differential equations, that is, a vector bundle with a connection.
\begin{dfn}
The \emph{small Dubrovin connection} $(H^\sideA,\nabla^\sideA)$ of $X$ is a flat
meromorphic connection $\nabla^\sideA$ on a trivial, holomorphic vector
bundle $H^\sideA$ over $\PP^1 \times \overline{\calK}$ with fiber
$H^*(X;\CC)$. The connection is given by
\begin{eqnarray}\label{eq:DubrovinConn}
  \nabla^\sideA_{\del_{q_i}} (T_j) &:= &\frac{1}{z} T_i \ast T_j \\
\nabla^\sideA_{\del_z} (T_j) &:= &-\frac{1}{z^2} c_1(X) \ast T_j + \frac{1}{z} \frac{\deg(T_j)}{2} T_j
\end{eqnarray}
where we denote by $z$ the coordinate centered at $0\in \CC\subseteq \PP^1$.
\schluss\end{dfn}
Notice however that this convention from quantum cohomology literature leads to some slight clash of notation.
Namely, the variable $z$ from above (a coordinate on $\PP^1$) is different from the variable $z$ used for univariate hypergeometric equations in Section \ref{sec:Introduction} as well as in Formula \eqref{eq:OneDimHyperOp}. In order to be consistent with the literature, we stick to these conventions and hope that it does not lead to confusion.

It is an easy but instructive exercise to check that the flatness of the connection $\nabla^\sideA$ implies the associativity and commutativity of the small quantum product.

\begin{exa}
The small Dubrovin connection of the first Hirzebruch surface is given by
\begin{align*}
\nabla^\sideA = \de &+ \left(\begin{matrix}0 & 0 & 0 & q_1q_2 \\ 1 & -q_1 &
  0 & 0 \\ 0 & q_1 & 0 & 0 \\ 0 & 0 & 1 & 0 \end{matrix} \right)
\frac{\de q_1}{zq_1} +
\left(\begin{matrix}0 & 0 & q_2 & q_1q_2 \\ 0 & 0 & 0 & q_2 \\ 1 & 0 &
  0 & 0 \\ 0 & 1 & 1 & 0 \end{matrix} \right) \frac{\de q_2}{zq_2}\\
&+
\left(\begin{matrix}0 & 0 & -2q_2 & -3q_1q_2 \\ -1 & q_1 & 0 & -2q_2
  \\ -2 & -q_1 & 0 & 0 \\ 0 & -2 & -3 & 0 \end{matrix} \right)
\frac{\de z}{z^2} + \left(\begin{matrix} 0 & 0 & 0 & 0 \\ 0 & 1 & 0 & 0
  \\ 0 & 0 & 1 & 0 \\ 0 & 0 & 0 & 2 \end{matrix} \right) \frac{\de z}{z}
\end{align*}
\schluss\end{exa}

\subsection{Landau--Ginzburg models}\label{subsec-LGmodels}

Let $\Sigma_X$ be the fan of the toric smooth projective Fano variety $X$
defined on the $d$-dimensional vector space $N \otimes_\ZZ \RR$ ($N\cong \ZZ^d$ being a lattice),
with $\Sigma_X(1)$ the set of one-dimensional cones whose primitive
elements in $N$ form the columns of the matrix $B\in\ZZ^{d\times n}$.
Denote by $M =
\Hom_\ZZ(N,\ZZ)$ the dual of $N$ which is identified with the group of
torus-invariant principal divisors and by $\Div_T(X)$ the group of
torus-invariant Weil divisors. There is the following
(split) exact sequence
\begin{equation}\label{eq-exseqtoric}
0 \longrightarrow M \longrightarrow  \Div_T(X) \longrightarrow  H^2(X,\ZZ) \longrightarrow 0
\end{equation}
Applying $(-) \otimes_\ZZ \CC^*$ one obtains the (split) exact sequence
\[
1 \longrightarrow \underbrace{M \otimes_\ZZ \CC^*}_{=\overline \TT}
\overset{\boldb}\longrightarrow
\Div_T(X) \otimes_\ZZ \CC^* \overset{\boldc}\longrightarrow
\underbrace{H^2(X,\ZZ)\otimes_\ZZ \CC^*}_{=\calK} \longrightarrow 1
\]
of algebraic tori, where $\boldb$ is the monomial map encoded by the
transpose of $B$, $\calK$ is as in Subsection \ref{subsec-GW+Dub}, and
$\overline \TT$ as in \eqref{eq:Embedding_g}. Recall that the
standard basis $e_1,\ldots,e_{d}$ of $M$ gives coordinates
$\boldt=(t_1,\ldots,t_{d})$ on $\overline{\TT}$.

The canonical basis of torus-invariant
divisors $D_1,\ldots,D_n$ for $\Div_T(X)$
corresponding to the one-dimensional cones induces
an isomorphism $\Div_T(X) \otimes_\ZZ\CC^* \simeq
(\CC^*)^n$.
Let $W: \Div_T(X) \otimes_\ZZ \CC^* = (\CC^*)^n \to \CC$ be the
function given by summing the coordinates.

\begin{dfn}\label{dfn-LGmodelK}
The Landau--Ginzburg model associated to the smooth, toric, Fano variety $X$ is the map
\[
(W,\boldc): \Div_T(X)\otimes_\ZZ \CC^* \to \CC \times \calK.
\]
\schluss\end{dfn}

If we view $\calK$ as an abstract algebraic torus, defining
the morphism $(W,\boldc)$ requires only the matrix $B$ (that is, the
generators of $\Sigma_X(1)$), but not the full data of
the fan $\Sigma_X$.  We shall later wish to (partially) compactify
$\calK$, as we have done before (see Formula \eqref{eq:CompactK}). For this, we need to equip $\calK$ with the coordinate system
$\{q_i\}_{i=1,\ldots,r}$, corresponding to the basis
$\{T_i\}_{i=1,\ldots,r}$ on $H^2(X;\CC)$. The compactification is
designed to contain the point $q_1=\ldots=q_r=0$, since there the quantum product collapses to the cup product.
This will be the case if the basis $\{T_i\}_{i=1,\ldots,r}$ of $H^2(X;\RR)$ consists of nef classes (this choice has already been made above at the beginning of Subsection \ref{subsec-GW+Dub}). Hence, fixing such a good coordinate system $\{q_i\}_{i=1,\ldots,r}$ on $\calK$ depends on the geometry of the toric variety $X_\Sigma$ and not just on the ray generators given by the matrix $B$ (see  \cite[Section 3.1]{ReiSe} for a more detailed discussion).

Since \eqref{eq-exseqtoric} splits,  we can find a section of the map
$\Div_T(X) \to H^2(X,\ZZ)$ which then induces a section
\begin{gather}\label{eqn-smap}
\bolds\colon \calK \to \Div_T(X) \otimes_\ZZ \CC^*.
\end{gather}
Again, $\bolds$, seen as a monomial map from $(\CC^*)^r$ to $(\CC^*)^n$, will depend on the fan structure of $\Sigma_X$ via the choice of coordinates on $\calK$. From now on, we will always fix such  coordinates and consider $\calK$ as the concrete $r$-dimensional torus $(\CC^*)^r$. The isomorphism
\[
(\boldb,\bolds):  \overline{\TT} \times \calK \to \Div_T(X) \otimes_\ZZ \CC^*
\]
gives a different presentation of the Landau--Ginzburg model, namely
as a family of Laurent polynomials
\begin{align}\label{eqn-psi}
\psi := (F,\pr_2): \overline{\TT} \times \calK &\to \CC \times \calK \\
(\pointt_1,\ldots, \pointt_{d}, \pointq_1, \ldots, \pointq_r) &\mapsto (\sum_{j=1}^n \pointq^{\bolds_j} \pointt^{\boldb_j},\pointq_1,\ldots,\pointq_r)\nonumber
\end{align}
where $S=(\bolds_1,\ldots,\bolds_n)\in\ZZ^{r\times n}$ and
$B=(\boldb_1,\ldots,\boldb_n)\in\ZZ^{d\times n}$ represent the maps
$\bolds$ and $\boldb$ respectively.

\begin{exa}
  We continue Example \ref{exa-Hirzebruch}.
The exact sequence \eqref{eq-exseqtoric} is given by
\setlength\arraycolsep{2pt}
\[
\begin{tikzcd}[ampersand replacement=\&,]
0 \to \ZZ^2 \arrow[rr,"{{\left(\tiny\begin{matrix} 1 & 0 \\ 0 & 1 \\ -1
      & -1 \\ 0 & -1\end{matrix} \right)}}"]\&\&
\ZZ^4 \arrow[rr,"{{\tiny \left(\begin{matrix} 1& 0 & 1 & -1 \\ 0 & 1 &
        0 & 1 \end{matrix} \right)}}"]\&\&
\ZZ^2 \to 0
\end{tikzcd}
\]
where we have chosen the basis $T_1 =[D_1], T_2 = [D_2]$ as a basis in
$H^2(X;\ZZ)$, as we did in Example \ref{exa-Hirzebruch}.  The
Landau--Ginzburg model
is given on the level of coordinate functions by
\begin{eqnarray*}
(W,\boldc)\colon \Div_T(X)\otimes_\ZZ\CC^*=(\CC^*)^4 &\longrightarrow&
  \CC\times (\CC^*)^2=\CC\times \calK\\
 (\pointw_1+\cdots+\pointw_4,\frac{\pointw_1 \pointw_3}{\pointw_4},\pointw_2 \pointw_4) &\mapsfrom& (\pointt,\pointq_1,\pointq_2).
\end{eqnarray*}
The corresponding family of Laurent polynomials is
\begin{align*}
\psi\colon \overline \TT\times\calK=(\CC^*)^2 \times (\CC^*)^2 &\longrightarrow \CC \times (\CC^*)^2=\CC\times\calK\\
(\pointt_1,\pointt_2,\pointq_1,\pointq_2) &\mapsto (\pointq_1 \pointt_1 + \pointq_2 \pointt_2 + \frac{1}{\pointt_1 \pointt_2} + \frac{1}{\pointt_2}, \pointq_1, \pointq_2),
\end{align*}
where we have chosen the section $\bolds: \calK \to \Div_T(X)
\otimes_\ZZ \CC^*$ as the one induced from the map
\[
\begin{tikzcd}[ampersand replacement=\&,]
 H^2(X;\ZZ)\cong \ZZ^2 \arrow[rr,"{{\left(\tiny\begin{matrix} 1 & 0 \\ 0 & 1 \\ 0 & 0 \\ 0 &0\end{matrix} \right)}}"]\&\&
\ZZ^4 \cong \Div_T(X).
\end{tikzcd}
\]
\schluss\end{exa}

\medskip

It was conjectured by Givental (see, \emph{e.g.} \cite{Giv7}) that
oscillating integrals over Lefschetz thimbles  with respect to the
Landau--Ginzburg model give flat sections of the Dubrovin
connection. An algebraic replacement of these oscillating integrals, localized and  partially Fourier--Laplace transformed Gau\ss--Manin systems
of the Landau--Ginzburg model.

We briefly explain this version of the ordinary Fourier--Laplace transformation functor (see Formula \eqref{eq:FL} above).
In the following,
$\calO_{\CC_t\times\CC_\tau\times Y}\cdot \exp({-t\tau})$ denotes a
free rank $1$ module with twisted differential given by the product
rule.
\begin{dfn}
  Given a smooth variety $Y$ and a holonomic $\calD_{\CC\times
    Y}$-module $\calN$, the \emph{localized, partial Fourier--Laplace
    transform} of $\calN$ is the sheaf
\begin{equation}\label{eq:FLloc}
\FL_Y^{\loc} \calN := (j_z)_+ j_{\tau}^+ (p_2)_+ \left(p_1^+ \calN \otimes_\calO \calO_{\CC_t\times\CC_\tau\times Y}\cdot \exp({-t\tau}) \right)[-1]
\end{equation}
where $p_1: \CC_t \times \CC_\tau \times Y \to \CC_t\times Y $
and $p_2: \CC_t \times \CC_\tau \times Y \to  \CC_\tau\times Y $ are
the indicated projections, and where $j_\tau : \CC^*_\tau \times Y \to
\CC_\tau\times Y $ and $j_z: \CC^*_\tau \times Y \to
(\PP^1_\tau\backslash\{0\})\times Y =\CC_z\times  Y$ are the canonical
open embeddings with the understanding that $z = 1/\tau$.
\schluss\end{dfn}
The name ``localized'' comes  from the fact that by using the direct image $(j_z)_+$, the action of $z$ is invertible on the resulting module (and so is the action of $\tau$).

The localized, partially Fourier--Laplace transformed Gau\ss--Manin system of the Landau--Ginzburg model $\psi$ is then defined as
\[
\calG^\psi := \FL^{\loc}_\calK \calH^0(\psi_+ \calO_{\overline{\TT} \times \calK}).
\]
It is an exercise (using the definition of the direct image functor, see, \emph{e.g.} \cite[Sections 1.3, 1.5]{HTT}) to show that the module of global sections $G^\psi$ of $\calG^\psi$ has the following presentation in terms of relative differential forms
\[
G^\psi \simeq H^0 \left( \Omega^{\bullet+d}_{\overline{\TT} \times
  \calK / \calK}[z^\pm], z\frakd- \de F \wedge \right),
\]
where $\frakd$ is the differential on the complex $\Omega^{\bullet+d}_{\overline{\TT} \times \calK / \calK}$.
Following an idea from singularity theory (see \cite{Brie,SM,Sa2}), one defines the \emph{Fourier-Laplace transformed Brieskorn lattice} by
\begin{eqnarray}\label{eqn-FLBrieskorn}
G_0^\psi &:=& H^0 \left( \Omega^{\bullet+d}_{\overline{\TT} \times \calK
  / \calK}[z], z\frakd- \de F \wedge \right)
\subseteq G^\psi.
\end{eqnarray}
We will see below, using GKZ-systems, that $G^\psi_0$ is $\calO_{\CC \times \calK}$-free.
In order to connect $\calG^\psi$ to a GKZ-system we observe that the family of Laurent polynomials $\psi$ is a pullback of a larger family
\begin{align*}
\varphi\colon \overline{\TT} \times \CC^n &\to \CC \times \CC^n\\
((\pointt_1,\ldots,\pointt_d),(\pointx_1,\ldots,\pointx_n)) &\mapsto (-\sum_{j=1}^n \pointx_j \pointt^{\boldb_j},(\pointx_1,\ldots,\pointx_n))
\end{align*}
by the map
\begin{align}\label{eqn-iota}
\begin{tikzcd}[ampersand replacement=\&,]                                             \iota\colon \CC \times \calK \arrow[hook,rr,"\id\times(-\bolds)"] \&\&
  \CC\times \Div_T(X)\otimes_\ZZ\CC^* \arrow[r,"\simeq"] \&
  \CC\times (\CC^*)^n \arrow[hook,r,"\can"] \&
  \CC \times \CC^n
\end{tikzcd}
\end{align}
where
$\bolds:\calK\hookrightarrow \Div_T(X)\otimes_\ZZ\CC^*\cong (\CC^*)^n$
is as in \eqref{eqn-smap} and the middle map is the identification
induced from the standard basis on $M$.

In Theorem \ref{theo:4TermSeq} we have connected
the Gau\ss--Manin  system of $\varphi$ to a GKZ system via the 4-term sequence
\[
0 \rightarrow H^{d-1}(\overline{\TT};\CC) \otimes_\calO \calO_{\CC^{n+1}} \rightarrow \calH^0(\varphi_+ \calO_{\overline{\TT} \times \CC^n}) \rightarrow \calM_A(0) \rightarrow H^{d}(\overline{\TT};\CC) \otimes_\calO \calO_{\CC^{n+1}} \rightarrow 0,
\]
where $A\in \ZZ^{(d+1)\times(n+1)}$ is the homogenization of the
matrix $B$ constructed from the ray generators of the fan $\Sigma_X$.
Since the outer two terms are free $\calO_{\CC^{n+1}}$-modules, they
are in the kernel of the localized partial Fourier--Laplace transform. Indeed,
on the level of global sections, $\FL^{\loc}_Y$ is the
composition the localization at $\del_t$
with the
ordinary Fourier--Laplace transformation $\FL_Y$, and
$\CC[t]=D_t/D_t\cdot \del_t$ naturally localizes to zero.
Thus, the localized partial Fourier--Laplace transform being the composition of two
exact functors, the previous display implies
\[
\calG^\varphi = \FL^{\loc}_{\CC^n}\calH^0(\varphi_+ \calO_{\TT\times \CC ^{n-1}}) \simeq \FL^{\loc}_{\CC^n} (\calM_A(0)).
\]
The module of global sections of $\FL^{\loc}_{\CC^n} (\calM_A(0))$
is the cyclic left module $D_{\CC \times \CC^n}[z^\pm]/I$ over the
ring
\[
D_{\CC \times \CC^n}[z^\pm]:=\CC[z^\pm,x_1,\ldots,x_n]\langle
\del_z,\del_{x_1}, \ldots, \del_{x_n}\rangle,
\]
where $I$ is generated by the operators $\widehat{E}_0,
(\widehat{E}_i)_{i = 1,\ldots,d}, (\widehat{\Box}_{\boldu})_{\boldu
  \in \ker(B)}$ from Equation \eqref{eq:zGKZ-ideal}.
We like to compare this computation to a presentation for the
Fourier-Laplace transformed Brieskorn lattice
$\calG^\varphi_0\subseteq\calG^\varphi$ for the map $\varphi$ instead of
$\psi$.  For this, we use again the Rees ring $R_{\CC \times \CC^n} =
\CC[z,x_1,\ldots,x_n]\langle z^2 \del_z, z\del_{x_1},\ldots, z
\del_{x_n}\rangle$ from Equation \eqref{eq-Rees}. The module of global
sections of the Fourier-Laplace transformed Brieskorn lattice
$G^\varphi_0$ can then be described as $R_{\CC \times \CC^n}/
H_B^z(0)$, recalling from Section \ref{sec:Hodge} that $H_B^z(0)$ is the
left $R_{\CC\times \CC^n}$-ideal generated by the operators
$\widehat{E}_0, (\widehat{E}_i)_{i = 1,\ldots,d},
(\widehat{\Box}_{\boldu})_{\boldu \in \ker(B)}$.

Using techniques borrowed from \cite{Adolphson-duke94} one can show:
\begin{lem}\cite[Lemma 2.12]{ReiSe}
The restriction of the Fourier-Laplace transformed Brieskorn lattice
$\calG^\varphi_0$ to the Zariski open subset $\CC\times
(\CC^*)^n\subseteq \CC\times \CC^n$ is a free $\calO_{\CC \times
  (\CC^*)^n}$-module. \qed
\end{lem}

One can prove by base change that the Fourier-Laplace transformed
Brieskorn lattice $\calG^\varphi_0$ is the inverse image of
$\calG^\psi_0$ under the map $\iota$ in \eqref{eqn-iota}. We therefore
arrive at the following result where, for $\boldu\in\ker(B)$, we read
it as an element of $H_2(X;\CC)$ via the dual of the sequence \eqref{eq-exseqtoric}:

Parallel to $\shRees_{\CC\times\CC^n}$ from \eqref{eq-Rees}, we define
\[
R_{\CC \times {\calK}} :=
\CC[z,q^\pm_1,\ldots,q^\pm_r]\langle z^2\del_z, z \del_{q_1}, \ldots, z
\del_{q_r} \rangle
\]
and denote by $\shRees_{\CC \times
  {\calK}}$ the associated sheaf on $\CC \times
{\calK}$.
\begin{prp}\label{prp-locBrieskorn}
The localized Fourier-Laplace transformed Brieskorn lattice
$\calG^\psi_0$ is $\calO_{\CC\times\calK}$-free.
As a sheaf over $\shRees_{\CC \times \calK}$, it is
isomorphic to the cyclic module $\shRees_{\CC \times \calK}/ \calJ$
where the left ideal $\calJ$ is generated by (here, $\boldu$ runs
through $\ker(B)$ and $\{q_a\}_{a=1\ldots,r}$ are coordinates on $\calK$ as always)
\begin{eqnarray*}
\widetilde{E}&:=& z^2\del_z + \sum_{a=1}^r c_1(X)_a zq_a \del_{q_a}\\
\widetilde{\Box}_\boldu &:= & \!\!\left(\prod_{a : T_a(\boldu) > 0}\!\! q_a^{T_a(\boldu)}\right)\!\! \prod_{j: u_j < 0}\!\! \prod_{\nu = 0}^{-u_j-1}(\sum_{a=1}^r [D_i]_a zq_a \del_{q_a}- \nu z)\\
&&-   \left(\prod_{a : T_a(\boldu) < 0}\!\! q_a^{-T_a(\boldu)}\right)\!\! \prod_{j: u_j > 0} \prod_{v=0}^{u_j-1} (\sum_{a=1}^r [D_i]_a z q_a \del_{q_a} -\nu z)
\end{eqnarray*}
where $[D_i] = \sum_{a=1}^r [D_i]_a T_a$ and $c_1(X) = \sum_{a=1}^r
c_1(X)_a T_a$.
\end{prp}

Set
\[
R^{\log}_{\CC \times \overline{\calK}} :=
\CC[z,q_1,\ldots,q_r]\langle z^2\del_z, zq_1 \del_{q_1}, \ldots, z q_r
\del_{q_r} \rangle
\]
and denote by $\shRees^{\log}_{\CC \times
  \overline{\calK}}$ the associated sheaf on $\CC \times
\overline{\calK}$. Then the following statements on some cyclic
$\shRees^{\log}_{\CC \times \overline{\calK}}$-modules are proved in
\cite{ReiSe} using methods from toric geometry, including the notions of
primitive collections and relations (see, \emph{e.g.},
\cite{CoxRenesse,CoxBook}).
\begin{prp}
Let $\calJ^{\log} \subseteq \shRees^{\log}$ be the left ideal generated by
$\widetilde{E}$ and $\widetilde{\Box}_\boldu$ from Proposition \ref{prp-locBrieskorn}. Then
\begin{itemize}
\item $\shRees^{\log}_{\CC \times \overline{\calK}}/\calJ^{\log}$ is $\calO_{\CC \times \overline{\calK}}$-free.
\item $(\shRees^{\log}_{\CC \times \overline{\calK}}/\calJ^{\log})_{\mid \CC \times \calK} \simeq \shRees_{\CC \times \calK}/\calJ$.
 \end{itemize}
\end{prp}

In order to construct an object which matches the small Dubrovin
connection coming from the Gromov--Witten invariants of $X$ we have to
go one step further. Recall that the small Dubrovin connection
\eqref{eq:DubrovinConn} is a family of vector bundles on $\PP^1$, parameterized by $\overline{\calK}$, equipped with
a certain connection operator. As of yet, starting from the
Landau--Ginzburg model $\psi$ from \eqref{eqn-psi} of $X$, we have
constructed a vector bundle $\shRees^{\log}_{\CC \times
  \overline{\calK}}/\calJ^{\log}$ on $\CC\times \overline{\calK}$ with
a differential structure, and it is easily verified that the behavior
along the poles $(\{0\}\times \overline{\calK})\cup(\CC\times
(\overline{\calK}\backslash \calK))$ of the connection operators on
both bundles
are of the same type. If we want to compare $\shRees^{\log}_{\CC
  \times \overline{\calK}}/\calJ^{\log}$ to the small Dubrovin connection, it thus remains to extend this bundle (together with its connection operator) over the divisor
$\{\infty\}\times \overline \calK$ to all of $\PP^1\times \overline\calK$. This is of course
always possible if no other condition is imposed. However, if we want to reconstruct the Dubrovin connection, this extension
needs to satisfy two strong conditions simultaneously: the resulting
object must be a family of trivial $\PP^1$-bundles \emph{and} the
connection must have a logarithmic pole at infinity. Fulfilling both
requirements is not always possible, and goes under the name
(Riemann-Hilbert-)Birkhoff problem; for a modern account see
\cite[Chapter IV]{Sa3}. However, under the current circumstances, a
solution to the Birkhoff problem can be found locally near the
boundary $\overline{\calK}\backslash\calK$, as the following result
shows.

\begin{thm}\label{theo:SolBirkhoff}(\cite[Proposition 3.10]{ReiSe})
There exists a Zariski open neighborhood $U$ of $0 \in
\overline{\calK}$ and sections $Q_0,\ldots, Q_s$ of
$(\shRees^{\log}_{\CC \times \overline{\calK}}/\calJ^{\log})_{\mid \CC
  \times U}$ which extend $(\shRees^{\log}_{\CC \times
  \overline{\calK}}/\calJ^{\log})_{\mid \CC \times U}$ as a (trivial)
holomorphic vector bundle over $\PP^1 \times U$, called $H^\sideB$, such
that the associated connection $\nabla^\sideB$ has a logarithmic pole along
the normal crossing divisor $(\{z=\infty\} \times U) \cup (\PP^1_z \times
(\overline{\calK} \setminus \calK))$.
\end{thm}

With all these preparations, we can state the following result, which can be considered as
the Hodge theoretic mirror statement for smooth toric Fano varieties.
\begin{thm}\label{theo:MirrorCompactCase}(\cite[Proposition 4.10]{ReiSe})
Let, as before, $X$ be a smooth projective toric Fano variety, $(H^\sideA,\nabla^\sideA)$ the small Dubrovin connection and $(H^\sideB,\nabla^\sideB)$ the solution to the Birkhoff problem from Theorem \ref{theo:SolBirkhoff}.
Then there is an isomorphism of holomorphic bundles over $\PP^1\times U$ with meromorphic connections
\[
(H^\sideA,\nabla^\sideA)_{\mid \PP^1 \times  U} \simeq (H^\sideB,\nabla^\sideB).
\]
\end{thm}
We remark that in \cite[Proposition 4.10]{ReiSe} a similar result
for the more general case of weak Fano toric manifolds is given, albeit with a weaker
conclusion:
the extension $H^\sideB$ there
only exists on an analytic open subset of $\calK$ (see the remark
after \cite[Proposition 3.10]{ReiSe}).

\begin{exa}
When $X$ is the Hirzebruch $F_1$ surface, the Fourier-Laplace transformed
Brieskorn lattice of the Landau--Ginzburg model is given by
\[
G^\psi_0 \simeq \CC[z,q^\pm_1,q^\pm_2]\langle z^2\del_z, z\del_{q_1}, z\del_{q_2}\rangle/ J
\]
where the left ideal $J$ is generated by the operators
\begin{align*}
  \tilde{E}&= z^2 \del_z + z q_1 \del_{q_1} + 2 z q_2 \del_{q_2},
  \displaystyle &&
\widetilde{\Box}_{\boldu_1} = (z q_1 \del_{q_1})^2 +q_1 (z q_1 \del_{q_1}) - q_1 (z q_2 \del_{q_2}),\\
\displaystyle
\widetilde{\Box}_{\boldu_2} &= (z q_1 \del_{q_1})^2(z q_2 \del_{q_2}) - q_1 q_2,&&
\displaystyle
\widetilde{\Box}_{\boldu_3} = -(z q_1 \del_{q_1})(z q_2 \del_{q_2}) + (z q_2 \del_{q_2})^2 - q_2,
\end{align*}
where $\boldu_1=(1,0,1,-1)$, $\boldu_2=(1,1,1,0)$,
$\boldu_3=(0,1,0,1)$ generate the integer kernel of $B$.

The logarithmic extension is equal to $\CC[z,q_1,q_2]\langle z^2
\del_z, z q_1 \del_{q_1}, z q_2 \del_{q_2}\rangle / J^{\log}$ where
$J^{\log}$ is  generated by the same operators as $J$.

The basis which solves the (Riemann-Hilbert)-Birkhoff problem is $Q_0
= 1, Q_1 = z q_1 \del q_1, Q_2 = z q_2 \del_{q_2}, Q_3 = (z q_1
\del_{q_1})(z q_2 \del_{q_2})$. These sections are identified with the
sections $T_0, T_1, T_2, T_3$  of $H^\sideA$ under the mirror isomorphism from Theorem \ref{theo:MirrorCompactCase}.
\schluss\end{exa}

\subsection{Reduced quantum $\calD$-modules and intersection
  cohomology}

In this section, we are going to discuss a mirror statement that
concerns weak Fano smooth complete intersections inside smooth
projective toric, possibly non-Fano, varieties.
From the point of view of physics, this
is an even more important class of examples than the one considered
previously since it includes Calabi--Yau manifolds that are
subvarieties of toric manifolds, although they are not toric
themselves. The most prominent example, namely, the quintic in $\PP^4$
(where the first enumerative predictions using the mirror symmetry
principle were made, see \cite{Candelas}) is of this type. We will
discuss a non-affine version of the Landau--Ginzburg models introduced
above. The mirror statement that we aim for will relate (part of) the
quantum cohomology of the complete intersection subvariety to the
lowest weight filtration step of a GKZ-system. It follows from the
results in Section \ref{sec:HodgeWeight} that the lowest weight
filtration step is a single intersection cohomology $\calD$-module
which arises as the image under a natural morphism from the holonomic
dual of the GKZ system to the GKZ system itself. In the cases we
discuss here this holonomic dual is isomorphic to a GKZ system with
the same matrix $A$ but different parameter vector $\beta$. Hence the
intersection cohomology $\calD$-module can be described as the image
of a morphism between two GKZ-systems by a contiguity morphism. Our
main reference in this section is \cite{ReiSe2}. We start with setting
the notation.

\begin{ntn}\label{ntn-QDM+IC}
  As before, $X$ will be a smooth
projective toric variety of Picard rank $r$ attached to the fan $\Sigma_X$ of dimension
$d$, whose primitive rays form the columns of the matrix
$B$.
In contrast to the previous case  we do in this subsection not make any
positivity assumption on $X$ here. Let $\calO(L_1), \ldots,
\calO(L_c)$ be globally generated line bundles; since $X$ is toric,
this amounts to asking that each $L_i$ be nef---their classes should lie in the
nef cone in $H^2(X, \RR)$. We shall assume also that
\begin{equation}
\label{eq:NefCondition}
-K_X-L_1-\ldots-L_c\quad\textup{is nef}.
\end{equation}
If $D_1,\ldots,D_n$ are the torus invariant divisors on $X$ we can
write
\begin{equation}\label{eq-dij}
L_j=\sum_{i=1}^nd_{ij}D_j
\end{equation}
for suitable non-negative integers $d_{ij}$.
Set
\[
\calE:=\calO(L_1)\oplus \ldots \oplus \calO(L_c),
\] and consider a
generic global section $\gamma\in \Gamma(X,\calE)$. Our assumptions
imply that
\[
Y := \gamma^{-1}(0)\subset X
\]
is a smooth complete intersection subvariety for which $-K_Y$ is nef; we
call this property \emph{weak Fano}.
\schluss\end{ntn}

In this paragraph we
briefly review a variant of the above quantum product that is designed
to encode enumerative information about stable maps to $Y$. The first
point is that one can generalize the definition of Gromov--Witten
invariants \eqref{eq:GWinv} to the \emph{twisted}
(three-point) \emph{GW-invariants}; these are also maps from
$H^*(X,\QQ)^{\otimes 3}\rightarrow \QQ$, but Chern classes of
certain tautological bundles (on the moduli space of stable maps)
derived from $\calE$ come into play. We denote by $\langle
I_{0,3,\beta}\rangle(\gamma_1,\gamma_2,\widetilde{\gamma}_3)\in\QQ$
the value of such a three point twisted GW-invariant for given
cohomology classes $\gamma_1,\gamma_2,\gamma_3\in H^*(X,\QQ)$  (see, \emph{e.g.}
\cite[Section 4.1]{ReiSe2} for a more detailed discussion, including an
explanation for the process $\gamma_3\leadsto\widetilde{\gamma}_3$). Then one
defines in complete analogy to Formula \eqref{eq-smallquantum} the
twisted (small) quantum product by
\begin{eqnarray}\label{eqn-*tw}
\gamma_1 \stackrel{\tw}{\ast}  \gamma_2 := \sum_{a=0}^s \sum_{\beta \in H_2(X, \ZZ)} q^\beta \langle I_{0,3,\beta}\rangle(\gamma_1, \gamma_2, \widetilde{T}_a) T^a\,,
\end{eqnarray}
where, as before, $q$ are coordinates on $\calK$ and
$q^\beta:=\exp({\delta(\beta)})$ for $\beta\in H_2(X;\CC)$.

We now follow the definition of the small Dubrovin connection,
Equation \eqref{eq:DubrovinConn}, and define the \emph{twisted quantum
  $\calD$-module}, denoted by $\QDM(X,\calE)$, as the vector bundle on
$\PP^1\times \calK$ with fiber $H^*(X;\CC)$ together with the
connection given by
\begin{eqnarray*}
\nabla^{\tw}_{\del_{q_i}} T_j &:= &\frac{1}{z} T_i \stackrel{\tw}{\ast} T_j
\\ \nabla^{\tw}_{z\del_z} T_j &:= &-\frac{1}{z}
(t_0T_0+c_1(X)-c_1(\calE)) \stackrel{\tw}{\ast} T_j +
\frac{\deg(T_j)-\dim(X)+\textup{rk}(\calE)}{2} T_j
\end{eqnarray*}
Notice that, unlike in the Fano case discussed in Subsection
\ref{subsec-LGmodels}, the
convergence of the twisted quantum product is not automatic.
We will therefore  later restrict to some analytic neighborhood $U\subset \calK$ of the
point $q_1=\ldots=q_r=0$ in $\overline{\calK}$, on which $\stackrel{\tw}{\ast}$ is convergent.

As we are interested in enumerative information about maps to $Y:=\gamma^{-1}(0)$,
the cohomology space $H^*(X;\CC)$ is not a well suited object for a
quantum cohomology theory of $Y$. We therefore consider the Gysin morphism
\begin{eqnarray*}
m_{\calE}\colon H^*(X) &\longrightarrow& H^*(X) \\
\alpha &\longrightarrow & c_{top}(\calE) \cup \alpha
\end{eqnarray*}
and define the \emph{reduced cohomology} of $(X,\calE)$ to be
\[
\overline{H^*(X)}:= H^*(X)/\ker(m_{\calE}).
\]
One checks that the twisted quantum
$\calD$-module $\QDM(X,\calE)$ has a quotient bundle
$\overline{\QDM}(X,\calE)$ with fiber $\overline{H^*(X)}$, and that
the connection $\nabla^{\tw}$ on $\QDM(X,\calE)$ descends to
$\overline{\QDM}(X,\calE)$.
We call
this  vector bundle on $\PP^1\times \calK$ with connection
$(\overline{\QDM}(X,\calE), \nabla^{\tw})$ the \emph{reduced quantum $\calD$-module} (see
\cite[Definition 4.3]{ReiSe2} for more details).

\medskip

We proceed by describing the relevant Landau--Ginzburg models
attached to the given data $(X,\calE)$. Denote by $\calE^\vee$ the
dual bundle of $\calE$, and by
\[
\VV:=\VV(\calE^\vee)\stackrel{\pi}{\longrightarrow} X
\]
its total
space. Then $\VV$ is a (non-compact) toric variety, whose fan

\[
\Sigma_\VV \subseteq (N\oplus \ZZ^c)\otimes_\ZZ \RR
\]
is given as
follows: The set of rays of $\Sigma_\VV$ are the columns of the matrix
\begin{equation}\label{eq:MatrixBprime}
B'=(\boldb'_1,\ldots,\boldb'_{n+c}):=
\left(
\begin{array}{c|c}
B & 0_{n,c} \\ \hline
(d_{ji}) & \textup{Id}_c
\end{array}
\right)\in \ZZ^{(d+c)\times(n+c)},
\end{equation}
where $B$ is the $d\times n$-matrix constructed from the primitive rays
in $\Sigma_X$ and
where $d_{ji}$ are as in \eqref{eq-dij}.
Then the fan  $\Sigma_\VV$ consists
of all cones
\[
\RR_{\geq 0} \boldb'_{i_1}+\ldots +\RR_{\geq 0} \boldb'_{i_k}+
\RR_{\geq 0}\boldb'_{j_1}+\ldots+\RR_{\geq 0}\boldb'_{j_\ell}
\]
such that $\RR_{\geq 0} \boldb_{i_1}+\ldots +\RR_{\geq 0}
\boldb_{i_k}\in \Sigma_X$ and  $j_1,\ldots, j_\ell\in
\{n+1,\ldots,n+c\}$. Notice that we have $H^2(\VV;\ZZ)\cong
H^2(X,\ZZ)\cong \ZZ^r$ and that $\Div_T(\VV)\cong \ZZ^{n+c}$.
Similarly to the discussion in Section \ref{subsec-LGmodels} we then
consider a family of Laurent polynomials associated to these toric
data.

\begin{dfn}(\cite[Definition 6.3.]{ReiSe2})\label{dfn-affineLG}
Let $(X,\calE)$ be as in Notation \ref{ntn-QDM+IC} and consider the
complexified K\"ahler moduli space $\calK \cong H^2(X;\ZZ)\otimes_\ZZ
\CC^*\cong H^2(\VV;\ZZ)\otimes_\ZZ\CC^*$ of both $X$ and
$\VV$. Write $\overline{\TT}_\VV:=(\CC^*)^{d+c}$ for the
  $(d+c)$-dimensional torus. Then the \emph{affine Landau--Ginzburg model} of
  $(X,\calE)$ is the morphism

\begin{eqnarray}\label{eqn-affLG}
\psi=(F,\pr_2)\colon \overline{\TT}_\VV\times
\calK^\circ&\longrightarrow &\CC\times  \calK^\circ \\
(\pointy,\pointq) &\longmapsto& \left(-\sum_{j=1}^{n} \pointq^{\bolds'_j} \cdot \pointy^{\boldb'_j}
+\sum_{j=n+1}^{n+c} \pointq^{\bolds'_i} \cdot \pointy^{\boldb'_i},\pointq\right),
\end{eqnarray}
where
\[
\calK^\circ \subseteq \calK
\]
is a Zariski open subset on which the Laurent polynomials
$\psi(-,\mathboldq)$  satisfy a non-degeneracy condition (see \cite[Section 3.2]{ReiSe2}) and
where $(\bolds'_1,\ldots,\bolds'_{n+c}) \in \ZZ^{r \times (n+c)}$ is a
section of the projection $\Div_T(\VV) \twoheadrightarrow
H^2(X,\ZZ)$.
\schluss\end{dfn}

One can establish a mirror symmetry theorem for the twisted quantum
$\calD$-module which involves the affine Landau--Ginzburg model, very
much in the same spirit (without looking at logarithmic extensions
over the boundary $\overline{\calK}\backslash \calK$ though, and also
neglecting the extension to families of bundles over $\PP^1$) as
Theorem \ref{theo:MirrorCompactCase} above (see \cite[Theorem 6.13,
  6.16]{ReiSe2} and also \cite{Mo15}). However, in order to
reconstruct the reduced quantum $\calD$-module
$\overline{\QDM(X,\calE)}$, we are forced to look at a
compactification of the morphism $\psi$. In order to define it,
consider the map $g_{B'}:
\overline\TT_\VV=(\CC^*)^{d+c}\hookrightarrow \PP^{n+c}$ (see Formula
\eqref{eq:Embedding_g} above). Then define
\begin{equation}\label{eqn-Zcirc}
  Z^\circ:=\overline{\Gamma}_F
\end{equation}
to be the closure in $\PP^{n+c}\times
\CC \times \calK^\circ$ of the graph $\Gamma_F \subseteq
\overline\TT_\VV\times \CC\times \calK$ of the function
$F:\overline{\TT}_\VV\times \calK^\circ \rightarrow \CC$ defined in
\eqref{eqn-affLG}. Notice that $Z^\circ$ is a partial compactification
of $\overline{\TT}_\VV\times \calK^\circ$, that is, quasi-projective but
in general not smooth.
\begin{dfn}\label{dfn-Psi}
Let $(X,\calE)$ be as above. Then we call the restriction
\[
\Psi: Z^\circ \longrightarrow \CC\times \calK^\circ
\]
of the projection
\[
\pr:\PP^{n+c}\times \CC \times \calK^\circ
\rightarrow\CC\times \calK^\circ
\]
the \emph{non-affine Landau--Ginzburg model} of $(X,\calE)$.
\schluss\end{dfn}
Clearly, $\Psi$ is a projective morphism, and hence should be
considered as a partial compactification
of the affine Landau--Ginzburg model $\psi$.

In a rather similar way to the case of Landau--Ginzburg models of projective
toric varieties, we obtain the following description of the relevant Gau\ss--Manin cohomologies
by GKZ-type systems. As a matter of notation, consider the
the matrix $A'\in \ZZ^{1+d+c,1+n+c}$ obtained by homogenizing the matrix $B'$ defined in
equation \eqref{eq:MatrixBprime}, that is
$$
A'=
\left(
\begin{array}{c|c}
1 & 1_{1,n+c} \\ \hline
0_{d+c,1} & B'
\end{array}
\right)
=
\left(
\begin{array}{c|c|c}
1 & 1_{1,n} & 1_{1,c}\\ \hline
0_{d,1} & B & 0_{n,c} \\ \hline
0_{c,1} & (d_{ji}) & \textup{Id}_c
\end{array}
\right).
$$
We choose the parameter vector
\[
\gamma:=(-c,\underbrace{0,\ldots,0}_{d
  \textup{ copies}},\underbrace{-1,\ldots,-1}_{c \textup{ copies}})\in \ZZ^{1+d+c}.
\]
With these definitions, we have the contiguity morphism (see Section \ref{subsec:Contiguity})
$$
\begin{tikzcd}
c_{\gamma,0}:\calM_{A'}(\gamma)
\ar{rr}{\partial_{n+1}\cdot\ldots\cdot\partial_{n+c}} &&
\calM_{A'}(0),
\end{tikzcd}
$$
due to the special shape of the matrix $A'$.
Notice that here we use the coordinates $(x_0,x_1,\ldots,x_{n+c})$ on $\CC\times\CC^{n+c}$ and $\partial_0,\partial_1,\ldots,\partial_{n+c}$ for the corresponding partials.

We can now formulate the following statement about the non-affine Landau-Ginzburg.
\begin{thm}[{\cite[Lemma 6.4 and Proposition 6.7]{ReiSe2}}]
There is an isomorphism of $\calD_{\CC\times \calK^\circ}$-modules
\[
\FL^{\loc}_{\calK^\circ} \calH^0\psi_+ \calO_{\overline{\TT}_\VV\times\calK^\circ} \cong \iota^+ \FL^{\loc}_{\CC^{n+c}}
\calM_{A'}(0)
\]
where we denote (with a slight abuse of notation) by $\iota: \CC\times
\calK^\circ\hookrightarrow \CC\times \CC^{n+c}$ the embedding already
used above (see Equation \eqref{eqn-iota}). Moreover, there is an
isomorphism of $\calD_{\CC\times \calK^\circ}$-modules
\[
\FL^{\loc}_{\calK^\circ} \calH^0 \pr_+ \IC(\underline{\CC}_{\overline{\TT}_\VV\times\calK^\circ}) \cong
\iota^+
\FL^{\loc}_{\CC^{m+c}} \image\left(c_{\gamma, 0}:\calM_{A'}(\gamma)\longrightarrow\calM_{A'}(0)\right).
\]
\end{thm}
Notice  that by
definition, the intersection cohomology module
$\IC(\underline{\CC}_{\overline{\TT}_\VV\times\calK^\circ})$ to the
constant sheaf on $\overline{\TT}_\VV\times\calK^\circ$ becomes a
$\calD_{\PP^{n+c}\times\CC\times\calK^\circ}$-module via Kashiwara
equivalence (using the locally closed embedding $\overline{\TT}_\VV\times\calK^\circ \cong \Gamma_F\hookrightarrow \overline{\Gamma}_F\hookrightarrow \PP^{n+c}\times\CC\times\calK^\circ$);
this is the reason for using the direct image by $\pr$
from Definition \ref{dfn-Psi}. Since it has support on the subvariety
$Z^\circ$, the corresponding perverse sheaf under the Riemann--Hilbert
correspondence is the (zeroth perverse cohomology of the) direct image
under the morphism $\Psi$ applied to the intersection complex of
$Z^\circ$.

\medskip

Finally, we want to state a mirror statement close in spirit to
Theorem \ref{theo:MirrorCompactCase} which concerns the
reduced quantum $\calD$-module. For this, we first need an extension
of the localized partial Fourier--Laplace transformation functor
$\FL^{\loc}_Y$ as defined in Formula \eqref{eq:FLloc} to a functor
acting on the category of \emph{filtered $\calD$-modules}. Without
giving the actual details (see, \emph{e.g.}  \cite[Appendix A]{SaYu}
or \cite[Definition 6.2]{ReiSe-Hodge}), let us just state that
starting from a filtered $\calD_Y$-module $(\calM,F_\bullet)$, this
version of the Fourier--Laplace transformation yields an
$\shRees$-module, where again $\shRees$ is the sheaf of Rees rings, as
discussed in Section \ref{sec:HodgeWeight} (see Formula
\eqref{eq-Rees}). We denote this $\shRees$-module by
$\FL^{\loc}_{\CC\times Y}(\calM,F_\bullet)$.

Moreover, in order to properly state the mirror theorem for nef
complete intersections, we have to take into account the so-called
\emph{mirror map}, which was not present in Theorem
\ref{theo:MirrorCompactCase} since we restricted our attention to the
Fano case there. For a sufficiently small $\varepsilon\in\RR_+$, write
$\Delta^*_\varepsilon:=\{\pointt\in (\CC^*)^r\,|\,
0<|\pointt|<\varepsilon\} \subseteq \calK^\circ$. Then the mirror map
is a morphism
\[
\textup{Mir}: \Delta^*_\varepsilon \longrightarrow H^0(X;\CC)\times U
\]
that has been defined in \cite{Giv7, CG}. Here, $U\subseteq \calK$ is
the set on which the twisted quantum product $\ast^\tw$ is defined
(converges).

With these preparations, our final mirror theorem can be stated as follows.
\begin{thm}(\cite[Conjecture 6.15]{ReiSe2}, \cite[Theorem 6.5, Theorem 6.6]{ReiSe-Hodge})
\label{thm:mirrCI}
  We have an isomorphism of $\shRees_{\CC\times \Delta^*_\varepsilon}$-modules
\begin{equation}\label{eq:FL_IC}
\FL^{\loc}_{\calK^\circ} (\calH^0 \pr_+
\IC(\underline{\CC}_{\overline{\TT}_\VV\times\calK^\circ}),
F^\Hodge_\bullet)\vert_{{\CC\times\Delta^*_\varepsilon}} \stackrel{\cong}{\longrightarrow} \left(\id_{\CC}\times \textup{Mir}\right)^*\overline{\QDM}(X,\calE).
\end{equation}
\end{thm}
This result depends in an essential way on the computation of the Hodge filtration on GKZ-systems, that is, on Theorem \ref{thm:HodgeOnGKZ}, since the expression of the Hodge filtration as the shifted order filtration on the modules $\calM_{A'}(\beta)$ for various parameters $\beta$ allows us to describe explicitly the left hand side of \eqref{eq:FL_IC}.

Notice that, by the very definition of the Dubrovin connection, the
restriction of the (reduced) quantum $\calD$-module to $\CC\times
\Delta^*_\varepsilon$ has the structure of an $\shRees_{\CC\times
  \Delta^*_\varepsilon}$-module.  A consequence of Theorem \ref{thm:mirrCI} is the
following Hodge theoretic property of the reduced quantum
$\calD$-module.
\begin{cor}\label{cor:HodgeRedQDmod}(\cite[Theorem 6.6]{ReiSe-Hodge})
Suppose $X,\calE,Y$ are as in Notation \ref{ntn-QDM+IC}.
  Then the reduced quantum $\calD$-module $\overline{\QDM}(X,\calE)$ underlies
a smooth pure polarizable twistor $\calD$-module on $\calK^\circ$ (in the sense of
\cite{Mo12}); that is, a (pure) non-commutative Hodge structure in the
sense of \cite{HS1,HS4,KKP}.
\end{cor}

\begin{exa}
  We discuss a concrete example taken from \cite[Section 1]{ReiSe2}:
  a $(2,3)$-intersection in $\PP^5$ (so, $Y\subseteq \PP^5$
  is the intersection of zero loci of generic sections of
  $\calL_1=\calO_{\PP^5}(2H)$ and $\calL_2=\calO_{\PP^5}(3H)$, where
  $H$ is the hyperplane class). The
  adjunction formula shows that this is a Fano
  variety. The (fan of the) total space of the bundle
  $\calE=\calL_1\oplus\calL_2$ has ray generators corresponding to the
  columns of the matrix
\[
B'=
\begin{pmatrix}
  1 & 0 & 0 & 0 & 0 & -1 & 0 & 0  \\
  0 & 1 & 0 & 0 & 0 & -1 & 0 & 0 \\
  0 & 0 & 1 & 0 & 0 & -1 & 0 & 0 \\
  0 & 0 & 0 & 1 & 0 & -1 & 0 & 0 \\
  0 & 0 & 0 & 0 & 1 & -1 & 0 & 0 \\
  0 & 1 & 1 & 0 & 0 &  0 & 1 & 0  \\
  0 & 0 & 0 & 1 & 1 &  1 & 0 & 1
\end{pmatrix} \in \ZZ^{7\times 8}.
\]
Then $\overline{\TT}_\VV=(\CC^*)^7$, $\calK^\circ=\CC^*$ and
the quasi-projective subvariety $Z^\circ$ of $\PP^8\times \CC \times\CC^*
=\Proj(\CC[w_0,\ldots,w_8])\times\Spec(\CC[\lambda,q^\pm])$ is given by
\[
Z^\circ =\left\{
\begin{array}{c}
w_0w_7^2w_8^3-w_1w_2w_3w_4w_5w_6=0,\\
\lambda w_0+w_1+\ldots+
w_5+qw_6+w_7+w_8=0
\end{array}
\right\}\subseteq\PP^8\times \CC \times\CC^*.
\]

The affine and the non-affine Landau--Ginzburg models of $(\PP^5,\calE)$ are given by
\begin{eqnarray*}
\psi\colon (\CC^*)^7\times \CC^*  & \longrightarrow & \CC \times \CC^* \\
(\pointt_1,\ldots, \pointt_7, \pointq) & \longmapsto & \left(-\pointt_1-\pointt_2\pointt_6 -\pointt_3\pointt_6 -\pointt_4\pointt_7-
 \pointt_5\pointt_7- \pointq\frac{\pointt_7}{\pointt_1\cdot\ldots\cdot
   \pointt_5}- \pointt_6- \pointt_7, \pointq\right)
\end{eqnarray*}
and
\begin{eqnarray*}
\Psi\colon Z^\circ & \longrightarrow & \CC \times \CC^*\\
(\pointw_0:\ldots:\pointw_8,\pointl ,\pointq) & \longmapsto & \left(\pointl, \pointq\right)
 \end{eqnarray*}

It follows from the calculations presented in \cite[Section 1]{ReiSe2}
that we have the following explicit representations of the
$\calD$-modules mentioned above: First define the operators
$P_1,P_2\in D_{\CC^*}$:
\[
\begin{array}{rcl}
P_1&=&q\cdot(3q\partial_q+1)(3q\partial_q+2)(3q\partial_q+3)(2q\partial_q+1)(2q\partial_q+2)+(q\partial_q)^6 \\ \\
&=&(q\partial_q)^2\cdot
\underbrace{\left(6q\cdot(3q\partial_q+1)(3q\partial_q+2)(2q\partial_q+1)+(q\partial_q)^4\right)}_{Q^{(2,3)}} =: (q\partial_q)^2\cdot Q^{(2,3)}\\ \\
P_2&=&q\cdot(3q\partial_q)(3q\partial_q+1)(3q\partial_q+2)(2q\partial_q)(2q\partial_q+1)+(q\partial_q)^6 \\ \\
&=& \underbrace{\left(6q\cdot(3q\partial_q+1)(3q\partial_q+2)(2q\partial_q+1)+(q\partial_q)^4\right)}_{Q^{(2,3)}}\cdot(q\partial_q)^2=:Q^{(2,3)}\cdot(q\partial_q)^2
\end{array}
\]
Then we have (we denote by $\tau$ the Fourier--Laplace dual
variable of $\lambda$, and consider the restriction to $\{\tau=1\}$ for simplicity)
\[
H^0\left(\CC^*,[\FL^{\loc}_{\calK^\circ} \calH^0\psi_+ \calO_{\overline{\TT}_\VV\times\calK^\circ}]_{|\tau=1}\right) \cong D_{\CC^*}/(P_2)
\]
and
\[
H^0\left(\CC^*,[\FL^{\loc}_{\calK^\circ} \calH^0 \pr_+ \IC(\underline{\CC}_{\overline{\TT}_\VV\times\calK^\circ})]_{|\tau=1}\right)
\cong \image(D),
\]
where $D$ is the left $\calD_{\CC^*}$-linear map
\begin{equation}\label{eq:MapExampleIntro}
\begin{array}{rcl}
  D:\CC[q^\pm]\langle\partial_q\rangle/(P_1) & \longrightarrow & \CC[q^\pm]\langle\partial_q\rangle/(P_2)\\
  Q & \longmapsto & Q\cdot (q\partial_q)^2.
\end{array}
\end{equation}
The map $D$ is well defined, its kernel is generated by $Q^{(2,3)}$
and we see that
\[
\image(D) \cong \frac{\CC[q^\pm]\langle\partial_q\rangle/(P_1) }{\ker(D)} \cong \CC[q^\pm]\langle\partial_q\rangle/(Q^{(2,3)}).
\]
The operator $Q^{(2,3)}$ is confluent, univariate and hypergeometric
(compare Subsection \ref{subsec-dim-red}) with a regular singularity
at $q=0$ and irregular singularity at $q=\infty$.

\medskip

Notice that if instead we consider a $(2,4)$-complete intersection
$Y\subset \PP^5$, then $Y$ is a Calabi--Yau manifold, and we
have
\[
H^0\left([\FL^{\loc}_{\calK^\circ} \calH^0 \pr_+
  \IC(\underline{\CC}_{\overline{\TT}_\VV\times\calK^\circ})]_{|\tau=1}\right)
\cong D_{\CC^*}/(Q^{(2,4)}),
\]
where
\[
Q^{(2,4)}=8 q \cdot (2q\partial_q+1)(4q\partial_q+1)(4q\partial_q+2)(4q\partial_q+3)-(q\partial_q)^4
\]
is a homogeneous, hence, regular (non-confluent) hypergeometric operator, with singularities
at $q=0,2^{-10}, \infty$. In this case, the Hodge theoretic result Corollary \ref{cor:HodgeRedQDmod}
simply states that $\calD_{\CC_q^*}/\calD_{\CC_q^*}\cdot Q^{(2,4)}$ underlies a pure polarized variation of Hodge structures; this is consistent with \cite[Corollary 8.1]{Si2}
and \cite[Prop. 1.13]{DeligneRigid} (see the discussion on page \pageref{page:SimpsonFedorov} above).
\schluss\end{exa}

Finally, let us remark that unlike in the previous example(s), it is
in general not easy to give a cyclic description of the intersection
cohomology $\calD$-module $\FL^{\loc}_{\calK^\circ} \calH^0 \pr_+
\IC(\underline{\CC}_{\overline{\TT}_\VV\times\calK^\circ})$. In other words, even
though
we know that it has a description as an (Fourier--Laplace transform of an)
image  of a contiguity morphism, it is not clear how to describe the
kernel of this morphism and how to give a presentation of the image as
a quotient of $\calD$ (see also \cite[Section 6]{MM17} for some
examples and conjectures).

\section*{Acknowledgements}

We would like to thank the referee for correcting various misprints
and suggesting some improvements regarding the presentation.

\section*{Table of Symbols}

Single letters (by alphabet):
\begin{itemize}
\item $A\in\ZZ^{d\times n}$, with columns $\bolda_1,\ldots,\bolda_n$
  that span $\ZZ A=\ZZ^d$ and permit a linear functional
  having positive values on them.\ref{cnv-basicA} but also
  \ref{cnv-newA} for notation in last two sections
\item $B$ a $d\times n$ submatrix of $A$ in final two sections,
  Convention \ref{cnv-newA}
\item $D_1,\ldots,D_n$ torus invariant divisors on $X$, Subsection \ref{subsec-LGmodels}
\item $j$ counts columns (and hence $x_j,\del_j,\bolda_j$), $i$ counts
  rows (hence $E_i$).
\item $\calK$ the complexified K\"ahler moduli space, the image of
  $H^2(X;\CC)$ under the exponential map, hence the quotient
  by the integer cohomology lattice scaled by $2\pi\sqrt{-1}$, Subsection \ref{subsec-GW+Dub}
\item $\overline \calK$ partial compactification of $\calK$,
    Subsection \ref{subsec-GW+Dub}
\item $[n]=\{1,2,\ldots,n\}$
\item $\mathboldq$ coordinates on $\calK$ inherited from chosen nef
  basis on $H^2(X;\CC)$,  \ref{subsec-GW+Dub}
\item $r=\dim_\CC H^2(X;\CC)$
\item $\TT$ the $d$-torus, Subsection \ref{subsec-torusaction}, but
  see Convention \ref{cnv-newA} and \eqref{eqn-h-diagram} for the
  final sections
\item $\overline\TT$ the quotient torus modulo $0$-th component of
  $\TT$ in final two sections
\item $U$ complement of $Z$
\item $\VV$ total space of tautological bundle $\calO_{\PP^n}(-1)$
\item $X$ smooth projective toric variety to fan $\Sigma_X$, often but
  not always Fano, Subsection \ref{subsec-GW+Dub},
\item $Y$ complete intersection in $X$ of codimension $c$,
\item $Z$ tautological hypersurface in $\PP^n\times\CC^{n+1}$
\item $Z^\circ$ the closure in $\PP^{n+c}\times \CC \times
  \calK^\circ$ of the graph of the function defined in
  \eqref{eqn-affLG}, see \eqref{eqn-Zcirc}
\item[]
\end{itemize}

Compounds (by alphabet of first occurring letter):]
\begin{itemize}
\item $\mathfrak{A}_A$ the admissible parameters, Definition
  \ref{dfn-admissible}
\item $\conv(S)$ the convex hull of $S$, before Definition \ref{96}
\item $c_{\beta,\beta'}\colon M_A(\beta)\to M_A(\beta')$ contiguity
  operators, Subsection \ref{subsec:Contiguity}
\item $\Div_T(X)$ equivariant divisor group of toric variety $X$,
  isomorphic to actual divisor group, generated by rays of fan
  $\Sigma_X$, \eqref{eq-exseqtoric}
\item $E_i$ Euler operators, Definition \ref{def:GKZ}
\item $F^\Hodge$ the Hodge filtration on the mixed Hodge module
  $\calM$, Subsections \ref{subsec-setupHodge}, \ref{sec:HodgeWeight},
  \eqref{eq:Fedorov}, \eqref{eqn:SaYu}
\item $\FL(\calM)$ the Fourier--Laplace transform, \eqref{eqn-exp}
\item $F^\ord$ the order filtration on rings of differential
  operators
\item $G^\psi$, $\calG^\psi$, $G^\psi_0$ Fourier-Laplace transformed
  Brieskorn lattice and variations, \eqref{eqn-FLBrieskorn} and
  following page
\item $h_A\colon \TT\to\CC^n$ the monomial map induced by $A$,
  Subsection \ref{subsec-FGKZ}
\item $(H^\sideA,\nabla^\sideA)$ small Dubrovin connection,
  \eqref{eq:DubrovinConn}
\item $H_{A,i}(N;\beta)$ the $i$-th Euler--Koszul homology of the
  toric module $N$ for the parameter $\beta$
\item $H_A(\beta)$ the hypergeometric ideal,
  \ref{def:GKZ}
\item $\widehat\calM$ the Fourier--Laplace transform of the module
  $\calM$
\item $M_A(\beta)$ the hypergeometric module,
  \ref{def:GKZ}
\item $\qdeg_A(N)$ the quasi-degrees  of
  an $A$-graded module, Definition \ref{19}
\item $\Rees$, $\shRees$ the twisted Rees ring/sheaf of differential
  operators on various spaces, Definition \ref{eq-Rees}, Proposition \ref{prp-locBrieskorn}
\item $\Radon(\calM)$ the Radon transform, Proposition \ref{prop-DE}
\item $S_A$ the semigroup ring $\CC[\NN A]$, Subsection \ref{subsec-torusaction}
\item $S^L_A$ the $L$-graded ring of $S_A$, Theorem \ref{13}
\item $\sRes(A)$ the strongly resonant parameters for $A$,
  Definition \ref{def-sRes}
\item $\tdeg_A(N)$ the true degrees  of
  an $A$-graded module, Definition \ref{19}
\item $\TT_\VV$ the $(n+c)$-torus, Definition \ref{dfn-affineLG}
\item $(W,\boldc)$ Landau--Ginzburg model on $\calK$, Definition
  \ref{dfn-LGmodelK}, \eqref{eqn-psi}
\item $W_k\calM$ the weight filtration on the mixed Hodge module
  $\calM$, Subsections \ref{subsec-setupHodge} and \ref{subsec-weight}
\item $X_A$ affine toric variety and spectrum of $S_A$, closure of
  $\TT$-orbit through $(1,\ldots,1)$, Subsection \ref{subsec-FGKZ}
\item[]
\end{itemize}

Greek letters and other symbols:
\begin{itemize}
\item $\Box_\boldu=\bolddel^{\boldu_+}-\bolddel^{\boldu_-}$ for $\boldu\in\ker A$,.
\item $\ast^{\tw}$ twisted quantum product,  \eqref{eqn-*tw}
\item $\Delta^L_A$ the $(A,L)$-polyhedron, the convex hull of the
  origin and all $\bolda_j^L$, $\Delta_A$ special case to
  $L=\boldzero$, Definition \ref{96} and Subsection \ref{sec:HodgeWeight}
\item $\Delta^*_\eps$ small ball around origin in $\calK^\circ$
\item $\Sigma^L_A$ initial complex of ideal for generic weight
  $L$, Definition \ref{dfn-Sigma}
\item $\Sigma_X$ fan of $X$
\item $\varphi\colon \overline\TT\to \CC^n$ family of Laurent
  polynomials, Theorem \ref{theo:4TermSeq}
\item $\Phi^L_A$ the $(A,L)$-umbrella, Definition \ref{96}
\item $\psi$ affine Landau--Ginzburg model on $\TT_\VV\times\calK$,
  Definition \ref{dfn-affineLG}
\item $\Psi$ non-affine Landau--Ginzburg model on $\CC\times \calK^\circ$,
  Definition \ref{dfn-Psi}

\end{itemize}

\bibliographystyle{amsalpha}

\def\cprime{$'$}
\providecommand{\bysame}{\leavevmode\hbox to3em{\hrulefill}\thinspace}
\providecommand{\MR}{\relax\ifhmode\unskip\space\fi MR }
\providecommand{\MRhref}[2]{  \href{http://www.ams.org/mathscinet-getitem?mr=#1}{#2}
}
\providecommand{\href}[2]{#2}

\end{document}